# HOMOGENEOUS FUNCTIONS ON LIGHT CONES: THE INFINITESIMAL STRUCTURE OF SOME DEGENERATE PRINCIPAL SERIES REPRESENTATIONS

ROGER E. HOWE AND ENG-CHYE TAN

## 1. INTRODUCTION

One of the classic and most pleasing calculations in the representation theory of semisimple Lie groups is Bargmann's [Ba] (see also [Na]) description of the structure of representations of $\mathrm{SL}(2, \mathbb{R})$. Bargmann constructs a basis of eigenvectors for the compact subgroup $\mathrm{SO}(2) \subseteq \mathrm{SL}(2, \mathbb{R})$ and then explicitly describes the action of the Lie algebra of $\mathrm{SL}(2, \mathbb{R})$ on this basis. From the straightforward and explicit formulas that result, it is easy to make conclusions about unitarity, irreducibility, isomorphism, and other natural questions concerning these representations.

Unfortunately Bargmann's example stands almost alone; the literature tends to give the impression that once beyond the civilized confines of $\mathrm{SL}(2, \mathbb{R})$, one immediately finds oneself in a dense thicket of difficulties, where explicit calculation is difficult and of little value. The purpose of this paper is to point to an infinite class of examples in which the pleasant characteristics of Bargmann's work can be retained; one can describe a model for representations in which the action of the Lie algebra is given by explicit simple formulas, which permit direct determination of important structural properties of the representations.

The examples we will treat are spaces of functions that live on light cones and satisfy appropriate homogeneity conditions. These spaces are modules for an appropriate isometry group—$\mathrm{O}(p, q)$, $\mathrm{U}(p, q)$, or $\mathrm{Sp}(p, q)$. In case $q = 1$, the spherical principal series are among our examples, and Kostant's celebrated results [Ko] on the complementary series can be read off from our formulas, excluding the one exceptional rank one group. (Essentially this calculation was done, in less systematic fashion, in [Jo, JW].)

We are hardly the first to study these representations. Significant parts of our results, including Dixmier's [Di] study of the De Sitter group, Hirai's notes on $\mathrm{O}(p, 1)$ and $\mathrm{U}(p, 1)$ [Hi1, Hi2], Kostant's work [Ko], the calculations of Johnson and Wallach [JW] for the spherical principal series, the studies of Schlichtkrull [Sc], Sekiguchi [Se], Strichartz [St1, St2], Rallis and Schiffmann [RS], Rossmann [Ro], Molcanov [Mo], Faraut [Fa] on functions on hyperboloids, the calculations of Klimyk and various colleagues (see [KG] and the references therein) and others [Br, Cl, Cw, FR, GN, HC1, JV,

Received by the editors January 27, 1992.
1991 *Mathematics Subject Classification.* Primary 22E46; Secondary 17B10.
*Key words and phrases.* Complementary series, composition series, degenerate principal series representations, $K$-type diagrams, light cones, unitary representations.
The first author was partially supported by NSF Grant DMS #9103608.

Ni, Sa, Sp1, Tk, Vi] are already in the literature. The paper [Br] has diagrams that resemble some of ours. The unitary duals of the rank one groups have been classified (see [BS, Cl, Hi1, Hi2] and the references therein). In fact, most or all of our results may be in the literature, at least implicitly. The main purpose of the present paper is to make our results explicit and to emphasize how concrete and how simple a picture one may draw of these representations. We use virtually no technology. There is no clear reason why these calculations could not have been done in the 1940s or early 1950s. If they had been, the development of semisimple representation theory might have been substantially different. For instance, besides the results of Kostant mentioned above, it is easy to read off from our results the unitary duals of $O(2, 1)$, $O(4, 1)$, and $U(2, 1)$. The last of these was not known until the mid 1960s [Hi2].

Here is a summary of the paper. In §2 we set the context and describe our representations in the case of the indefinite orthogonal groups $O(p, q)$. We describe the $K$-spectrum of the representations, where $K = O(p) \times O(q)$ is the maximal compact subgroup of $O(p, q)$, and we compute the infinitesimal action of the Lie algebra $\mathfrak{o}(p, q)$ in terms of the $K$-spectrum. On the basis of the formulas for the $\mathfrak{o}(p, q)$ action, we determine composition series for our representations. It turns out that the structure of the composition series involves some integrality conditions which make it depend in a significant way on the parities of $p$ and $q$. Also, we must distinguish the cases when $p$ or $q$ is less than or equal to 2, since these involve some degeneration in the $K$-spectrum that affects the representation structure.

In §3 we apply our formulas to determine which of our representations carry invariant Hermitian forms, and on which subquotients these forms can be definite, that is, which constituents of our representations are unitary. In particular, we determine the extent of the complementary series and what happens at the end of the complementary series.

In §§4 and 5 we adapt the discussion of §§2, 3 to the unitary groups $U(p, q)$ and the indefinite symplectic groups $Sp(p, q)$. The overall features of the picture for these groups is controlled by $O(p, q)$, as is made clear by Lemmas 4.1 and 5.3, but differences in the $K$-spectrum in the degenerate cases when $q = 1$ lead to interesting new phenomena, notably considerable variation in the length of the complementary series, including Kostant's results concerning the spherical complementary series of $U(p, 1)$ and $Sp(p, 1)$.

In order not to distract from the coherence and directness of the phenomena that we present, we omit the details of various calculations in §§2–5. These are worked out in §6. The final §7 discusses how our results fit into the general theory of representations of semisimple Lie groups.

To conclude this introduction, we note that those readers who are familiar with the theory of reductive dual pairs [Ho3] will find our results can be profitably interpreted in terms of this theory, which can be used to predict much of what we present. We emphasize again, however, that our elementary and computational approach is independent of this or any general theory (except that we use [Ho3] as a convenient reference for facts about spherical harmonics that must somewhere be in the classical literature).

## 2. The $\mathrm{O}(p,q)$-module structure of homogeneous functions on light cones: $K$-types and composition series

Let $\mathbb{R}^{p+q} \simeq \mathbb{R}^p \oplus \mathbb{R}^q$ denote the vector space of tuples $w = (x,y)$ where

$$x = \begin{bmatrix} x_1 \\ x_2 \\ \vdots \\ x_p \end{bmatrix} \in \mathbb{R}^p \quad \text{and} \quad y = \begin{bmatrix} y_1 \\ y_2 \\ \vdots \\ y_q \end{bmatrix} \in \mathbb{R}^q.$$

To fix ideas, we will take $p \geq q$. We consider $\mathbb{R}^{p+q}$ to be endowed with the indefinite inner product

$$(w, w')_{p,q} = (x, x')_p - (y, y')_q = \sum_{j=1}^{p} x_j x_j' - \sum_{j=1}^{q} y_j y_j'.$$

We will use the notation

$$r_{p,q}^2(w) = (w, w)_{p,q}$$

and likewise for $r_p^2(x)$, $r_q^2(y)$.

Let $\mathrm{O}(p,q)$ denote the isometry group of $(\cdot, \cdot)_{p,q}$; that is, $\mathrm{O}(p,q)$ is the subgroup of linear transformations $g \in \mathrm{GL}(p+q, \mathbb{R})$ such that

$$(gw, \, gw')_{p,q} = (w, \, w')_{p,q}, \qquad w, \, w' \in \mathbb{R}^{p+q}.$$

Clearly the group $\mathrm{O}(p)$, the isometry group of $(\cdot, \cdot)_p$ on $\mathbb{R}^p$, can be regarded as a subgroup of $\mathrm{O}(p,q)$ by letting it act on the $x$ component of $w = (x,y)$. Similarly, $\mathrm{O}(q)$, the isometry group of $(\cdot, \cdot)_q$ on $\mathbb{R}^q$, can be considered a subgroup of $\mathrm{O}(p,q)$. Combining these, we get an embedding of $\mathrm{O}(p) \times \mathrm{O}(q) \subseteq \mathrm{O}(p,q)$, as the subgroup of $\mathrm{O}(p,q)$, which stabilizes the decomposition $\mathbb{R}^{p+q} \simeq \mathbb{R}^p \oplus \mathbb{R}^q$ implicit in the notation $w = (x,y)$. The group $\mathrm{O}(p) \times \mathrm{O}(q)$ is clearly a compact subgroup of $\mathrm{O}(p,q)$. It is well known and not hard to see that $\mathrm{O}(p) \times \mathrm{O}(q)$ is a maximal compact subgroup of $\mathrm{O}(p,q)$.

Let $M_{m,n}(F)$ denote the space of $m \times n$ matrices over a field $F$; if $m = n$, we shall abbreviate $M_{m,m}(F)$ by $M_m(F)$. We denote by $\mathfrak{o}(p,q)$ the Lie algebra of $\mathrm{O}(p,q)$, the space of $(p+q) \times (p+q)$ matrices $T \in M_{p+q}(\mathbb{R})$, which are generators of one-parameter subgroups $\{e^{sT} \mid s \in \mathbb{R}\}$ of $\mathrm{O}(p,q)$. By differentiating the relation

$$(e^{sT}w, e^{sT}w')_{p,q} = (w, w')_{p,q},$$

we obtain the description

(2.1)

$$\mathfrak{o}(p,q) = \big\{ T \in M_{p+q}(\mathbb{R}) \, \big| \, (Tw, w')_{p,q} + (w, Tw')_{p,q} = 0$$

$$\text{for all } w, \, w' \in \mathbb{R}^{p+q} \big\}$$

$$= \left\{ \begin{bmatrix} A & B \\ B^{\mathrm{t}} & C \end{bmatrix} \, \middle| \, A = -A^{\mathrm{t}} \in M_p(\mathbb{R}); \; C = -C^{\mathrm{t}} \in M_q(\mathbb{R}); \; B \in M_{p,q}(\mathbb{R}) \right\}.$$

Here $A^{\mathrm{t}}$ denotes the transpose of the matrix $A$.

Consider the light cone

$$(2.2) \qquad \begin{aligned} X_{p,q}^0 = X^0 &= \left\{ w \in \mathbb{R}^{p+q} - \{0\} \,\Big|\, r_{p,q}^2(w) = 0 \right\} \\ &= \left\{ (x,y) \in \mathbb{R}^{p+q} - \{0\} \,\Big|\, r_p^2(x) = r_q^2(y) \right\}. \end{aligned}$$

Observe that we exclude the origin from $X^0$.

The action of $O(p,q)$ on $\mathbb{R}^{p+q}$ preserves the light cone $X^0$. It is well known (in particular, it is a direct consequence of Witt's Theorem [Ja]) that $X^0$ consists of a single orbit for $O(p,q)$, that is, the action of $O(p,q)$ on $X^0$ is transitive. Thus if $w_1$ is a specified point in $X^0$ and $Q_1 \subseteq O(p,q)$ is the subgroup that stabilizes $w_1$, then the map $g \to gw_1$ defines a bijection (in fact, a diffeomorphism)

$$\alpha : O(p,q)/Q_1 \simeq X^0.$$

Let $P_1 \subseteq O(p,q)$ be the stabilizer of the line $\mathbb{R}w_1$ through $w_1$. Then $P_1$ is a parabolic subgroup [Kn] of $O(p,q)$, and $\alpha(P_1) = \mathbb{R}^\times w_1$; more precisely, the map $p \to \tilde{\alpha}(p)$ where $pw_1 = \tilde{\alpha}(p)w_1$ defines a group isomorphism

$$(2.3) \qquad \tilde{\alpha} : P_1/Q_1 \simeq \mathbb{R}^\times.$$

The second description of $X^0$ in (2.2) allows us to analyze the action of $O(p) \times O(q)$ on $X^0$. Let

$$\mathbb{S}^{p-1} = \{x \in \mathbb{R}^p \mid (x,x)_p = 1\}$$

be the unit sphere in $\mathbb{R}^p$. It is a homogeneous space for $O(p)$. We have on $X^0$ an analog of "polar coordinates".

**Scholium 2.1.** *The mapping*

$$(2.4) \qquad \beta : X^0 \to \mathbb{S}^{p-1} \times \mathbb{S}^{q-1} \times \mathbb{R}_+^\times$$

*defined by*

$$w = (x,y) \to \left( \frac{x}{r_p^2(x)^{1/2}}, \frac{y}{r_q^2(y)^{1/2}}, r_p^2(x)^{1/2} \right)$$

*is a diffeomorphism. The level sets of $r_p^2(x)$ are $O(p) \times O(q)$ orbits, equivalent via a scalar dilation to $\mathbb{S}^{p-1} \times \mathbb{S}^{q-1}$.*

The action of $O(p,q)$ on $\mathbb{R}^{p+q}$ can be used to define an action of $O(p,q)$ on functions on $\mathbb{R}^{p+q}$ by the familiar formula

$$(2.5) \qquad \rho(g)(f)(w) = f(g^{-1}w), \qquad g \in O(p,q), \ w \in \mathbb{R}^{p+q},$$

for a function $f$ on $\mathbb{R}^{p+q}$. Precisely the same formula defines an action on functions on any set $X \subseteq \mathbb{R}^{p+q}$ that is invariant under $O(p,q)$. Furthermore, the restriction mapping $r_X : f \to f_1$ taking a function $f$ to its restriction to $X$, clearly commutes with the $O(p,q)$ action. This means that if we have a function $f$ on $X$, we can compute $\rho(g)(f)$ by extending $f$ to some function $\tilde{f}$ on all of $\mathbb{R}^{p+q}$, computing $\rho(g)(\tilde{f})$, then restricting back to $X$. This compatibility justifies at least partially our abuse of notation in failing to specify on what set $f$ is defined in formula (2.5). We will be interested in taking $X = X^0$, the light cone.

If we differentiate the action (2.5) along one-parameter subgroups of $\mathrm{O}(p, q)$, we obtain a representation of the Lie algebra $\mathfrak{o}(p, q)$ on functions via differential operators. Straightforward computation shows that the space of operators representing $\mathfrak{o}(p, q)$ is spanned by the first-order operators:

$$
\begin{aligned}
& k_{jk}^1 = x_j \frac{\partial}{\partial x_k} - x_k \frac{\partial}{\partial x_j}, && 1 \le j < k \le p, \\
(2.6) \qquad & k_{lm}^2 = y_l \frac{\partial}{\partial y_m} - y_m \frac{\partial}{\partial y_l}, && 1 \le l < m \le q, \\
& p_{jl} = x_j \frac{\partial}{\partial y_l} + y_l \frac{\partial}{\partial x_j}, && 1 \le j \le p, \ 1 \le l \le q.
\end{aligned}
$$

Here the $k_{jk}^1$ describe the action of the Lie subalgebra $\mathfrak{o}(p) \subseteq \mathfrak{o}(p, q)$ corresponding to the subgroup $\mathrm{O}(p)$ of $\mathrm{O}(p, q)$. Similarly the $k_{lm}^2$ describe the action of $\mathfrak{o}(q)$. The $p_{jl}$ span a complement to $\mathfrak{o}(p) \oplus \mathfrak{o}(q)$ in $\mathfrak{o}(p, q)$, corresponding to the off-diagonal matrix $B$ in formula (2.1). We denote the span of the $p_{jl}$ by $\mathfrak{p}$.

We propose to study spaces of homogeneous functions on the light cone $X^0$. Thus for $a \in \mathbb{C}$ denote by $S^a(X^0)$ the space

$$
(2.7) \qquad S^a(X^0) = \left\{ f \in C^\infty(X^0) \,\middle|\, f(tw) = t^a f(w), \ w \in X^0, \ t \in \mathbb{R}_+^\times \right\}
$$

of smooth functions on $X^0$ homogeneous of degree $a$ under dilations by positive scalars. Since $\mathrm{O}(p, q)$ commutes with scalar dilations, it is clear that $S^a(X^0)$ will be invariant under the action (2.5) of $\mathrm{O}(p, q)$. Indeed there is no need to consider only positive dilations. We could further define two subspaces

$$
S^{a\pm}(X^0) = \left\{ f \in S^a(X^0) \,\middle|\, f(-w) = \pm f(w) \right\}.
$$

It is clear that

$$
(2.8) \qquad S^a(X^0) = S^{a+}(X^0) \oplus S^{a-}(X^0)
$$

and that each of $S^{a+}(X^0)$ and $S^{a-}(X^0)$ are individually invariant under $\mathrm{O}(p, q)$. Thus it would be possible, and might seem advisable, to study the $S^{a+}(X^0)$ and the $S^{a-}(X^0)$ separately. However, we will find it convenient and enlightening to lump $S^{a\pm}(X^0)$ together in the space $S^a(X^0)$.

Before proceeding to explicit computations, we observe that the spaces $S^{a\pm}(X^0)$ are naturally identifiable as $\mathrm{O}(p, q)$ modules to certain induced representations, of a type known informally as "degenerate principal series". Recall the homomorphism $\tilde{\alpha} : P_1 \to \mathbb{R}^\times$ of formula (2.3). For $a \in \mathbb{C}$, and a choice $\pm$ of sign, define a quasicharacter ($=$ homomorphism to $\mathbb{C}^\times$) of $P_1$ by

$$
(2.9) \qquad \psi_a^\pm(p) = \begin{cases} \tilde{\alpha}(p)^a & \text{if} \quad \tilde{\alpha}(p) > 0, \\ \pm\,|\tilde{\alpha}(p)|^a & \text{if} \quad \tilde{\alpha}(p) < 0. \end{cases}
$$

Define a map

$$
(2.10) \qquad \delta : S^{a\pm}(X^0) \to C^\infty(G)
$$

by the formula

$$
\delta(f)(g) = \rho(g)(f)(w_1) = f(g^{-1}w_1).
$$

Here $w_1 \in X^0$ is the point used to define $P_1$.

Straightforward formal checking shows that

$$(2.11) \qquad \begin{aligned} \delta(f)(pg) &= \psi_a^{\pm}(p)^{-1}\delta(f)(g)\,, \qquad p \in P_1,\ g \in G, \\ \delta(f)(gh) &= \delta(\rho(h)f)(g)\,. \end{aligned}$$

These formulas, plus the definition of $P_1$, imply that $\delta$ defines an equivalence of $\mathrm{O}(p,q)$ representations between $S^{a\pm}(X^0)$ and the representation $\mathrm{ind}_{P_1}^G(\psi_a^{\pm})^{-1}$ induced from the character $(\psi_a^{\pm})^{-1}$ of $P_1$ [Kn]. (Here we are using unnormalized induction.) If $q = 1$ and we choose $+$ from $\pm$, these representations are the spherical principal series of $\mathrm{O}(p,1)$ [Kn].

Now let us investigate the structure of the $S^a(X^0)$ as $\mathrm{O}(p,q)$ modules. First consider the action of the compact subgroup $\mathrm{O}(p) \times \mathrm{O}(q)$. The decomposition (2.4) of $X^0$ shows that any function in $S^a(X^0)$ is determined by its restriction to $\beta^{-1}(\mathbb{S}^{p-1} \times \mathbb{S}^{q-1} \times \{1\})$. More exactly, restriction to $\beta^{-1}(\mathbb{S}^{p-1} \times \mathbb{S}^{q-1} \times \{1\})$ of elements of $S^a(X^0)$ defines an isomorphism of $\mathrm{O}(p) \times \mathrm{O}(q)$ modules.

The classical theory of spherical harmonics [Ho3, Vi, Zh] describes $C^{\infty}(\mathbb{S}^{p-1})$ as a representation of $\mathrm{O}(p)$. It is the (topological) direct sum of the finite-dimensional spaces $\mathcal{H}^m(\mathbb{R}^p)$ consisting of the (restrictions to $\mathbb{S}^{p-1}$ of the) harmonic polynomials of degree $m$, $m \in \mathbb{Z}_+$;

$$(2.12) \qquad C^{\infty}(\mathbb{S}^{p-1}) \simeq \sum_{m \geq 0} \mathcal{H}^m(\mathbb{R}^p).$$

The space $\mathcal{H}^m(\mathbb{R}^p)$ consists of polynomials of degree $m$ which are annihilated by the Laplace operator

$$\Delta_p = \sum_{j=1}^{p} \frac{\partial^2}{\partial x_j^2}.$$

Let $\mathcal{P}^m(\mathbb{R}^p)$ denote the space of polynomials of degree $m$ on $\mathbb{R}^p$. Then the Laplace operator maps $\mathcal{P}^m(\mathbb{R}^p)$ to $\mathcal{P}^{m-2}(\mathbb{R}^p)$; that is,

$$(2.13) \qquad \Delta_p:\ \mathcal{P}^m(\mathbb{R}^p) \to \mathcal{P}^{m-2}(\mathbb{R}^p),$$

with kernel $\mathcal{H}^m(\mathbb{R}^p)$. It is part of the theory of spherical harmonics that the mapping (2.13) is surjective; hence

$$\begin{aligned} \dim \mathcal{H}^m(\mathbb{R}^p) &= \dim \mathcal{P}^m(\mathbb{R}^p) - \dim \mathcal{P}^{m-2}(\mathbb{R}^p) \\ &= \binom{m+p-1}{p-1} - \binom{m+p-3}{p-1}. \end{aligned}$$

If $p > 1$, this is always nonzero, but if $p = 1$, it is zero for $m \geq 2$. Thus although the sum (2.12) always is nominally over all of $\mathbb{Z}_+$, for $p = 1$, it is effectively only over $m = 0, 1$. Each $\mathcal{H}^m(\mathbb{R}^p)$ defines an irreducible representation of $\mathrm{O}(p)$, and the representations defined by two different $\mathcal{H}^m(\mathbb{R}^p)$ are inequivalent.

Combining the decompositions (2.12) for $p$ and for $q$, we obtain a decomposition

$$(2.14) \qquad C^{\infty}(\mathbb{S}^{p-1} \times \mathbb{S}^{q-1}) \simeq \sum_{m,n \geq 0} \mathcal{H}^m(\mathbb{R}^p) \otimes \mathcal{H}^n(\mathbb{R}^q)$$

of the functions on $\mathbb{S}^{p-1} \times \mathbb{S}^{q-1}$ into irreducible inequivalent subspaces for $\mathrm{O}(p) \times \mathrm{O}(q)$. The above remarks about restriction to $\beta^{-1}(\mathbb{S}^{p-1} \times \mathbb{S}^{q-1} \times \{1\})$ show that equation (2.14) also yields a description of the $\mathrm{O}(p) \times \mathrm{O}(q)$-module structure of $S^a(X^0)$. The main point to observe is that we can adjust the homogeneity of functions by multiplication by powers of $r_p^2$ (or $r_q^2$) without changing the way they transform under $\mathrm{O}(p) \times \mathrm{O}(q)$. Thus we can define an embedding

$$(2.15) \qquad j_{a,m,n} = j_a : \mathcal{H}^m(\mathbb{R}^p) \otimes \mathcal{H}^n(\mathbb{R}^q) \to S^a(X^0)$$

by the formula

$$j_a(h_1 \otimes h_2)(x,y) = h_1(x)h_2(y)r_p^2(x)^b,$$

where $h_1 \in \mathcal{H}^m(\mathbb{R}^p)$, $h_2 \in \mathcal{H}^n(\mathbb{R}^q)$, and $b$ is chosen so that the total degree of homogeneity is $a$, that is, $m + n + 2b = a$. In formula (2.15) we are implicitly restricting $h_1 h_2 r_p^{2b}$ to $X^0$. The following statement summarizes our discussion of the $\mathrm{O}(p) \times \mathrm{O}(q)$ action.

**Lemma 2.2.** *As an $\mathrm{O}(p) \times \mathrm{O}(q)$ module, the space $S^a(X^0)$ decomposes into a direct sum*

$$(2.16) \qquad S^a(X^0) \simeq \sum_{m,n \geq 0} j_a(\mathcal{H}^m(\mathbb{R}^p) \otimes \mathcal{H}^n(\mathbb{R}^q))$$

*of mutually inequivalent irreducible representations.*

We will refer to the spaces $j_a(\mathcal{H}^m(\mathbb{R}^p) \otimes \mathcal{H}^n(\mathbb{R}^q))$ as the **$K$-types** of $S^a(X^0)$. We remind the reader that in case $p$ (or $q$) equals 1, the sum over $m$ (or $n$) is limited to $m$ (or $n$) equal to 0 or 1.

Lemma 2.2 tells us the structure of $S^a(X^0)$ as an $\mathrm{O}(p) \times \mathrm{O}(q)$ module. Since the $K$-types all define distinct $\mathrm{O}(p) \times \mathrm{O}(q)$ subrepresentations, any $\mathrm{O}(p) \times \mathrm{O}(q)$-invariant subspace of $S^a(X^0)$ will be a sum of the $K$-types it contains. This will in particular be true of any $\mathrm{O}(p,q)$-invariant subspace of $S^a(X^0)$. Thus, to understand the $\mathrm{O}(p,q)$-module structure of $S^a(X^0)$, we need to see how $\mathrm{O}(p,q)$ transforms one $K$-type to another. For this purpose, it is helpful to represent the $K$-types graphically, as points in the plane $\mathbb{R}^2$. Precisely, we will represent $j_a(\mathcal{H}^m(\mathbb{R}^p) \otimes \mathcal{H}^n(\mathbb{R}^q))$ as the point with integer coordinate $(m,n)$ in the positive quadrant $(\mathbb{R}_+)^2$ in $\mathbb{R}^2$. Thus the set of all $K$-types is matched to the set of all integer points in the positive quadrant, as illustrated by Diagram 2.17 (see next page). (This picture is for $p, q \geq 2$, which we will assume for the present. We will treat the degenerate case $q = 1$ later.) We have used two different symbols to mark the points to remind the reader of the decomposition of $S^a(X^0)$ into $S^{a\pm}(X^0)$. The dots indicate $K$-types in $S^{a+}(X^0)$ (the parity of $m + n$ is even), and the crosses indicate $K$-types in $S^{a-}(X^0)$ (the parity of $m+n$ is odd). Sometimes we will indicate a $K$-type $j_a(\mathcal{H}^m(\mathbb{R}^p) \otimes \mathcal{H}^n(\mathbb{R}^q))$ simply by its associated coordinates $(m, n)$.

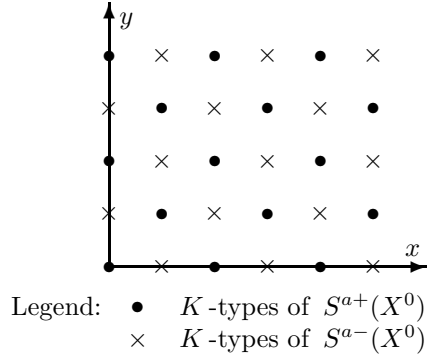

Legend:    •    $K$-types of $S^{a+}(X^0)$
           ×    $K$-types of $S^{a-}(X^0)$

DIAGRAM 2.17

To go beyond the $K$-type structure of $S^a(X^0)$, we will investigate the effect of applying elements from the subspace $\mathfrak{p}$ of the Lie algebra $\mathfrak{o}(p,q)$ to the $K$-types. Since the $K$-types are individually invariant under $\mathfrak{o}(p) \oplus \mathfrak{o}(q)$, any collection of $K$-types such that any one of them is mapped by $\mathfrak{p}$ into the span of the others will define (i.e., its span will be) an $\mathfrak{o}(p,q)$ submodule of $S^a(X^0)$. An early theorem of Harish-Chandra [HC2] guarantees that the closure of this span in $S^a(X^0)$ will be an $\mathrm{O}(p,q)$ submodule. Thus our computation of the action of $\mathfrak{p}$ on individual $K$-types is in principle (and will turn out to be in practice) sufficient for understanding the submodule structure of $S^a(X^0)$.

The operators of $\mathfrak{p}$ are described in the formula (2.6). A computation, which will be carried out in §6.1, establishes the following rules.

**Lemma 2.3.** *For each pair $(m,n)$, there are maps*

$$T_{m,n}^{\pm,\pm} : \ \mathfrak{p} \otimes (\mathcal{H}^m(\mathbb{R}^p) \otimes \mathcal{H}^n(\mathbb{R}^q)) \to \mathcal{H}^{m\pm1}(\mathbb{R}^p) \otimes \mathcal{H}^{n\pm1}(\mathbb{R}^q),$$

*which are independent of $a$ and which are nonzero as long as the target space is nonzero, such that the action of $z \in \mathfrak{p}$ on the $K$-type $j_a(\mathcal{H}^m(\mathbb{R}^p) \otimes \mathcal{H}^n(\mathbb{R}^q))$ is described by the formula*

$$
\begin{aligned}
(2.18) \qquad \rho(z)j_a(\phi) = & (a-m-n)j_a(T_{mn}^{++}(z \otimes \phi)) \\
& + (a-m+n+q-2)j_a(T_{mn}^{+-}(z \otimes \phi)) \\
& + (a+m-n+p-2)j_a(T_{mn}^{-+}(z \otimes \phi)) \\
& + (a+m+n+p+q-4)j_a(T_{mn}^{--}(z \otimes \phi)),
\end{aligned}
$$

*where $\phi \in \mathcal{H}^m(\mathbb{R}^p) \otimes \mathcal{H}^n(\mathbb{R}^q)$.*

*Proof.* See §6.1.  □

In terms of the graphical representation (see Diagram 2.17) of the $K$-types, we may picture this formula as follows. Starting from a given point $(m,n)$, the action of $\mathfrak{p}$ potentially can take us to any of the diagonally adjacent points $(m \pm 1, \ n \pm 1)$. Whether we can, in a given $S^a(X^0)$, actually proceed from $(m,n)$ to a given one of the points $(m \pm 1, n \ \pm 1)$ depends on whether the coefficient of the operator $T_{mn}^{\pm\pm}$ (e.g., $a - m - n$ for $T_{mn}^{++}$) in the formula (2.18) is nonzero or not. We can represent the situation graphically as in Diagram 2.19.

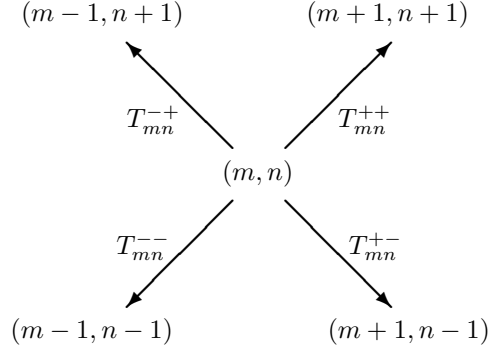

DIAGRAM 2.19

The four corner points in Diagram 2.19 are the $K$-types we may hope to reach from $(m, n)$ by the action of $\mathfrak{p}$; whether we actually can get to a given one of the points depends on whether the *transition coefficient* $A_{p,q,a}^{\pm\pm}(m, n)$ given by

$$
\begin{aligned}
(2.20) \qquad A_{p,q,a}^{++}(m, n) &= a - m - n, \\
A_{p,q,a}^{+-}(m, n) &= a - m + n + q - 2, \\
A_{p,q,a}^{-+}(m, n) &= a + m - n + p - 2, \\
A_{p,q,a}^{--}(m, n) &= a + m + n + p + q - 4,
\end{aligned}
$$

is nonzero. We shall omit the subscripts in $A_{p,q,a}^{\pm\pm}$ if it is clear from context.

*Remark.* We see that from a given $(m, n)$ by arbitrarily repeated action of $\mathfrak{p}$ we can never hope to get to all other points $(m', n')$—only to those such that $m' + n'$ has the same parity as $m + n$. This reflects the fact, already observed, that $S^a(X^0)$ breaks up into the two invariant subspaces $S^{a\pm}(X^0)$.

We see that if $a$ is not an integer, then the transition coefficients $A^{\pm\pm}$ are never zero for any $(m, n)$ in $(\mathbb{Z}_+)^2$. From this and Diagram 2.19, we see that, starting from any $K$-type $j_a(\mathcal{H}^m(\mathbb{R}^p) \otimes \mathcal{H}^n(\mathbb{R}^q))$, we may by sufficiently many applications of elements from $\mathfrak{p}$ eventually reach any other $K$-type corresponding to $(m', n')$, providing $m' + n'$ has the same parity as $(m, n)$. Combining this with the descriptions (2.17), we can make our first conclusion concerning reducibility of the $S^a(X^0)$.

**Proposition 2.4.** *If $a$ is not an integer, then each of $S^{a\pm}(X^0)$ is an irreducible* $\mathrm{O}(p, q)$ *module.*

Consider now what happens when $a$ is an integer. Suppose $a \geq 0$, and consider a $K$-type $(m, n)$ such that $m + n = a$. This means $A^{++} = 0$, so we cannot move from $(m, n)$ to $(m + 1, n + 1)$ by acting by $\mathfrak{p}$. We can perhaps move to $(m', n') = (m + 1, n - 1)$ or to $(m', n') = (m - 1, n + 1)$, but from these points we again cannot move to the upper right, as seen in terms of Diagram 2.19 since again $m' + n' = a$. We can see that the action of $\mathfrak{p}$, no matter how many times repeated, will never take us above the line $m + n = a$. Thus the set of $K$-types $(m, n)$ with $m + n \le a$ and $m + n \equiv a \pmod 2$ form an $\mathfrak{o}(p, q)$ submodule of $S^a(X^0)$.

We interpret this using our graphical representation of $K$-types as follows. We think of the line $a = x + y$ in $\mathbb{R}^2$ as being a barrier, which moves as $a$ varies. When this moving line passes through integer points in the positive quadrant, it blocks attempts to move past it from those $K$-types (of the correct parity); thus, trapping all $K$-types below it in an $\mathfrak{o}(p, q)$ submodule. This is illustrated by Diagram 2.21.

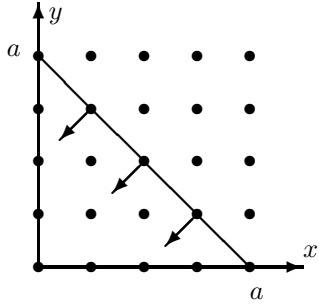

Diagram 2.21

Similar reasoning applies to each of the transition coefficients $A^{\pm\pm}$. Each of them can be thought of as defining a line in $\mathbb{R}^2$, namely, the set of $(x, y) \in \mathbb{R}^2$ that satisfy $A^{\pm\pm}(x, y) = 0$. We will abuse notation and use the symbol $A^{\pm\pm}$ to refer to the line or barrier defined by the zero locus of the linear functional $A^{\pm\pm}(x, y)$ defined in (2.20). These lines move as $a$ varies, and when they pass through integer points (which they all do simultaneously, when $a$ is an integer), they trap the $K$-types on one side and so give rise to a submodule. The position of these barriers for a typical $a > 0$ is given in the Diagram 2.22.

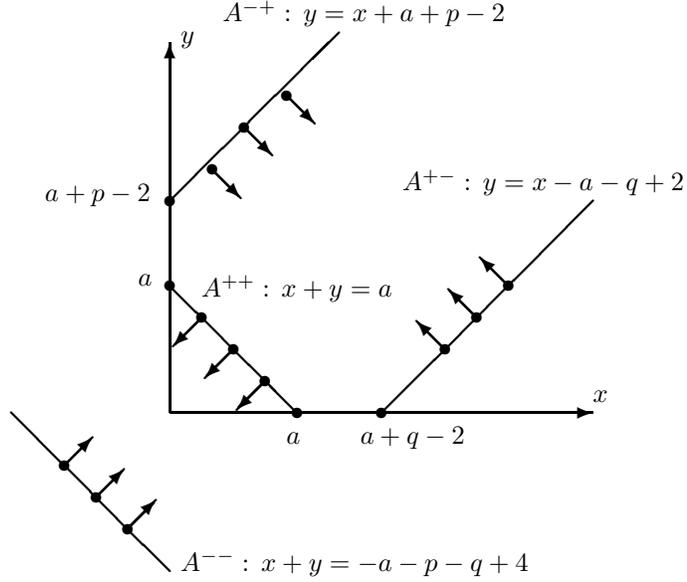

DIAGRAM 2.22

The arrows on each barrier point in the direction of the submodule defined by the barrier. Observe that the barriers $A^{++}$ and $A^{--}$ are parallel and the associated arrows point in opposite directions. Likewise the barriers $A^{+-}$ and $A^{-+}$ are parallel to each other and perpendicular to the other two. As $a$ varies, the barriers $A^{++}$ and $A^{-+}$ move so that they intersect the $x$-axis at points $q-2$ units apart, and the barriers $A^{++}$ and $A^{-+}$ always intersect the $y$-axis at points $p-2$ units apart. Thus $A^{++}$ never intersects $A^{+-}$ or $A^{-+}$ inside the positive quadrant unless $p$ or $q$ equals 2. Also, the barriers $A^{++}$ and $A^{--}$ never intersect the positive quadrant at the same time. They coincide when $a = 2 - \frac{p+q}{2}$. The barriers $A^{+-}$ and $A^{-+}$ also coincide at this value of $a$. For $0 > a > -p - q + 4$, neither of the barriers $A^{++}$ or $A^{--}$ intersects the positive quadrant.

All submodules of $S^a(X^0)$ are accounted for by the barriers $A^{\pm\pm}$ and by the decomposition (2.8). The interactions between these effects yield several patterns of reducibility according to the parities of $p$ and $q$. The details are quite easy to work out. We will summarize the results, with diagrams illustrating typical situations.

The main point to consider is which of the barriers $A^{\pm\pm}$ affect $S^{a+}(X^0)$ and which affect $S^{a-}(X^0)$. For example, suppose $(m, n) \in (\mathbb{Z}_+)^2$ satisfies $A^{++}(m, n) = 0$ while another point $(m', n') \in (\mathbb{Z}_+)^2$ satisfies $A^{+-}(m', n') = 0$. Then

$$\begin{aligned}
(m+n) - (m'+n') &= a - (2m' - a - q + 2) \\
&= 2(a - m' - 1) + q.
\end{aligned}$$

Hence the parity of $(m+n) - (m'+n')$ is the same as the parity of $q$. Thus if $q$ is even, the barriers $A^{++}$ and $A^{+-}$ both affect the same module of the pair $S^{a\pm}$, but if $q$ is odd, one of the barriers affects $S^{a+}$, and the other affects $S^{a-}$. Similarly, the pair of barriers $A^{++}$ and $A^{-+}$ affects the same module if $p$ is even and different ones if $p$ is odd. Also $A^{+-}$ and $A^{-+}$ affect the same or different modules according as $p + q$ is even or odd.

We describe the various possibilities with a collection of diagrams with commentaries. In all the diagrams below, the arrows on the barriers point toward the region of $K$-types that define a submodule.

**Case OO:** $p$ **odd,** $q$ **odd.** In this situation, the barriers $A^{+-}$ and $A^{-+}$ affect the same one of the $S^{a\pm}(X^0)$ and $A^{++}$ or $A^{--}$ affects the other. For the record, we note that $A^{+-}$ and $A^{-+}$ affect $S^{a\varepsilon}(X^0)$ when $\varepsilon = (-1)^{a+1}$, while $A^{++}$ and $A^{--}$ affect $S^{a\varepsilon}(X^0)$ when $\varepsilon = (-1)^a$. We describe the submodule structure for various ranges of values of $a$.

(i) $a \geq 0$.   Here the module affected by the barriers $A^{+-}$ and $A^{-+}$ decomposes into three components: one irreducible submodule and two irreducible quotients, all infinite dimensional, illustrated by Diagram 2.23.

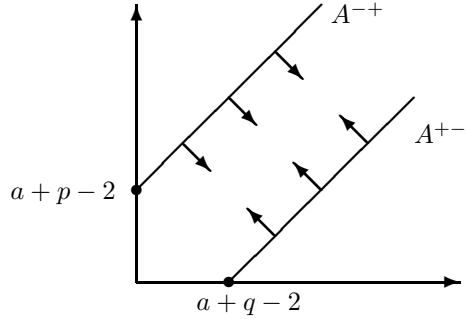

DIAGRAM 2.23

The other of $S^{a\pm}(X^0)$ decomposes into two components: one finite-dimensional submodule and one infinite-dimensional quotient, illustrated by Diagram 2.24.

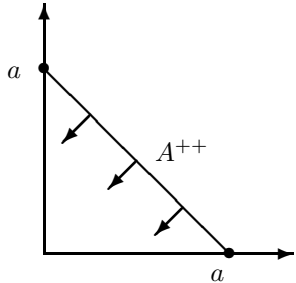

DIAGRAM 2.24

(ii) $0 > a > -[(p+q)/2] + 2$.   In this range, the structure of the module affected by the barriers $A^{+-}$ and $A^{-+}$ remains qualitatively the same as when $a \geq 0$, but the other module is irreducible.

(iii) $a = -[(p+q)/2] + 2$. For this value of $a$, the two barriers $A^{+-}$ and $A^{-+}$ coalesce, and there is a submodule consisting of $K$-types supported on the line $y = x + (p-q)/2$ and two quotient modules consisting of the $K$-types on one side or the other of this line (see Diagram 2.25).

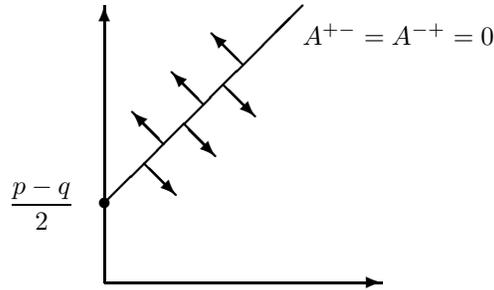

DIAGRAM 2.25

(iv) $a = -(p+q)/2 + 1$. For this value of $a$, the module affected by $A^{+-}$ and $A^{-+}$ (which is $S^{a+}(X^0)$ if $(p-q)/2$ is odd and $S^{a-}(X^0)$ if $(p-q)/2$ is even) decomposes into the direct sum of two irreducible modules. This is because there are no $K$-types in the region between the two barriers (see Diagram 2.26).

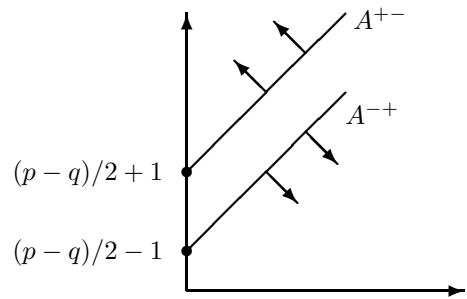

DIAGRAM 2.26

This is the only value of $a$ that gives a decomposable module. The other module remains irreducible.

(v) $a = -[(p+q)/2]$. This situation is dual to case (iii). The constituent with a single line of $K$-types is now a quotient, with barriers trapping the $K$-types on either side of this line in a submodule. See Diagram 2.27.

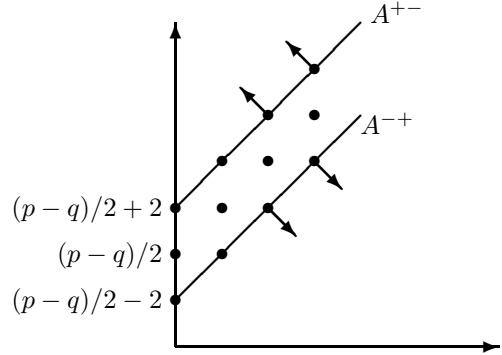

DIAGRAM 2.27

Notice that if the module depicted by Diagram 2.26 in case (iv) is $S^{a+}(X^0)$, the one affected here is $S^{a-}(X^0)$, and vice versa. The other module remains irreducible.

(vi) $-[(p+q)/2] > a \geq -(p+q)+4$.    Here the picture is dual to case (ii). The module affected by $A^{+-}$ and $A^{-+}$ continues to have three constituents, but now it has two submodules and one quotient (see Diagram 2.28). The other module remains irreducible.

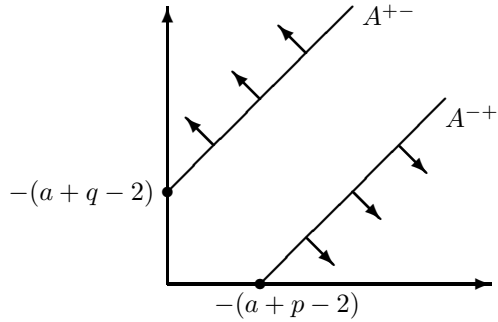

DIAGRAM 2.28

(vii) $a < -(p+q)+4$.    The module affected by $A^{+-}$ and $A^{-+}$ has qualitatively the same structure as in case (vi). The other module now has two constituents: an infinite-dimensional subrepresentation and a finite-dimensional quotient (see Diagram 2.29).

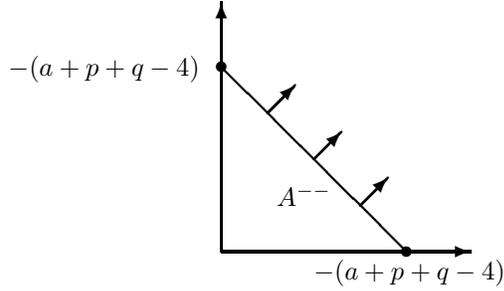

DIAGRAM 2.29

**Case OE: $p$ odd, $q$ even.** In case OE, the barriers $A^{+-}$ and $A^{++}$ affect one of $S^{a\pm}(X^0)$, and the barriers $A^{-+}$ and $A^{--}$ affect the other one. For the record, $A^{+-}$ and $A^{++}$ affect $S^{a\varepsilon}(X^0)$ if $\varepsilon = (-1)^a$, while $A^{-+}$ and $A^{--}$ affect $S^{a\varepsilon}(X^0)$ if $\varepsilon = (-1)^{a+1}$. Here is the description of submodule structure in various ranges of the value of $a$.

(i) $a > 0$. Here the module affected by $A^{++}$ and $A^{+-}$ has three constituents: one finite-dimensional subrepresentation, one infinite-dimensional quotient, and one subquotient in between, as illustrated by Diagram 2.30.

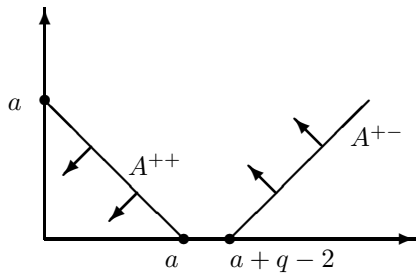

DIAGRAM 2.30

The module affected by $A^{-+}$ has two constituents (see Diagram 2.31).

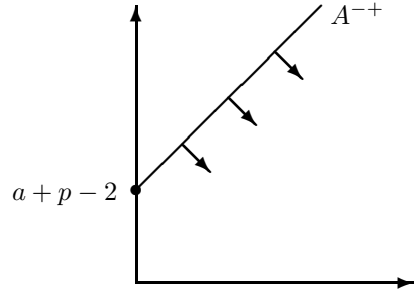

DIAGRAM 2.31

(ii) $0 > a > -(p+q)+4$.    In this range of $a$-values, the module affected by $A^{-+}$ retains the same qualitative structure as it has in case (i), but the other module is no longer affected by $A^{++}$, so it also consists of only two constituents (see Diagram 2.32).

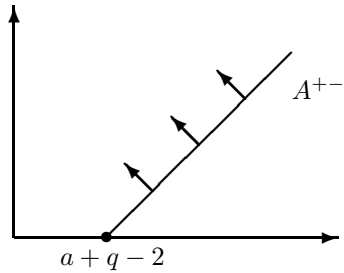

DIAGRAM 2.32

(iii) $a < -(p+q)+4$.    The situation is dual to that of case (i). The module affected by $A^{+-}$ continues qualitatively as in case (ii), while the module affected by $A^{-+}$ is now also affected by $A^{--}$ and has three constituents, including a finite-dimensional quotient (see Diagram 2.33).

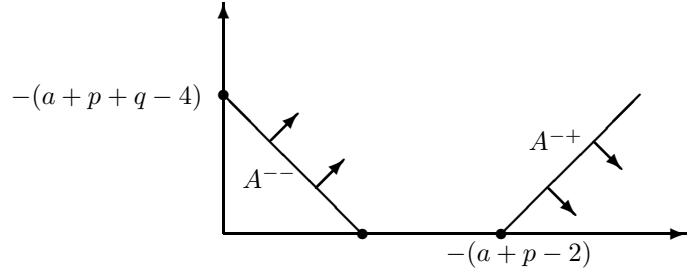

DIAGRAM 2.33

**Case EO:** $p$ **even,** $q$ **odd.** This case, *mutatis mutandis*, is essentially the same as Case OE.

**Case EE:** $p$ **even,** $q$ **even.** When $p$ and $q$ are both even, all the possible barriers affect the same one of $S^{a\pm}(X^0)$ (i.e., $S^{a\varepsilon}(X^0)$ if $\varepsilon = (-1)^a$), and the other is irreducible. The discussion and diagrams below will treat only the reducible summand.

(i) $a > 0$.  Here there are four constituents: one finite-dimensional subrepresentation, two quotients, and one subquotient in between, as illustrated by Diagram 2.34.

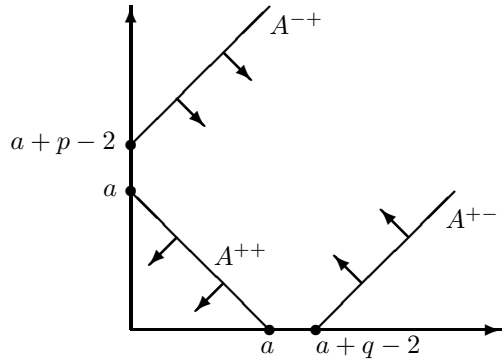

DIAGRAM 2.34

(ii) $0 > a > -(p+q)+4$.  Here, in the range where the barriers $A^{++}$ and $A^{--}$ have no effect, the description of submodules is the same as in cases OO(ii)–(vi). The most interesting things happen in the range $-[(p+q)/2]+2 \geq a \geq -[(p+q)/2]$.

(iii) $a < -(p+q)+4$.  The picture is dual to case (i). There is one finite-dimensional quotient, two infinite-dimensional subrepresentations, and one constituent in between (see Diagram 2.35).

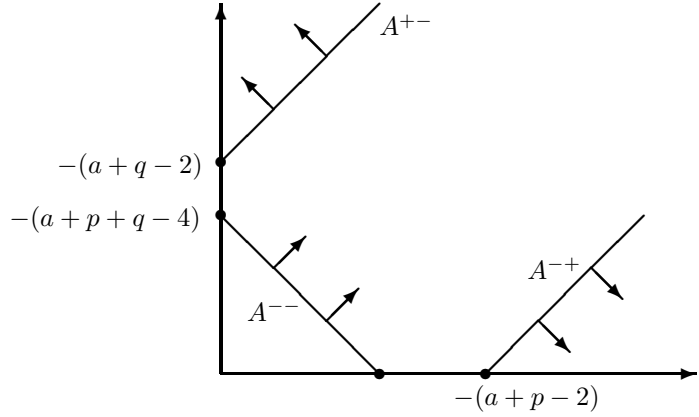

DIAGRAM 2.35

*Remarks.* (a) Although there are considerable variations in the submodule structure of $S^a(X^0)$ according to the parities of $p$ and $q$, there is always a total of five constituents for $a > 0$, one of which is finite dimensional.

(b) The extra complexity when $p + q$ is even of the submodule structure in the range $-[(p+q)/2] + 2 \geq a \geq -[(p+q)/2]$ will acquire additional significance in §3 where unitarity is investigated.

(c) Although we have only described the constituents of $S^a(X^0)$ as "finite dimensional" or "infinite dimensional", the Diagrams 2.23, 2.28, 2.34, and 2.35 certainly suggest that the constituent "in the middle" is smaller than the constituents "in the wings". This difference in size could be made precise in terms of Gelfand-Kirillov dimension [Vo1] or $N$-spectrum [Ho2].

*Degenerate cases $p \geq q = 1$.* We will now indicate how the case of $q = 1$ differs from the picture developed above and also make some supplementary remarks about the case $q = 2$, which is also somewhat exceptional.

When $q = 1$, the light cone $X^0$ is not connected. We can divide it into two relatively open and closed pieces $X^{0\pm}$, called *nappes*, defined by

$$X^{0\pm} = \{(x, y) \in X^0 \mid y = \pm|x|\}.$$

(Note that $y$ is just a real number when $q = 1$.) The two nappes $X^{0\pm}$ of the light cone are each stabilized by a subgroup $\mathrm{O}^+(p, 1)$ of index 2 in $\mathrm{O}(p, 1)$. We see that any element of $S^{a\pm}(X^0)$ is determined by its restriction to $X^{0+}$. Furthermore the function sign on $X^0$, which takes the values $\pm 1$ on $X^{0\pm}$, is clearly invariant under $\mathrm{O}^+(p, 1)$ and is an eigenfunction of $\mathrm{O}(p, 1)$. Evidently multiplying by sign gives an $\mathrm{O}^+(p, 1)$-module isomorphism between the $S^{a\pm}(X^0)$. Thus when $q = 1$, analyzing $S^a(X^0)$ as an $\mathrm{O}(p, 1)$ module is essentially the same as analyzing $S^{a+}(X^0) \simeq S^a(X^{0+})$ as an $\mathrm{O}^+(p, 1)$ module. We will do this, following an abbreviated version of our agenda for the case of $\mathrm{O}(p, q)$, $p, q > 1$.

The map

(2.36) $$\beta : \; X^{0+} \to \mathbb{S}^{p-1} \times \mathbb{R}_+^\times, \qquad \beta(x,y) = (x/y, y),$$

where $x \in \mathbb{R}^p - \{\, 0 \,\}$ and $y = r_p^2(x)^{1/2} \in \mathbb{R}_+^\times$, is a diffeomorphism. The maximal compact subgroup of $\mathrm{O}^+(p,1)$ is $\mathrm{O}(p)$, which acts transitively on the set $\{(x,1) \in X^{0+}\} \simeq \mathbb{S}^{p-1}$. We can define embeddings

(2.37) $$j_{a,m} = j_a : \; \mathcal{H}^m(\mathbb{R}^p) \to S^a(X^{0+}), \qquad j_a(h)(x,y) = h(x)y^{a-m}.$$

Note that since $y > 0$ on $X^{0+}$, complex powers of $y$ make sense. From decomposition (2.36) and the theory of spherical harmonics, we know that the decomposition of $S^a(X^{0+})$ as an $\mathrm{O}(p)$ module is

$$S^a(X^{0+}) \simeq \sum_{m \geq 0} j_a(\mathcal{H}^m(\mathbb{R}^p)).$$

To complete the picture we want to compute the action of $\mathfrak{p}$ (see formulas (2.6)) on the $K$-types $j_a(\mathcal{H}^m(\mathbb{R}^p))$. The analogue of Lemma 2.3 is

$$\left( x_j \frac{\partial}{\partial y} + y \frac{\partial}{\partial x_j} \right) (h y^{a-m})$$
$$= (a-m) h_j^+ y^{a-m-1} + (a+m+p-2) h_j^- y^{a-m+1}$$
$$= (a-m) T_j^+(h) y^{a-m-1} + (a+m+p-2) T_j^-(h) y^{a-m+1},$$

where

$$h_j^+ = x_j h - (2m+p-2)^{-1} \frac{\partial h}{\partial x_j} r_p^2,$$

$$h_j^- = (2m+p-2)^{-1} \frac{\partial h}{\partial x_j}, \qquad h \in \mathcal{H}^m(\mathbb{R}^p).$$

We note that $h_j^-$ is harmonic, while $h_j^+$ is the projection of $x_j h$ into $\mathcal{H}^{m+1}(\mathbb{R}^p)$. The principles behind this calculation are discussed in §6.1. It follows easily that $S^a(X^{0+})$ is always irreducible except when $a$ is an integer, with either $a \geq 0$ or $a \leq -p+1$, in which case $S^a(X^{0+})$ has two constituents, one of which is finite dimensional. In what we hope will be a self-evident analogy with the situation for $q > 1$, we illustrate this by Diagram 2.38.

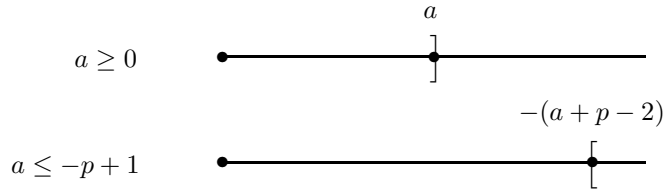

DIAGRAM 2.38

For $-p+1 < a < 0$, the representations $S^a(X^{0+})$ are irreducible.

Although the general picture presented in Diagrams 2.23–2.35 is valid when $q = 2$, this case is also somewhat exceptional. Although the $\mathcal{H}^n(\mathbb{R}^2)$ are nonzero for all $n \geq 1$, their dimensions, instead of growing with $n$, remain fixed at 2 (for $n = 0$, $\mathcal{H}^n(\mathbb{R}^q)$ is the trivial representation for all $q$ and has dimension 1). More significantly, whereas the representations $\mathcal{H}^n(\mathbb{R}^q)$, $n \geq 1$, $q \geq 3$, of O$(q)$ remain irreducible when restricted to the special orthogonal group SO$(q)$, the representations $\mathcal{H}^n(\mathbb{R}^2)$, $n \geq 1$, decompose into two one-dimensional representations when restricted to SO(2), which is isomorphic to $\mathbb{T}$, the unit circle. This suggests that one consider varying the picture presented above by restricting the representations $S^a(X^0)$ to the identity component O$^0(p,q)$ of O$(p,q)$. For $p, q \geq 3$ this has no effect. When $q = 2$ (and $p > 2$), however, the pictures presented above for O$(p,2)$ must be "folded out" along the $x$-axis to yield pictures appropriate for O$^0(p,q)$. Precisely the Diagrams 2.30, 2.31, and 2.33 now should look as in Diagrams 2.39, 2.40, and 2.41 respectively.

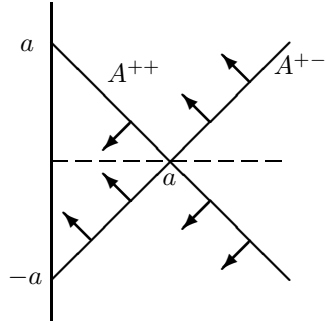

DIAGRAM 2.39

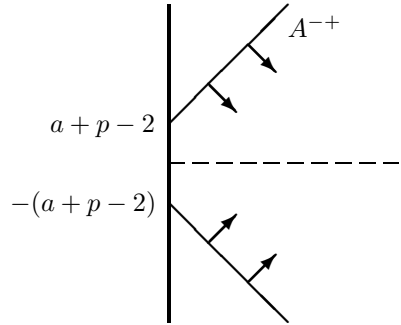

DIAGRAM 2.40

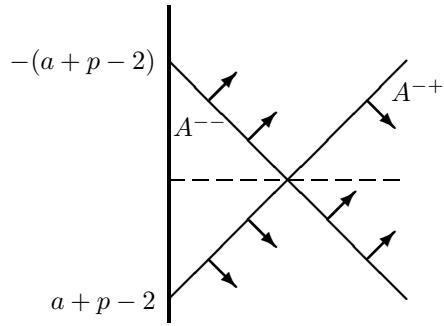

DIAGRAM 2.41

We observe that in Diagram 2.22, the gap between the points where $A^{++}$ and $A^{+-}$ intersect the $x$-axis is $q-2$, so when $q=2$, these points coincide, and so when reflected across the $x$-axis, these barriers make crossing straight lines. Similar remarks apply to $A^{--}$ and $A^{-+}$.

If $p$ is even, then the relevant diagrams are 2.34 and 2.35. The modified version of Diagram 2.34 is Diagram 2.42.

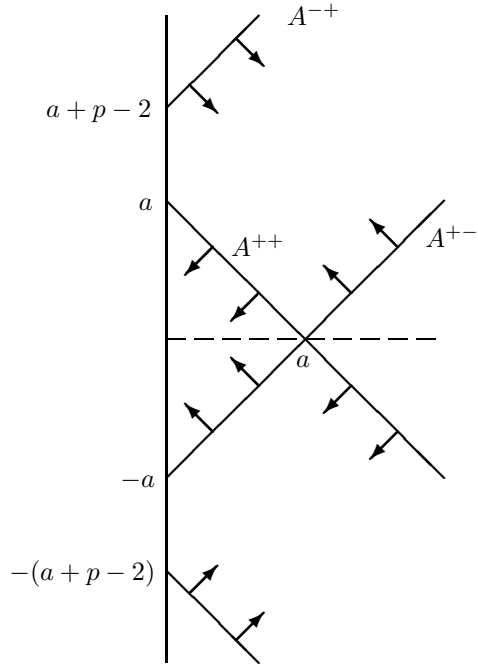

DIAGRAM 2.42

Thus in this case, we have six constituents in all, with *Hasse diagram* [Sta] as illustrated by Diagram 2.43.

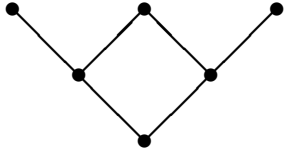

<div align="center">DIAGRAM 2.43</div>

We omit the diagram corresponding to Diagram 2.35, which would be dual to Diagram 2.42. The main effect of the restriction is that representations that have no $K$-types along the $x$-axis, that is, that contain no fixed vectors for O(2), break up into two summands symmetrically placed across the $x$-axis. However, representations whose $K$-type region includes part of the $x$-axis remain irreducible. Reinforcing this splitting is the fact that when $q = 2$, the barriers $A^{++}$ and $A^{+-}$ (and the barriers $A^{--}$ and $A^{-+}$) intersect on the $x$-axis, leaving no gap for $K$-types of the intermediate constituent in Diagram 2.30. Hence this constituent becomes two on restriction to $O^0(p, 2)$. The representations of $O^0(p, 2)$ whose $K$-types are restricted to only one side of the $x$-axis are of the type known as *holomorphic* or *highest weight modules* [EHW, Br].

If $p = 2$ and $q = 2$, then Diagram 2.42 should be further folded out along the $y$-axis. Since O(2, 2) is just isogenous to SL(2, $\mathbb{R}$) × SL(2, $\mathbb{R}$), we will not dwell on this example.

## 3. UNITARITY

In this section we will investigate which of the subquotients of the $S^a(X^0)$ described in §2 can define unitary representations, that is, possess an O($p, q$)-invariant positive-definite Hermitian form. As with the calculations of submodule structure, the results are primarily based on the formulas of Lemma 2.3 for the action of $\mathfrak{p} \subseteq \mathfrak{o}(p, q)$.

Since the $K$-types $j_a(\mathcal{H}^m(\mathbb{R}^p) \otimes \mathcal{H}^n(\mathbb{R}^q))$ are irreducible modules for the compact group O($p$) × O($q$), each one admits an O($p$) × O($q$)-invariant Hermitian inner product, and this is unique up to multiples. For this discussion we will fix an inner product $\langle \cdot, \cdot \rangle_{m,n}$ on each space $\mathcal{H}^m(\mathbb{R}^p) \otimes \mathcal{H}^n(\mathbb{R}^q)$ by the geometric recipe

$$(3.1) \qquad \langle \phi_1, \phi_2 \rangle_{m,n} = \int_{\mathbb{S}^{p-1} \times \mathbb{S}^{q-1}} \phi_1(z) \overline{\phi_2(z)} \, dz,$$

where $\phi_j \in \mathcal{H}^m(\mathbb{R}^p) \otimes \mathcal{H}^n(\mathbb{R}^q)$ and $dz$ indicates the O($p$)×O($q$)-invariant probability measure on $\mathbb{S}^{p-1} \times \mathbb{S}^{q-1}$. We will push forward this inner product by $j_a$ to obtain an inner product on the $K$-type $j_a(\mathcal{H}^m(\mathbb{R}^p) \otimes \mathcal{H}^n(\mathbb{R}^q))$. We call the collection of inner products so obtained the *standard inner products* on the $K$-types.

If $\langle \cdot, \cdot \rangle$ is an O($p, q$)-invariant Hermitian form on $S^a(X^0)$ (or on some constituent of $S^a(X^0)$), the restriction of $\langle \cdot, \cdot \rangle$ to any $K$-type must be a multiple of the standard inner product on that $K$-type:

$$(3.2) \qquad \langle j_a(f_1), j_a(f_2) \rangle = c_{m,n} \langle f_1, f_2 \rangle_{m,n},$$

where $f_j \in \mathcal{H}^m(\mathbb{R}^p) \otimes \mathcal{H}^n(\mathbb{R}^q)$ and $c_{m,n}$ is an appropriate real number, which must be positive if $\langle \cdot, \cdot \rangle$ is to be positive definite. Furthermore the $K$-types must be mutually orthogonal, so that $\langle \cdot, \cdot \rangle$ is determined by the numbers $c_{m,n}$.

According to an early theorem of Harish-Chandra [HC2], a positive-definite inner product $\langle \cdot, \cdot \rangle$ on $S^a(X^0)$ will be $\mathrm{O}(p,q)$ invariant if and only if the operators from $\mathfrak{o}(p,q)$ are formally skew adjoint with respect to $\langle \cdot, \cdot \rangle$. Equation (3.2) guarantees this will be so for operators from $\mathfrak{o}(p) \oplus \mathfrak{o}(q)$. Thus to determine invariance of $\langle \cdot, \cdot \rangle$ we just need to check skew-adjointness of the operators from the subspace $\mathfrak{p} \subseteq \mathfrak{o}(p,q)$. From formulas (2.18) and (3.2), we see that this amounts to the equations

$$(3.3) \qquad \begin{aligned} \text{(a)} \quad 0 =& (a - m - n)c_{m+1,n+1}\langle T_{m,n}^{++}(z \otimes \phi), \phi' \rangle_{m+1,n+1} \\ &+ (\overline{a} + m + n + p + q - 2)c_{m,n}\langle \phi, T_{m+1,n+1}^{--}(z \otimes \phi') \rangle_{m,n}, \\ \text{(b)} \quad 0 =& (a - m + n + q - 2)c_{m+1,n-1}\langle T_{m,n}^{+-}(z \otimes \phi), \phi'' \rangle_{m+1,n-1} \\ &+ (\overline{a} + m - n + p)c_{m,n}\langle \phi, T_{m+1,n-1}^{-+}(z \otimes \phi'') \rangle_{m,n}, \end{aligned}$$

where $z \in \mathfrak{p}$, $\phi \in \mathcal{H}^m(\mathbb{R}^p) \otimes \mathcal{H}^n(\mathbb{R}^q)$, $\phi' \in \mathcal{H}^{m+1}(\mathbb{R}^p) \otimes \mathcal{H}^{n+1}(\mathbb{R}^q)$, and $\phi'' \in \mathcal{H}^{m+1}(\mathbb{R}^p) \otimes \mathcal{H}^{n-1}(\mathbb{R}^q)$. Here $\overline{a}$ means the complex conjugate of $a$.

Equations (3.3) can be substantially simplified by an elementary observation, which will be justified in §6.

**Lemma 3.1.** *When* $\mathrm{Re}\, a = -[(p+q)/2] + 1$, *that is, for* $a = -[(p+q)/2] + 1 + is$, $s \in \mathbb{R}$, *the Hermitian inner product which equals the standard inner product on each $K$-type, is* $\mathrm{O}(p,q)$ *invariant.*

*Remark.* More generally, one can show that if $\phi_1$ in formula (3.1) is taken to belong to $S^\alpha(X^0)$ and $\phi_2$ is taken from $S^{\overline{\alpha'}}(X^0)$, where $\alpha + \alpha' = 2 - (p+q)$, then (3.1) defines an $\mathrm{O}(p,q)$-invariant pairing between $S^\alpha(X^0)$ and $S^{\overline{\alpha'}}(X^0)$, which are thus canonically mutual Hermitian duals. This duality explains the symmetry in structure among the $S^\alpha(X^0)$, as discussed in §2.

*Proof.* See §6.1. □

If we take $a = -[(p+q)/2] + 1$ in formula (3.3), we see that Lemma 3.1 implies

(a) $\langle T_{mn}^{++}(z \otimes \phi), \phi' \rangle_{m+1,n+1} = \langle \phi, T_{m+1,n+1}^{--}(z \otimes \phi') \rangle_{m,n}$,

(b) $\langle T_{mn}^{+-}(z \otimes \phi), \phi'' \rangle_{m+1,n-1} = \langle \phi, T_{m+1,n-1}^{-+}(z \otimes \phi'') \rangle_{m,n}$

for all possible choices of $z$, $\phi$, $\phi'$, $\phi''$ as in (3.3). Given these relations, a reinspection of formulas (3.3) shows they are equivalent to the purely numerical equations

$$\begin{cases} (a - m - n)c_{m+1,n+1} + (\overline{a} + m + n + p + q - 2)c_{m,n} = 0, \\ (a - m + n + q - 2)c_{m+1,n-1} + (\overline{a} + m - n + p)c_{m,n} = 0 \end{cases}$$

or

$$(3.4) \qquad \begin{cases} A^{++}(m,n)c_{m+1,n+1} + \overline{A^{--}(m+1,n+1)}c_{m,n} = 0, \\ A^{+-}(m,n)c_{m+1,n-1} + \overline{A^{-+}(m+1,n-1)}c_{m,n} = 0. \end{cases}$$

Let us explore the implications of equations (3.4). Consider first the case when the $S^{a\pm}(X^0)$ are irreducible. Then the quotient

$$\frac{m+n-a}{\overline{a}+m+n+p+q-2}$$

must at least be real for all possible choices of $m$ and $n$. Setting

$$\tilde{a} = a + \frac{p+q}{2} - 1,$$

we can rewrite this ratio as

$$(3.5) \qquad \frac{m+n+\delta-\tilde{a}}{m+n+\delta+\overline{\tilde{a}}} = \frac{(m+n+\delta-\operatorname{Re}\tilde{a})-i\operatorname{Im}\tilde{a}}{(m+n+\delta+\operatorname{Re}\tilde{a})-i\operatorname{Im}\tilde{a}},$$

where

$$\delta = \frac{p+q}{2} - 1.$$

It is easy to convince oneself that the ratio (3.5) can be real for large $m, n$ only if either $\operatorname{Re}\tilde{a} = 0$ or $\operatorname{Im}\tilde{a} = 0$. The case of $\operatorname{Re}\tilde{a} = 0$ is covered by Lemma 3.1. We will call the set of $a$ with $\operatorname{Re}\tilde{a} = 0$ the *unitary axis*. If $\operatorname{Im}\tilde{a} = \operatorname{Im}a \neq 0$, we know (see Proposition 2.4) that $S^{a\pm}(X^0)$ is irreducible; hence for it to be Hermitian, expression (3.5) would have to be real for all $m, n$. This leads us to the following conclusion.

**Lemma 3.2.** *In order for $S^a(X^0)$ to have any Hermitian constituents, either $a$ must be on the unitary axis, or $a$ must be real. Conversely, for $a$ on the unitary axis, $S^a(X^0)$ is unitary, and for $a$ real, all components of $S^a(X^0)$ are Hermitian (i.e., they possess an invariant Hermitian form).*

The final assertion follows from the observation that if $a$ is real then we can start at a given $K$-type $(m_0, n_0)$ in any constituent of $S^a(X^0)$, set $c_{m_0,n_0} = 1$, and use equations (3.4) to determine $c_{m,n}$ for other $K$-types in that constituent. It should be checked that this procedure is consistent, i.e., that it results in a well-defined $c_{m,n}$ on each $K$-type. For this, it suffices to check that the transitions around squares, like $(m,n) \longrightarrow (m+1, n+1) \longrightarrow (m+2, n)$ and $(m,n) \longrightarrow (m+1, n-1) \longrightarrow (m+2, n)$, predict the same value for the ratio of the start and finish (here $c_{m+2,n}/c_{m,n}$). This is a short computation. Here, as in the analysis of irreducibility, it is germane to note that the transition coefficients $A^{++}$ and $A^{--}$ are unchanged by the transitions $(m,n) \longrightarrow (m,n) \pm (1, -1)$; likewise $A^{+-}$ and $A^{-+}$ are unchanged under $(m,n) \longrightarrow (m,n) \pm (1, 1)$.

Our goal is now to investigate for which constituents of $S^a(X^0)$ the process described in the previous paragraph will yield positive values for all $c_{m,n}$ in the constituent. Looking again at equations (3.4), we can easily see and formulate the following rule.

**Criterion 1.** Let $a$ be real. A given constituent of $S^a(X^0)$ is unitary if and only if the ratios

$$(3.6) \qquad (i) \quad \frac{A^{++}(m,n)}{A^{--}(m+1,n+1)} \qquad \text{and} \qquad (ii) \quad \frac{A^{+-}(m,n)}{A^{-+}(m+1,n-1)}$$

are negative for all pairs of $K$-types $(m, n)$ and $(m + 1, n + 1)$ or $(m, n)$ and $(m + 1, n - 1)$ contained in the constituent.

We proceed to investigate the nature of this criterion in order to obtain an explicit list of the unitary constituents.

Consider again our graphical representation of the $K$-types as integer points in the positive quadrant, together with barrier lines corresponding to the zeroes of the linear functions $A^{\pm\pm}$. We illustrate the situation for a typical $a > 0$ in Diagram 3.7.

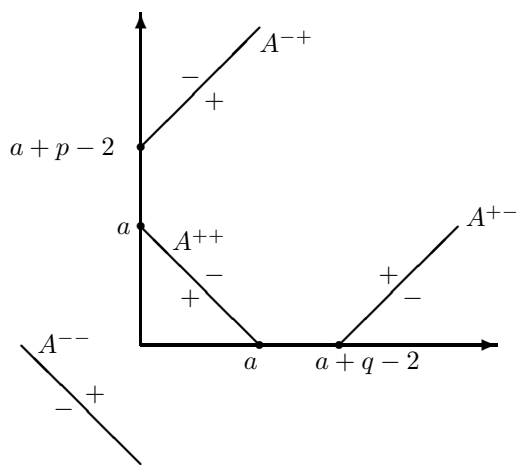

Diagram 3.7

In this diagram we have indicated, for each transition coefficient, the side of its barrier where it is positive and the side where it is negative. We see that, in the positions illustrated, both $A^{+-}$ and $A^{-+}$ are positive in the region between their barriers and likewise for $A^{++}$ and $A^{--}$. When $a < 0$ and the $A^{+-}$ barrier is above the $A^{-+}$ barrier, both $A^{+-}$ and $A^{-+}$ are negative in the region between barriers. In general, we may say that $A^{+-}$ and $A^{-+}$ have the same sign in the region between barriers, and have opposite signs in the regions outside the barriers, similarly for $A^{++}$ and $A^{--}$. Thus, if $(m, n)$ and $(m + 1, n - 1)$ both lie between the $A^{+-}$ and $A^{-+}$ barriers, or more precisely if $(m, n)$ lies above the $A^{+-}$ barrier and $(m+1, \ n-1)$ lies below the $A^{-+}$ barrier (or vice versa), then the ratio (3.6)(ii) is positive and unitarity cannot hold. On the other hand, if $(m, n)$ lies above the $A^{+-}$ barrier and $(m+1, \ n-1)$ also lies above the $A^{-+}$ barrier (or both lie below), then the transition from $(m, n)$ to $(m+1, \ n-1)$ is consistent with unitarity. Similar remarks apply to the pair $A^{++}$, $A^{--}$.

We can apply these considerations to determine the unitarity of constituents of the $S^a(X^0)$. As for the discussion of submodule structure in §2, it is convenient to distinguish cases according to the parities of $p$ and $q$. Before beginning we note that, using the notation $\tilde{a}$ as in formula (3.5), the same argument that establishes Lemma 3.1 shows that the modules $S^{\tilde{a}}(X^0)$ and $S^{-\tilde{a}}(X^0)$ are mutual Hermitian duals. This will be shown in §6.1. It follows that in considering questions of unitarity, it will suffice to take $\tilde{a} \geq 0$.

**Case OO: $p$ odd, $q$ odd.** From the unitary axis ($\tilde{a} = 0$) until $\tilde{a} = [(p+q)/2] - 1$ (i.e., $a = 0$), neither of the barriers $A^{++}$ or $A^{--}$ intersects the positive quadrant. Hence unitarity depends only on the barriers $A^{+-}$ and $A^{-+}$. Their position when $\tilde{a} = 0$ is illustrated in Diagram 2.26. Let $\varepsilon_0$ denote the choice of $\pm$ such that $S^{a\varepsilon_0}(X^0)$ is reducible at $\tilde{a} = 0$, and let $\varepsilon_1$ denote the other choice. For the record, we note that $\varepsilon_0$ is $+$ if $(p-q)/2$ is odd and is $-$ if $(p-q)/2$ is even.

At $\tilde{a} = 0$ none of the $K$-types of $S^{a\varepsilon_0}(X^0)$ lie between the barriers $A^{+-}$ and $A^{-+}$, which are located on the two lines $L_{\pm 1}$ defined by equations $y - x = [(p-q)/2] \pm 1$; however, there are $K$-types lying on these two lines. As $\tilde{a}$ increases from 0, the barrier $A^{+-}$ moves to the right and $A^{-+}$ moves left. Hence, for $0 < \tilde{a} < 2$, transitions between $L_1$ and $L_{-1}$ cross both barriers and so violate unitarity. For $\tilde{a} \geq 2$, transitions between $L_1$ and $L_{-1}$ cross neither barrier and so again violate unitarity. Hence the full module $S^{a\varepsilon_0}(X^0)$ is never unitary when $\tilde{a} > 0$, and at points of reducibility the constituent with $K$-types between the barriers cannot be unitary. The other two constituents, however, are easily seen to be unitary, since all their transitions take place on the same side of both barriers.

The module $S^{a\varepsilon_1}(X^0)$ is irreducible at $\tilde{a} = 0$ and remains so for $0 < \tilde{a} < 1$. For this range of $a$, the only $K$-types between the barriers are those on the line $L_0$ defined by $y - x = (p-q)/2$. One checks that transitions in either direction from $L_0$ are consistent with unitarity. At $\tilde{a} = 1$ the barriers coincide on $L_0$ (see Diagram 2.25), and the $K$-types on this line span a constituent of $S^{a\varepsilon_1}(X^0)$. There are no transitions transverse to the barriers from $L_0$, so this constituent is unitary; likewise the other constituents are unitary, so all three are unitary. For $\tilde{a} > 1$ the line $L_0$ is on the opposite side of the barriers $A^{+-}$ and $A^{-+}$, from where it was when $\tilde{a} < 1$, so the transitions from $L_0$ violate unitarity. Hence unitarity for $S^{a\varepsilon_1}(X^0)$ stops at $\tilde{a} = 1$. Also, for points of reducibility beyond $\tilde{a} = 1$, the same reasoning applies as to $S^{a\varepsilon_1}(X^0)$.

In the range $a \geq 0$, when there may be a finite-dimensional component, one easily checks that the trivial representation at $a = 0$ is unitary, but no other finite-dimensional representations are so.

In sum, the unitarity properties of the $S^{a\pm}(X^0)$ in this case are summed up in Diagrams 3.8 and 3.9. For Diagram 3.9, we note that (a) and (b) depict $S^{a\varepsilon}(X^0)$ when $\varepsilon = (-1)^{a+1}$, while (c) and (d) depict $S^{a\varepsilon}(X^0)$ when $\varepsilon = (-1)^a$.

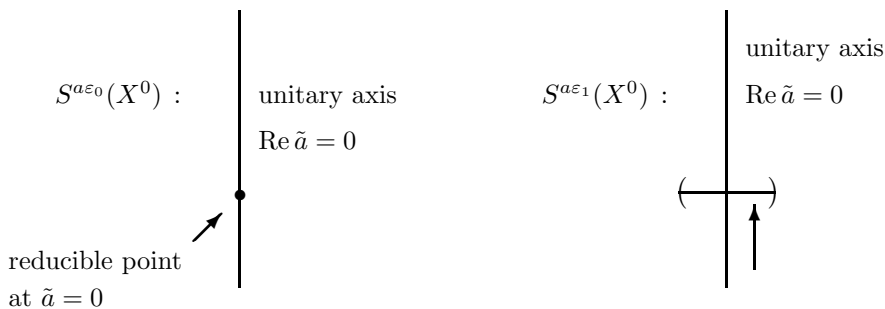

DIAGRAM 3.8. Parameter values for unitarity of the full modules $S^{a\pm}(X^0)$.

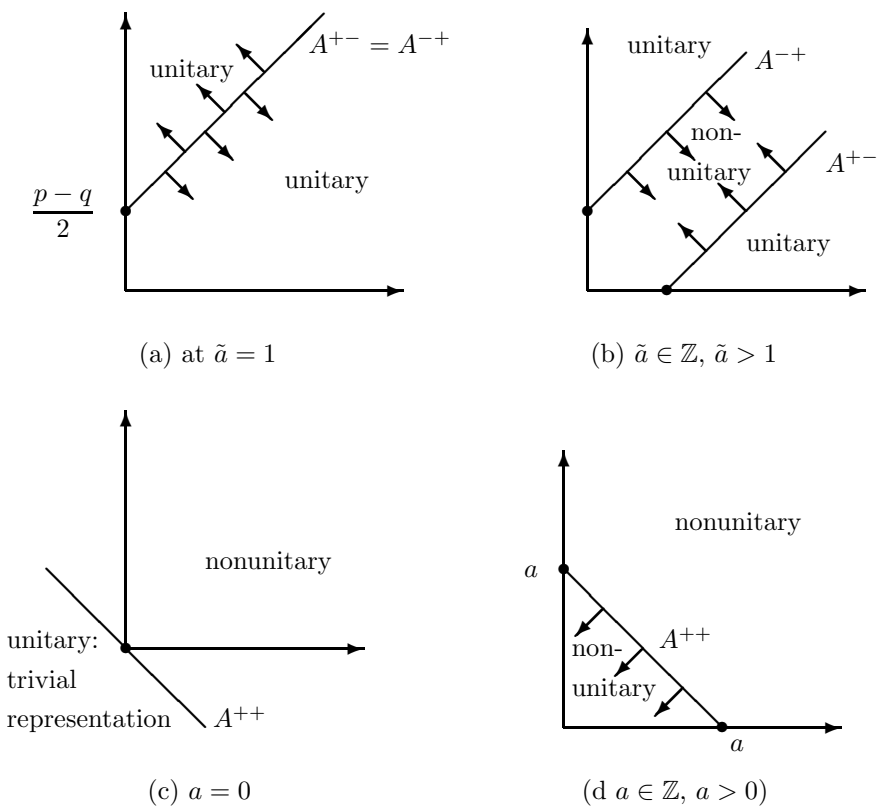

DIAGRAM 3.9. Unitarity at points of reducibility.

**Case OE: $p$ odd, $q$ even.** In this case, when $\tilde{a} = 0$, the barriers do not pass through integral points, so the representations $S^{-\delta\pm}(X^0)$ are irreducible (recall that $\delta = [(p+q)/2] - 1$). Each has $K$-types between the barriers: one on the line $L_{1/2}$ defined by $y - x = [(p-q)/2] + \frac{1}{2}$ and the other on the line $L_{-1/2}$, defined by $y - x = [(p-q)/2] - \frac{1}{2}$ (see Diagram 3.10). As $\tilde{a}$ increases from 0, the barriers approach these lines and coincide with them for $\tilde{a} = \frac{1}{2}$. One checks that for $0 < \tilde{a} < \frac{1}{2}$ transitions from these lines is consistent with unitarity, but as $\tilde{a}$ passes $\frac{1}{2}$, the barriers move past the lines $L_{\pm 1/2}$, and then transition from these lines does violate unitarity. Hence for $\tilde{a} > \frac{1}{2}$, the full modules $S^{\tilde{a}\pm}(X^0)$ can no longer be unitary. At points of reducibility for $\tilde{a} > \frac{1}{2}$, one can check that the constituent with $K$-types to the right of $A^{+-}$ and the one with $K$-types to the left of $A^{-+}$ is unitary. Again the only unitary finite-dimensional unitary constituent is the trivial representation. Diagrams 3.11 and 3.12 illustrate the main conclusions. For Diagram 3.12, we note that $A^{++}$ and $A^{+-}$ affect $S^{a\varepsilon}(X^0)$ if $\varepsilon = (-1)^a$, while $A^{-+}$ and $A^{--}$ affect $S^{a\varepsilon}(X^0)$ if $\varepsilon = (-1)^{a+1}$.

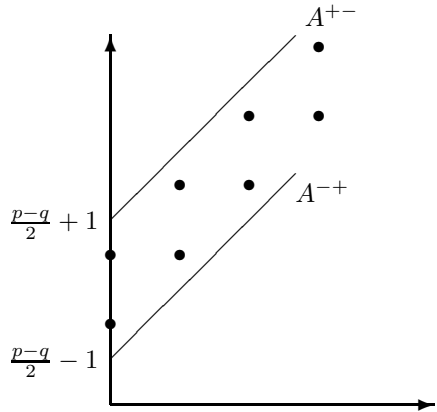

DIAGRAM 3.10. $\tilde{a} = 0$.

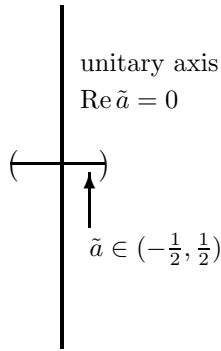

DIAGRAM 3.11. Parameter values of unitarity of the full modules $S^{a\pm}(X^0)$.

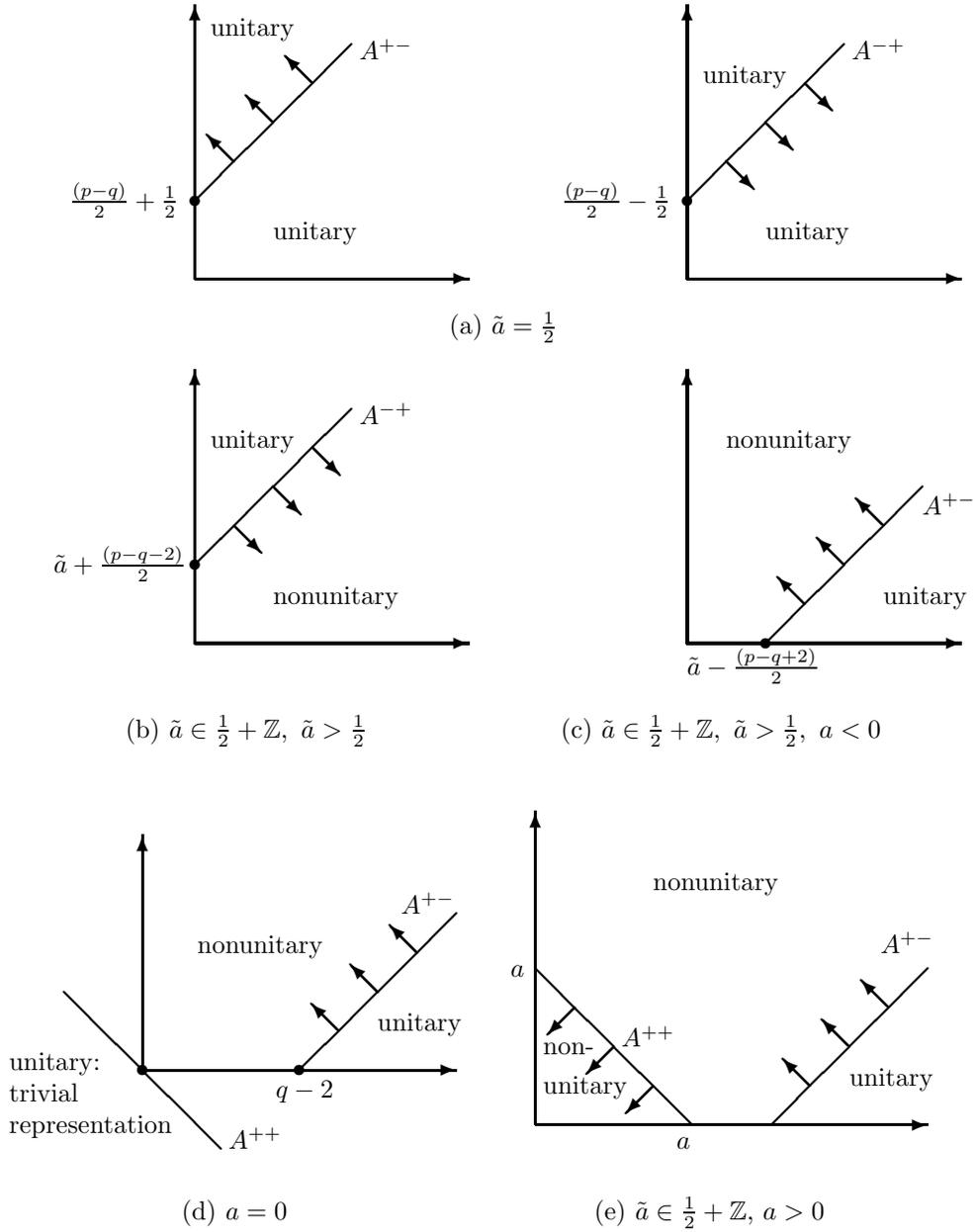

(a) $\tilde{a} = \frac{1}{2}$

(b) $\tilde{a} \in \frac{1}{2} + \mathbb{Z}, \ \tilde{a} > \frac{1}{2}$

(c) $\tilde{a} \in \frac{1}{2} + \mathbb{Z}, \ \tilde{a} > \frac{1}{2}, \ a < 0$

(d) $a = 0$

(e) $\tilde{a} \in \frac{1}{2} + \mathbb{Z}, \ a > 0$

DIAGRAM 3.12. Unitarity of constituents at points of reducibility.

**Case EO: $p$ even, $q$ odd.** This case is essentially the same as case OE. We will not give the details.

**Case EE: $p$ even, $q$ even.** Here the discussion is essentially the same as for Case OO, except that, when $a \geq 0$, at a point of reducibility, three barriers affect the same $S^{a\pm}(X^0)$, and the other one is irreducible. We draw the pictures illustrating this situation (see Diagrams 3.13(a) and (b)).

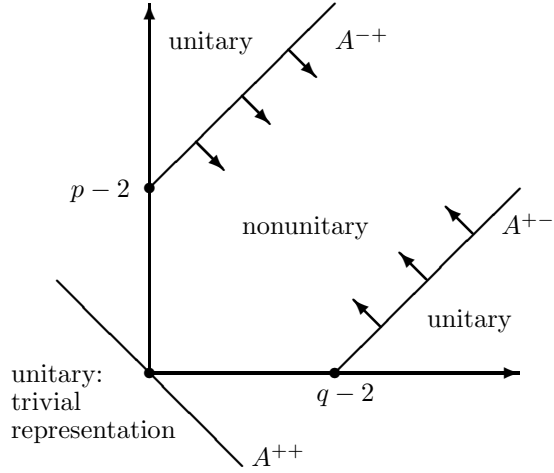

(a) $a = 0$

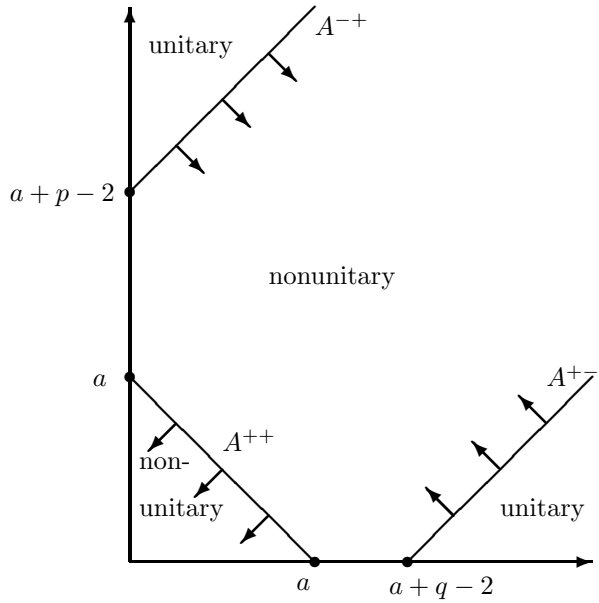

(b) $a \in \mathbb{Z},\ a > 0$

Diagram 3.13

*Remarks.* (a) It is interesting to note that both summands $S^{a\pm}(X^0)$ of $S^a(X^0)$ "see" all the barriers via their Hermitian forms, whether or not these barriers cause reducibility in one summand or the other. In fact, they see the barriers even when the barriers do not pass through integral points at all.

(b) The unitary representation with $K$-types distributed along a single line which occurs at $\tilde{a} = 1$ when $p + q$ is even [see Diagram 3.9(a)] is particularly interesting. It is paired with the trivial representation in the local correspondence for the dual pair $(\mathrm{O}(p,q), \mathrm{SL}(2,\mathbb{R}))$. It is the smallest possible nontrivial unitary representation of $\mathrm{O}(p,q)$ and may be regarded as a "quantization of the minimal nilpotent orbit" in the argot of geometric quantization [EPW]. It turns up in physics in connection with the quantum Coulomb problem [Ab, AFR, DGN, Fr, On].

(c) The representations $S^{a+}(X^0)$ are "spherical" representations, meaning that they contain the trivial representation of the maximal compact subgroup, and the $S^{a-}(X^0)$ are nonspherical representations. Notice that, for each of these families or "series" of representations, the only point on the unitary axis that can be reducible is the point on the real axis, namely, $\tilde{a} = 0$. When $p + q$ is even, exactly one of the two series is reducible at $\tilde{a} = 0$, while when $p + q$ is odd neither series is. The group $\mathrm{SL}(2,\mathbb{R})$ has two principal series of representations, analogous to the $S^{a\pm}(X^0)$: one spherical series and one nonspherical series. The spherical series is irreducible on the real axis, and the nonspherical series is reducible. Also, as was well known and has been shown again here, the unitary spherical principal series for $\mathrm{O}(n,1)$ is always irreducible. These examples perhaps created some expectation that spherical series were more likely to be irreducible than other series. In our example, however, when $p + q$ is even, either the spherical or the nonspherical series may be reducible at $\tilde{a} = 0$, according to the parity of $(p - q)/2$, while when $p + q$ is odd, neither series is reducible.

We note also that when $a$ is a positive integer, the trivial $K$-type is always contained in the central strip where nonunitarity reigns. Thus outside the unitary axis and its immediate neighborhood (the complementary series), there are only finitely many unitary spherical representations.

*Degenerate cases $p \geq q = 1$.* In closing, we record the analogous pictures for $\mathrm{O}(p,1)$.

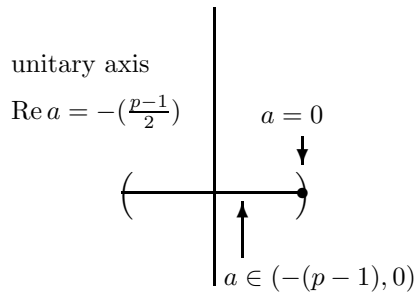

DIAGRAM 3.14. Parameter values for unitarity of the full modules $S^{a\pm}(X^0)$.

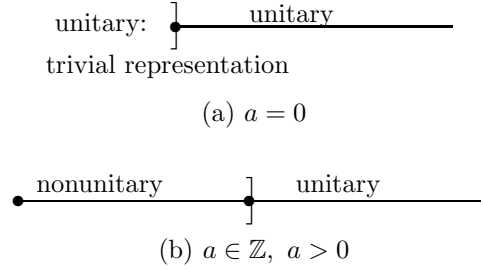

Diagram 3.15. Unitarity at points of reducibility.

The reason that the unitary representations can live so far off the unitary axis is that since the $K$-types now live along a single line, the $A^{+-}$ and $A^{-+}$ barriers no longer affect unitarity—only the barriers which cut off finite-dimensional representations, and these do not come into play until the finite-dimensional representations begin appearing at $a = 0$.

## 4. The U$(p,q)$-module structure of homogeneous functions on light cones

**4.1. The U$(p,q)$ modules $S^{\alpha,\beta}(X^0)$.** We now wish to adapt the analysis given for O$(p,q)$ in §§2 and 3 to the indefinite unitary groups U$(p,q)$.

Let $\mathbb{C}^{p+q} = \mathbb{C}^p \oplus \mathbb{C}^q$ denote complex $(p+q)$-dimensional space, parametrized by tuples $(z,w)$, with $z = (z_1, z_2, \ldots, z_p)$ in $\mathbb{C}^p$ and $w = (w_1, w_2, \ldots, w_q)$ in $\mathbb{C}^q$. We endow $\mathbb{C}^{p+q}$ with the Hermitian form

$$((z,w),(z',w'))_{p,q} = (z,z')_p - (w,w')_q = \sum_{j=1}^p z_j \overline{z_j'} - \sum_{j=1}^q w_j \overline{w_j'}$$

of signature $(p,q)$. Let U$(p,q) \subset$ GL$(p+q, \mathbb{C})$ denote the group of linear isometries of the form $(\cdot, \cdot)_{p,q}$. Observe that the subgroup U$(p) \times$ U$(q)$, which preserves the $z$ variables and the $w$ variables, is a compact subgroup of U$(p,q)$. (In fact, it is a maximal compact subgroup.) Let

$$(4.1) \qquad \begin{aligned} X^0 &= \left\{ \eta \in \mathbb{C}^{p+q} - \{0\} \,\middle|\, (\eta, \eta)_{p,q} = 0 \right\} \\ &= \left\{ (z,w) \in \mathbb{C}^{p+q} - \{0\} \,\middle|\, (z,z)_p = (w,w)_q \right\} \end{aligned}$$

be the light cone or variety of null vectors. By Witt's Theorem [Ja] or by elementary means, one can see that $U(p,q)$ acts transitively on $X^0$. We wish to study the action of $U(p,q)$ on appropriate spaces of functions on $X^0$.

If we regard $\mathbb{C}^{p+q}$ as $\mathbb{R}^{2p+2q}$ by means of coordinates $x_l, y_l$, where

$$(4.2) \qquad \begin{aligned} z_j &= x_{2j-1} + \boldsymbol{i}x_{2j}, & 1 \le j \le p, \\ w_k &= y_{2k-1} + \boldsymbol{i}y_{2k}, & 1 \le k \le q, \end{aligned}$$

then the real part of the form $(\cdot,\cdot)_{p,q}$ becomes the indefinite inner product of signature $(2p, 2q)$, the light cone $(4.1)$ is the same as the light cone for $\mathbb{R}^{2p+2q}$, and $U(p,q)$ is identified to a subgroup of $O(2p, 2q)$. Thus we are led to consider the action of $U(p,q)$ on the spaces $S^{a\pm}(X^0)$ defined by restriction of the action of $O(2p, 2q)$. These spaces, however, can clearly never be irreducible for $U(p,q)$, because $U(p,q)$ commutes not simply with $\mathbb{R}^\times$, the real scalar dilation operators used to define the $S^{a\pm}(X^0)$, but with $\mathbb{C}^\times$, the complex scalar dilations. (Indeed, $U(p,q)$ may be characterized as the subgroup of $O(2p, 2q)$, which commutes with $\mathbb{C}^\times$.) Thus it is appropriate to break up the spaces $S^a(X^0)$ into subspaces that are homogeneous for the action of $\mathbb{C}^\times$ (i.e., eigenspaces for $\mathbb{C}^\times$). We observe

$$(4.3) \qquad \mathbb{C}^\times \simeq \mathbb{R}_+^\times \cdot \mathbb{T} = \{ |t| \cdot t/|t| \,\big|\, t \in \mathbb{C}^\times \}.$$

Here $\mathbb{T} = \{ s \in \mathbb{C}^\times \mid |s| = 1 \}$ is the usual unit circle. Thus a character $\psi$ of $\mathbb{C}^\times$ has the form

$$(4.4) \qquad \psi(t) = |t|^a (t/|t|)^n, \qquad a \in \mathbb{C}^\times, \ n \in \mathbb{Z}.$$

Since $|t| = (t\bar{t})^{1/2}$ where $\bar{t}$ is the complex conjugate of $t$, we can formally write

$$(4.5) \qquad \psi(t) = \psi^{\alpha,\beta}(t) = t^\alpha \bar{t}^\beta,$$

with $\alpha = (a+n)/2$ and $\beta = (a-n)/2$ or $a = \alpha + \beta$, and $n = \alpha - \beta$. We see that a pair $(\alpha, \beta)$ can occur in $(4.5)$ if and only if $\alpha - \beta$ is an integer. For all such pairs, we define $S^{\alpha,\beta}(X^0) \subseteq S^{\alpha+\beta}(X^0)$ to be the $\psi^{\alpha,\beta}$ eigenspace for $\mathbb{C}^\times$. We clearly have a decomposition

$$(4.6) \qquad S^a(X^0) = \sum_{\substack{\alpha+\beta=a \\ \alpha-\beta\in\mathbb{Z}}} S^{\alpha,\beta}(X^0).$$

Further each summand in equation $(4.6)$ is obviously invariant under the group $U(p,q)$. We will study the representations of $U(p,q)$ defined by the $S^{\alpha,\beta}(X^0)$. We remark that they can be interpreted as induced representations, again "degenerate principal series", just as the $S^{a\pm}(X^0)$ were interpreted as induced representations of $O(2p, 2q)$ (see formulas $(2.9)$–$(2.11)$). When $n = 0$ (or when $\alpha = \beta$ and only then), the representations $S^{\alpha,\beta}(X^0)$ will be spherical, that is, will contain a vector invariant under the maximal compact subgroup $U(p) \times U(q) \subseteq U(p,q)$. When $q = 1$ and $\alpha = \beta$, we have the spherical principal series of $U(p,1)$.

To understand the $U(p,q)$-module structure of $S^{\alpha,\beta}(X^0)$, we proceed as we did for the $S^a(X^0)$: First we describe the action of the maximal compact subgroup $U(p) \times U(q)$ and then we see how the noncompact part of the Lie algebra $\mathfrak{u}(p,q)$ of $U(p,q)$ moves around the $U(p) \times U(q)$-isotypic subspaces.

**4.2. $K$-Structure of $S^{\alpha,\beta}(X^0)$.** Under the inclusion of $\mathrm{U}(p,q)$ in $\mathrm{O}(2p,2q)$, it is clear that $\mathrm{U}(p) \times \mathrm{U}(q) \subseteq \mathrm{O}(2p) \times \mathrm{O}(2q)$. Since the $\mathrm{O}(2p) \times \mathrm{O}(2q)$-module structure of the $S^a(X^0)$ is described in terms of spherical harmonics (see formulas (2.12)–(2.16)), we are led to consider the structure of the space $\mathcal{H}^m(\mathbb{C}^p) \simeq \mathcal{H}^m(\mathbb{R}^{2p})$ as a $\mathrm{U}(p)$ module.

Consider the algebra $\mathcal{P}_{\mathbb{R}}(\mathbb{C}^p) \simeq \mathcal{P}(\mathbb{R}^{2p})$ of polynomials on $\mathbb{C}^p$ considered as a real vector space. We may choose the complex coordinates $z_1, z_2, \ldots, z_p$ and their complex conjugates $\overline{z_1}, \overline{z_2}, \ldots, \overline{z_p}$ as generators for $\mathcal{P}_{\mathbb{R}}(\mathbb{C}^p)$. For integers $\alpha, \beta$, let $\mathcal{P}^{\alpha,\beta}(\mathbb{C}^p)$ be the space of polynomials that are homogeneous of degree $\alpha$ in the $z_j$ and homogeneous of degree $\beta$ in the $\overline{z_j}$. This notation is consistent with that of §4.1, in the sense that $\mathcal{P}^{\alpha,\beta}(\mathbb{C}^p)$ is the $\psi^{\alpha,\beta}$ eigenspace for the action of $\mathbb{C}^\times$ on $\mathcal{P}_{\mathbb{R}}(\mathbb{C}^p)$ by scalar dilations. Clearly we have a decomposition

$$(4.7) \qquad\qquad \mathcal{P}_{\mathbb{R}}^m(\mathbb{C}^p) \simeq \sum_{\alpha+\beta=m} \mathcal{P}^{\alpha,\beta}(\mathbb{C}^p)$$

of $\mathbb{R}$-homogeneous polynomials of degree $m$ into bihomogeneous polynomials of appropriate bi-degrees.

The Laplacian on $\mathbb{C}^p = \mathbb{R}^{2p}$ can be written

$$\Delta = \sum_{j=1}^{2p} \frac{\partial^2}{\partial x_j^2} = 4 \sum_{j=1}^{p} \frac{\partial^2}{\partial z_j \partial \overline{z_j}}.$$

From this formula it is clear that $\Delta$ maps $\mathcal{P}^{\alpha,\beta}(\mathbb{C}^p)$ to $\mathcal{P}^{\alpha-1,\beta-1}(\mathbb{C}^p)$. It follows that the harmonic polynomials have a decomposition analogous to (4.7), namely,

$$(4.8) \qquad \mathcal{H}^m(\mathbb{C}^p) = \sum_{\alpha+\beta=m} \mathcal{H}^m(\mathbb{C}^p) \cap \mathcal{P}^{\alpha,\beta}(\mathbb{C}^p) = \sum_{\alpha+\beta=m} \mathcal{H}^{\alpha,\beta}(\mathbb{C}^p),$$

where the second equation serves to define the $\mathcal{H}^{\alpha,\beta}(\mathbb{C}^p)$. Since $\mathrm{U}(p) \subseteq \mathrm{O}(2p)$ commutes with $\mathbb{C}^\times$, the spaces $\mathcal{H}^{\alpha,\beta}(\mathbb{C}^p)$ will be invariant under the action of $\mathrm{U}(p)$ on $\mathcal{H}^m(\mathbb{C}^p)$. Hence they define representations of $\mathrm{U}(p)$. It is known (see [Zh]; this is also a special case of the general results of [Ho3]) that the $\mathcal{H}^{\alpha,\beta}(\mathbb{C}^p)$ are irreducible and mutually inequivalent $\mathrm{U}(p)$ modules. We note that, in particular, the spaces $\mathcal{H}^{\alpha,\beta}(\mathbb{C}^p)$ are eigenspaces for the group $\mathbb{C}^\times \cap \mathrm{U}(n) \simeq \mathbb{T}$ of scalar unitary matrices; the eigencharacter by which $\mathbb{T}$ acts on $\mathcal{H}^{\alpha,\beta}(\mathbb{C}^p)$ is

$$(4.9) \qquad\qquad\qquad \chi_{\alpha-\beta} : z \to z^{\alpha-\beta}.$$

If $p \geq 2$, then $\mathcal{H}^{\alpha,\beta}(\mathbb{C}^p)$ is nonzero for all $\alpha, \beta \geq 0$; however, when $p = 1$, only $\mathcal{H}^{\alpha,0}(\mathbb{C}) = \mathbb{C}z^\alpha$ and $\mathcal{H}^{0,\beta}(\mathbb{C}) = \mathbb{C}\overline{z}^\beta$ are nonzero. This degeneracy in the $\mathcal{H}^{\alpha,\beta}(\mathbb{C}^p)$ when $p = 1$ leads to some interesting phenomena, which will be detailed later in the paper (see §4.5).

Next consider the action of $\mathrm{U}(p) \times \mathrm{U}(q)$ on the algebra $\mathcal{P}_{\mathbb{R}}(\mathbb{C}^{p+q})$. Taking the tensor product of equation (4.8) with its analog for $\mathbb{C}^q$ gives us a decomposition

$$(4.10) \qquad \mathcal{H}^m(\mathbb{C}^p) \otimes \mathcal{H}^n(\mathbb{C}^q) \simeq \sum_{\substack{m_1+m_2=m \\ n_1+n_2=n}} \mathcal{H}^{m_1,m_2}(\mathbb{C}^p) \otimes \mathcal{H}^{n_1,n_2}(\mathbb{C}^q).$$

The action of the scalar unitary matrices on the summand $\mathcal{H}^{m_1,m_2}(\mathbb{C}^p) \otimes \mathcal{H}^{n_1,n_2}(\mathbb{C}^q)$ is via the character $\chi_r$ where $r = m_1 - m_2 + n_1 - n_2$ [see definition (4.9)]. Using this, we can easily see that the maps $j_{a,m,n}$ of formula (2.15) split into maps

$$(4.11) \quad j_{\alpha,\beta,m_1,m_2,n_1,n_2} = j_{\alpha,\beta} \; : \; \mathcal{H}^{m_1,m_2}(\mathbb{C}^p) \otimes \mathcal{H}^{n_1,n_2}(\mathbb{C}^q) \to S^{\alpha,\beta}(X^0),$$

where

$$\alpha + \beta = a, \qquad \alpha - \beta = m_1 - m_2 + n_1 - n_2.$$

In this connection we observe that

$$(4.12) \qquad r_{2p}^2 = \sum_{j=1}^{2p} x_j^2 = \sum_{j=1}^{p} z_p \overline{z_p}$$

so that $r_{2p}^2$ is in $\mathcal{P}^{1,1}(\mathbb{C}^p)$ and multiplication by any power of $r_{2p}^2$ increases both components of bidegree equally and so does not affect their difference.

Combining formulas (2.16), (4.6), (4.10), and (4.11) yields the decomposition

$$
\begin{aligned}
(4.13) \qquad S^{\alpha,\beta}(X^0) &= \sum_{m,n \geq 0} \left( j_{\alpha+\beta}(\mathcal{H}^m(\mathbb{C}^p) \otimes \mathcal{H}^n(\mathbb{C}^q)) \cap S^{\alpha,\beta}(X^0) \right) \\
&= \sum_{m,n \geq 0} \sum_{\substack{m_1,m_2,n_1,n_2 \geq 0 \\ m_1+m_2=m \\ n_1+n_2=n \\ m_1-m_2+n_1-n_2=\alpha-\beta}} j_{\alpha,\beta}(\mathcal{H}^{m_1,m_2}(\mathbb{C}^p) \otimes \mathcal{H}^{n_1,n_2}(\mathbb{C}^q)) \\
&= \sum_{\substack{m_1,m_2,n_1,n_2 \geq 0 \\ m_1-m_2+n_1-n_2=\alpha-\beta}} j_{\alpha,\beta}(\mathcal{H}^{m_1,m_2}(\mathbb{C}^p) \otimes \mathcal{H}^{n_1,n_2}(\mathbb{C}^q))
\end{aligned}
$$

of $S^{\alpha,\beta}(X^0)$ into irreducible $\mathrm{U}(p) \times \mathrm{U}(q)$ modules.

From the final form of the decomposition, we see that the $\mathrm{U}(p) \times \mathrm{U}(q)$ components of $S^{\alpha,\beta}(X^0)$ form a three-parameter slice out of the full four-parameter family of $\mathrm{U}(p) \times \mathrm{U}(q)$ representations that appear in $\mathcal{P}_{\mathbb{R}}(\mathbb{C}^{p+q})$. We may refer to these irreducible $\mathrm{U}(p) \times \mathrm{U}(q)$ modules as "$K$-types". The middle form of the decomposition emphasizes the relation of these $\mathrm{U}(p) \times \mathrm{U}(q)$ modules and the $\mathrm{O}(2p) \times \mathrm{O}(2q)$ modules of the first form of the decomposition. We see from formula (4.10) that the $\mathrm{U}(p) \times \mathrm{U}(q)$ submodules in $\mathcal{H}^m(\mathbb{C}^p) \otimes \mathcal{H}^n(\mathbb{C}^q)$ can be parametrized by the pair $(m_1, n_1)$, which varies in the integral points in a rectangle (see Diagram 4.14).

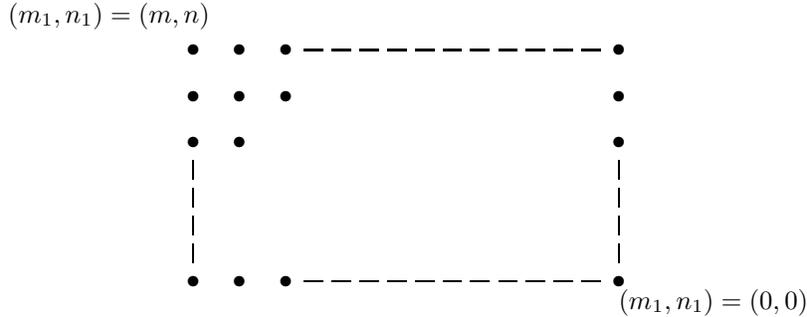

DIAGRAM 4.14

The subset of points such that $m_1 - m_2 + n_1 - n_2 = \alpha - \beta$ forms a diagonal line inside this rectangle (see Diagram 4.15).

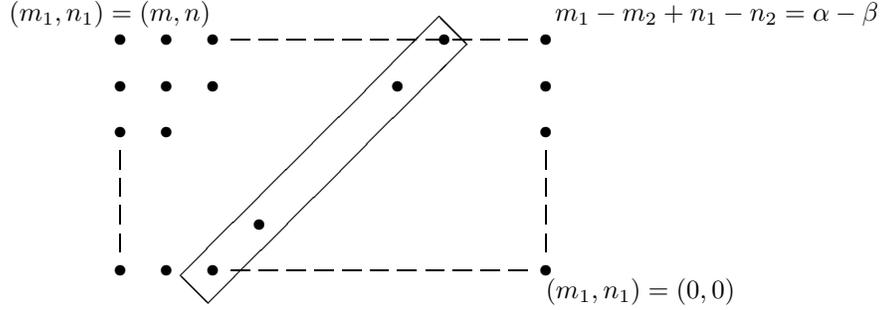

**Diagram 4.15**

Here we should observe that since

$$m_1 - m_2 + n_1 - n_2 = m_1 + m_2 + n_1 + n_2 - 2(m_2 + n_2)$$
$$= m + n - 2(m_2 + n_2),$$

the inner sum in the middle line of formula (4.13) will be empty unless

(i) $m + n \geq |\alpha - \beta|$,
(ii) $m + n \equiv \alpha - \beta \mod 2$.

When $q = 1$, there are further restrictions, described in §4.5. For now, we assume $p, q > 1$.

To exploit the relation between the situation under study and the results of §§2 and 3, we will represent formulas (4.13) graphically as follows. We think of the points $(m, n)$ for which the inner sum in the middle line of (4.13) is non-empty as being integral points in the plane. This gives us all points in a coset of the "diamond lattice" in the positive quadrant and lying above the line $x + y = |\alpha - \beta|$, as illustrated by Diagram 4.16.

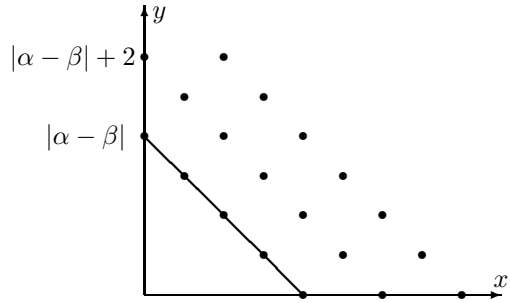

**Diagram 4.16**

We further think of the "line" of $\mathrm{U}(p) \times \mathrm{U}(q)$ modules in the inner sum of the middle line of (4.13) as points in a "fiber" lying over the point $(m, n)$ in Diagram

4.16. We call the region

$$\{(x,y) \mid x \geq 0, \ y \geq 0, \ x + y \geq |\alpha - \beta|\}$$

the *K-type region*.

**4.3. Action of the noncompact part of $\mathfrak{u}(p,q)$ on $S^{\alpha,\beta}(X^0)$.** Consider now how the noncompact part of the Lie algebra $\mathfrak{u}(p,q)$ acts on the $\mathrm{U}(p) \times \mathrm{U}(q)$ components, as described in the preceding section. Since $\mathfrak{u}(p,q) \subseteq \mathfrak{o}(2p,2q)$, it must act compatibly with Diagram 2.19, so the $\mathrm{U}(p) \times \mathrm{U}(q)$ modules in the fiber over $(m,n)$ must be taken to the fibers over $(m \pm 1, \ n \pm 1)$; but since each fiber consists of several $\mathrm{U}(p) \times \mathrm{U}(q)$ modules, we must calculate how the $\mathfrak{u}(p,q)$ action affects the various points in each fiber. To formulate this we observe that the complexified Lie algebra $\mathfrak{u}(p,q)_{\mathbb{C}}$ has a decomposition

$$(4.17) \qquad \mathfrak{u}(p,q)_{\mathbb{C}} \simeq (\mathfrak{u}(p) \oplus \mathfrak{u}(q))_{\mathbb{C}} \oplus \mathfrak{p}^+ \oplus \mathfrak{p}^-,$$

where $\mathfrak{p}^{\pm}$ are the eigenspaces for the adjoint action of the center of $\mathfrak{u}(p) \oplus \mathfrak{u}(q)$. Concretely we have

$$(4.18) \qquad \begin{aligned} \mathfrak{p}^+ &= \mathrm{span} \ \left\{ S_{ij} = z_i \frac{\partial}{\partial w_j} + \overline{w_j} \frac{\partial}{\partial \overline{z_i}} \right\}, \\ \mathfrak{p}^- &= \mathrm{span} \ \left\{ T_{ij} = w_j \frac{\partial}{\partial z_i} + \overline{z_i} \frac{\partial}{\partial \overline{w_j}} \right\}. \end{aligned}$$

The following lemma is a refinement to the current situation for Lemma 2.3. We leave the calculations justifying it to §6.2.

**Lemma 4.1.** *There exist nonzero maps*

$$\begin{aligned} &S^{1001}_{m_1,m_2,n_1,n_2}, \quad S^{0-101}_{m_1,m_2,n_1,n_2}, \quad S^{10-10}_{m_1,m_2,n_1,n_2}, \quad S^{0-1-10}_{m_1,m_2,n_1,n_2}, \\ &T^{0110}_{m_1,m_2,n_1,n_2}, \quad T^{-1010}_{m_1,m_2,n_1,n_2}, \quad T^{010-1}_{m_1,m_2,n_1,n_2}, \quad T^{-100-1}_{m_1,m_2,n_1,n_2}, \end{aligned}$$

*where*

$$\begin{aligned} S^{\varepsilon_1 \varepsilon_2 \eta_1 \eta_2}_{m_1,m_2,n_1,n_2} &: \mathfrak{p}^+ \otimes (\mathcal{H}^{m_1,m_2}(\mathbb{C}^p) \otimes \mathcal{H}^{n_1,n_2}(\mathbb{C}^q)) \\ &\qquad \to \mathcal{H}^{m_1+\varepsilon_1,m_2+\varepsilon_2}(\mathbb{C}^p) \otimes \mathcal{H}^{n_1+\eta_1,n_2+\eta_2}(\mathbb{C}^q), \\ T^{\varepsilon_1 \varepsilon_2 \eta_1 \eta_2}_{m_1,m_2,n_1,n_2} &: \mathfrak{p}^- \otimes (\mathcal{H}^{m_1,m_2}(\mathbb{C}^p) \otimes \mathcal{H}^{n_1,n_2}(\mathbb{C}^q)) \\ &\qquad \to \mathcal{H}^{m_1+\varepsilon_1,m_2+\varepsilon_2}(\mathbb{C}^p) \otimes \mathcal{H}^{n_1+\eta_1,n_2+\eta_2}(\mathbb{C}^q), \end{aligned}$$

*independent of $\alpha, \beta$, such that the action of $z \in \mathfrak{p}^+$ and $\tilde{z} \in \mathfrak{p}^-$ [see (4.18)] on the K-type $j_{\alpha,\beta}(\mathcal{H}^{m_1,m_2}(\mathbb{C}^p) \otimes \mathcal{H}^{n_1,n_2}(\mathbb{C}^q))$ is described by the formula*

$$\begin{aligned} \rho(z &+ \tilde{z}) j_{\alpha,\beta}(\phi) \\ &= A^{++}(m,n) j_{\alpha,\beta} \left\{ S^{1001}_{m_1,m_2,n_1,n_2}(z \otimes \phi) + T^{0110}_{m_1,m_2,n_1,n_2}(\tilde{z} \otimes \phi) \right\} \\ &\quad + A^{-+}(m,n) j_{\alpha,\beta} \left\{ S^{0-101}_{m_1,m_2,n_1,n_2}(z \otimes \phi) + T^{-1010}_{m_1,m_2,n_1,n_2}(\tilde{z} \otimes \phi) \right\} \\ &\quad + A^{+-}(m,n) j_{\alpha,\beta} \left\{ S^{10-10}_{m_1,m_2,n_1,n_2}(z \otimes \phi) + T^{010-1}_{m_1,m_2,n_1,n_2}(\tilde{z} \otimes \phi) \right\} \\ &\quad + A^{--}(m,n) j_{\alpha,\beta} \left\{ S^{0-1-10}_{m_1,m_2,n_1,n_2}(z \otimes \phi) + T^{-100-1}_{m_1,m_2,n_1,n_2}(\tilde{z} \otimes \phi) \right\}, \end{aligned}$$

*where $m = m_1 + m_2$, $n = n_1 + n_2$, $\phi \in \mathcal{H}^{m_1, m_2}(\mathbb{C}^p) \otimes \mathcal{H}^{n_1, n_2}(\mathbb{C}^q)$, and $A^{\pm\pm} = A^{\pm\pm}_{2p, 2q, \alpha+\beta}$ is as defined in (2.20).*

*Proof.* See §6.2. □

The action described in this lemma can be pictured as follows. Imagine the fibers over $(m, n)$ and $(m + 1, n + 1)$ as lines of points, as suggested by Diagram 4.15. Then the effect of the $S$ and $T$ operators of Lemma 4.1 can be illustrated as in Diagram 4.19.

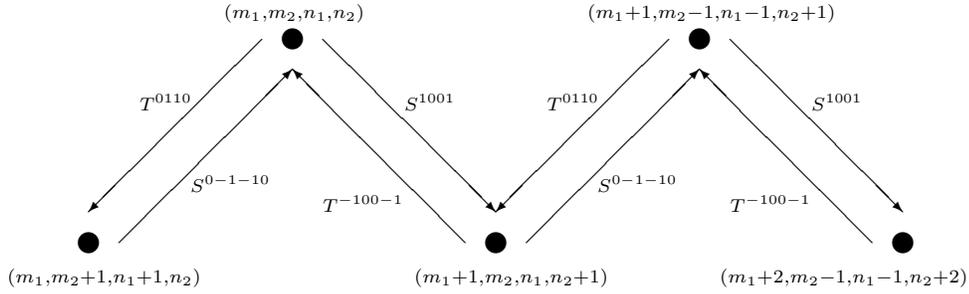

DIAGRAM 4.19

From Diagram 4.19, we see that, providing the transition coefficients $A^{++}$, $A^{--}$ are nonzero at $(m, n)$ and $(m + 1, \ n + 1)$ respectively, that is, providing that the fibers over $(m, n)$ and $(m + 1, \ n + 1)$ belong to the same $\mathfrak{o}(2p, 2q)$ constituent of $S^{\alpha+\beta}(X^0)$, application of an element of $\mathfrak{p}^+$ followed by another from $\mathfrak{p}^+$ moves a given $U(p) \times U(q)$ module to its neighbor in the same fiber. Computations in §6.2 will show that by successive applications of $(\mathfrak{p}^+)^2$ [or $(\mathfrak{p}^-)^2$] we can move from any point in a fiber to any other. A similar analysis applies to transitions between $(m, n)$ and $(m + 1, \ n - 1)$. This implies the result: The $U(p, q)$-module structure of $S^{\alpha, \beta}(X^0)$ is controlled by the $O(2p, 2q)$-module structure.

**Theorem 4.2.** *Assume $p, q \geq 2$. The image of $S^{\alpha, \beta}(X^0)$ in any $O(2p, 2q)$-irreducible constituent of $S^{\alpha+\beta}(X^0)$ is an irreducible $U(p, q)$ module.*

In other words, the extra degree of freedom in the $U(p) \times U(q)$ fibers has no effect whatsoever on the $U(p, q)$ composition series of the $S^{\alpha, \beta}(X^0)$, except in so far as some fibers may be empty. We may therefore read off the $U(p, q)$ constituents of the $S^{\alpha, \beta}(X^0)$ from the corresponding results from $O(2p, 2q)$ (see Diagrams 2.25–28, 2.34, and 2.35), providing we take into account Diagram 4.16 describing the nontrivial fibers for $S^{\alpha, \beta}(X^0)$.

**Corollary 4.3.** *$S^{\alpha, \beta}(X^0)$ is reducible if and only if $\alpha, \beta \in \mathbb{Z}$, that is, if and only if $\alpha + \beta \in \mathbb{Z}$ and $\alpha + \beta \equiv \alpha - \beta \mod 2$.*

*Proof.* For reducibility, we need one of the transition coefficients (2.20) to vanish for some $(m, n)$ corresponding to a nonempty fiber of $\mathrm{U}(p) \times \mathrm{U}(q)$ modules. For example, the coefficient $A^{+-}_{2p,2q,\alpha+\beta}(m, n)$ [see definition (2.20)] will vanish if and only if

$$\alpha + \beta = m - n - 2(q - 1).$$

Since we know $\alpha - \beta \in \mathbb{Z}$, and $m - n \equiv \alpha - \beta \mod 2$, we see that $\alpha + \beta \in \mathbb{Z}$ and $\alpha + \beta \equiv \alpha - \beta \mod 2$. This implies $\alpha, \beta \in \mathbb{Z}$. Conversely, if $\alpha, \beta \in \mathbb{Z}$, then we can reverse the process to find $(m, n)$ such that $A^{+-}(m, n) = 0$, implying reducibility. The same conclusion follows by considering $A^{-+}$, and $A^{++}$ and $A^{--}$ also cannot cause reducibility unless $\alpha, \beta \in \mathbb{Z}$. $\square$

To describe what happens when $\alpha$ and $\beta$ are integers, we simply have to combine the relevant diagrams from §2 with Diagram 4.16. The various possibilities are described in the following series of diagrams.

(i) $\alpha, \beta \in \mathbb{Z}_+$. When $\alpha, \beta$ are both nonnegative integers, the barriers $A^{+-}$, $A^{-+}$, and $A^{++}$ all intersect the $K$-type region and cause $S^{\alpha,\beta}(X^0)$ to break into four parts: one finite-dimensional subrepresentation, two quotients, and one constituent that is neither a quotient nor a subrepresentation (see Diagram 4.20).

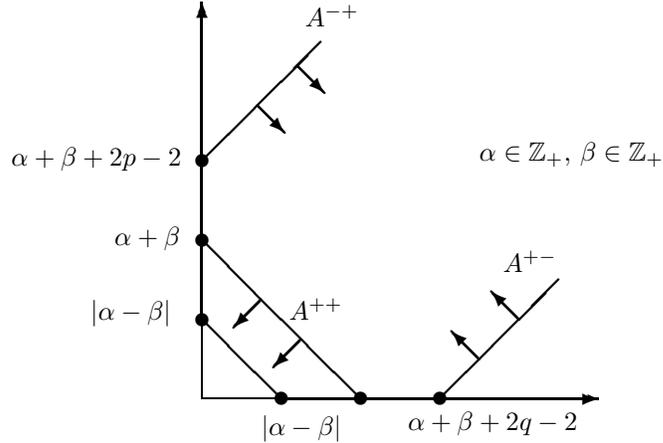

DIAGRAM 4.20. $\alpha \in \mathbb{Z}_+$, $\beta \in \mathbb{Z}_+$.

(ii) $\alpha$ or $\beta$ is negative, but $\alpha + \beta \geq 1 - p - q$. Here the $A^{++}$ barrier no longer intersects the region of $K$-types, and there are at most three constituents, all infinite dimensional. As $\alpha + \beta$ decreases while $\alpha - \beta$ remains fixed, the barriers $A^{+-}$ and $A^{-+}$ approach each other. When $\alpha$ and $\beta$ are nonnegative, the point where the $A^{+-}$ barrier leaves the $K$-type region, that is, the intersection of the $A^{+-}$ barrier with the boundary of the $K$-type region, is on the $x$-axis. Similarly, the barrier $A^{-+}$ passes through the boundary on the $y$-axis. As $\alpha + \beta$ decreases, $A^{+-}$ moves left and $A^{-+}$ moves right. When $\alpha + \beta$ becomes sufficiently small (precisely, when $\min(\alpha, \beta) \leq 1 - q$), the barrier $A^{+-}$ leaves the $K$-type region through the antidiagonal segment $\{x + y = |\alpha - \beta| \mid x \geq 0, \ y \geq 0\}$. As $\alpha + \beta$ continues to decrease, the point of exit of $A^{+-}$ moves up the antidiagonal segment and finally (when $\max(\alpha, \beta) \leq 1-q$) moves onto the $y$-axis. The barrier $A^{-+}$ moves in the opposite direction. Its point of exit crosses from the $y$-axis to the antidiagonal segment when $\min(\alpha, \beta) = 1-p$ and further crosses from the antidiagonal to the $x$-axis when $\max(\alpha, \beta) = 1-p$. The two barriers meet when $\alpha + \beta = 2-p-q$; their common point of exit from the $K$-type region is on the antidiagonal segment if and only if $|p-q| \leq |\alpha - \beta|$. The real point on the unitary axis (see §4.4) is at $\alpha + \beta = 1-p-q$. For a given value of $\alpha - \beta$, exactly one of the values $2-p-q$ and $1-p-q$ for $\alpha+\beta$ will be a point of reducibility. If the unitary point $\alpha+\beta = 1-p-q$ is reducible, then $S^{\alpha,\beta}(X^0)$ breaks into a sum of two irreducible representations, as illustrated in Diagram 2.26. An illustrative situation among the several possibilities detailed above is pictured in Diagram 4.21.

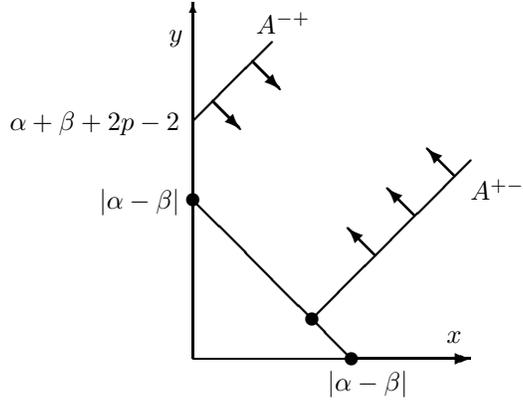

Diagram 4.21. $1-p \leq \min(\alpha, \beta) \leq 1-q$.

(iii) $\alpha + \beta < 1-p-q$. In analogy with the remark following Lemma 3.1, the modules $S^{\alpha,\beta}(X^0)$ and $S^{\alpha',\beta'}(X^0)$, where $\alpha+\alpha' = 1-p-q = \beta+\beta'$, are naturally dual, so the structure of $S^{\alpha,\beta}(X^0)$ when $\alpha + \beta < 1-p-q$ is deducible from that of $S^{\alpha',\beta'}(X^0)$ as described above.

**4.4. Unitarity.** The unitarity criteria (3.6) for $O(2p, 2q)$ apply directly here and the results are as summarized in Diagrams 3.7 and 3.13. The unitary axis is defined by $\mathrm{Re}(\alpha+\beta) = 1-p-q$. The possibilities for the values of $\alpha+\beta$ where the full module $S^{\alpha,\beta}(X^0)$ is unitary are as illustrated in Diagram 3.8 (note that $\tilde{a} = \alpha+\beta+p+q-1$, $\varepsilon_0$ is $+$ if $p-q$ is odd, and $\varepsilon_1$ is $-$ if $p-q$ is even). At a typical point of reducibility, as in Diagrams 4.20 and 4.21, the one or two constituents in the central strip are not unitary, while the constituents outside the barriers are unitary. Exceptions in the range $\mathrm{Re}(\alpha+\beta) \geq 1-p-q$ occur when $\alpha = \beta = 0$, when the finite-dimensional subrepresentation is the trivial representation, and when $\alpha + \beta = 2-p-q$, with $\alpha$, $\beta \in \mathbb{Z}$, in which case all three components, one of which has $K$-types only along the line $y-x = p-q$, are unitary. See Diagram 3.9(a).

**4.5. The case of** $\mathrm{U}(p,1)$, $p > 1$. For $\mathrm{U}(p,1)$, degeneration of the $K$-types leads to a situation substantially more complex than what is pictured in Diagram 4.16.

Since the module structure and unitarity result from an interaction between the $K$-types and the barriers, this change in $K$-type structure causes a corresponding change in the module structure and unitarity.

As noted in the discussion following formula (4.9), the only $\mathcal{H}^{a,b}(\mathbb{C})$ that are nonzero are $\mathcal{H}^{n,0}(\mathbb{C})$ and $\mathcal{H}^{0,n}(\mathbb{C})$, $n \in \mathbb{Z}_+$. Hence, instead of what is represented in Diagram 4.14 as a rectangle of $K$-types, we have only the $K$-types corresponding to the right and left sides of this rectangle, and in a given space $S^{\alpha,\beta}(X^0) \cap j_{\alpha+\beta}(\mathcal{H}^m(\mathbb{C}^p) \otimes \mathcal{H}^n(\mathbb{C}))$, we will have at most two irreducible components for $\mathrm{U}(p) \times \mathrm{U}(1)$.

This leads us to the following construction. Let us write

$$
(4.22) \qquad \mathcal{H}^n(\mathbb{C}) = \begin{cases} \mathcal{H}^{n,0}(\mathbb{C}), & n \geq 0, \\ \mathcal{H}^{0,-n}(\mathbb{C}), & n \leq 0. \end{cases}
$$

Further let us represent $j_{\alpha+\beta}\left(\mathcal{H}^{m_1,m_2}(\mathbb{C}^p) \otimes \mathcal{H}^n(\mathbb{C})\right)$ by the point $(m_1 + m_2, n)$ in $\mathbb{R}^2$. Here $m_1, m_2 \in \mathbb{Z}_+$ and $n \in \mathbb{Z}$. Since

$$
j_{\alpha+\beta}\left(\mathcal{H}^{m_1,m_2}(\mathbb{C}^p) \otimes \mathcal{H}^n(\mathbb{C})\right) \subseteq S^{\alpha,\beta}(X^0),
$$

where

$$
(4.23) \qquad \alpha - \beta = m_1 - m_2 + n,
$$

and the triple $(m_1, m_2, n)$ is determined by the pair $(m_1 + m_2, n)$ together with $\alpha - \beta$, we see that each point $(m_1 + m_2, \ n)$ represents at most one $K$-type in $S^{\alpha,\beta}(X^0)$. Furthermore, if the point $(m,n)$, $m \in \mathbb{Z}_+$, $n \in \mathbb{Z}$, does represent a $K$-type of $\mathrm{O}(2p, 2q)$, then we must have

$$
(4.24) \qquad \begin{aligned}
& m_1 = \tfrac{1}{2}(m - (n - (\alpha - \beta))) \geq 0, \\
& m_2 = \tfrac{1}{2}(m + (n - (\alpha - \beta))) \geq 0, \ \text{or} \\
& m \geq |n - (\alpha - \beta)|.
\end{aligned}
$$

Conversely, as long as $p > 1$, this condition plus $m + n \equiv \alpha - \beta \mod 2$ guarantees we can find $m_1$, $m_2$, and $n$ so that $j_{\alpha+\beta}\left(\mathcal{H}^{m_1,m_2}(\mathbb{C}^p) \otimes \mathcal{H}^n(\mathbb{C})\right) \subseteq S^{\alpha,\beta}(X^0)$ sits over $(m,n)$. Thus the $K$-type region for $S^{\alpha,\beta}(X^0)$ is as depicted in Diagram 4.25.

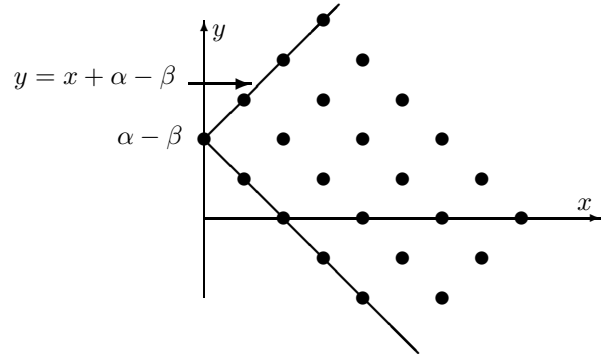

<div align="center">Diagram 4.25</div>

Note that the $x$-intercept is no longer $|\alpha - \beta|$ but now is just $\alpha - \beta$.

When we add barriers to this picture to determine reducibility, we find that each of the former barriers $A^{\pm,\pm}$ produces two barriers, one in its original position and a second reflected through the $x$-axis. Since when $q = 1$ the barriers $A^{++}$ and $A^{+-}$ intersect on the $x$-axis, the reflection of one continues the other in a straight line. The reflection of $A^{-+}$ is $A^{--}$ and vice versa. Thus, with barriers inserted, Diagram 4.25 is transformed into Diagram 4.26 (the shaded area is the $K$-type region).

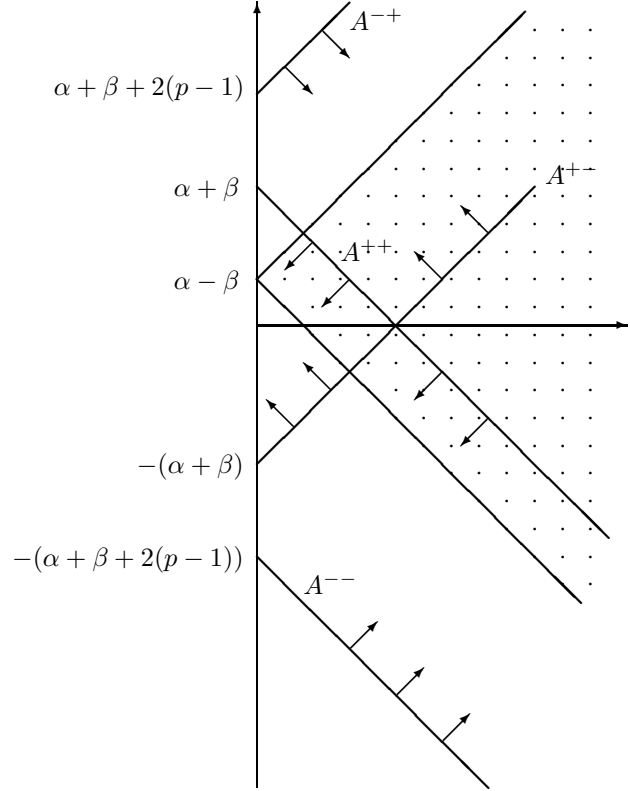

DIAGRAM 4.26

Study of Diagram 4.26 and similar ones yields the following conclusions when $\alpha$ and $\beta$ are integers.

**Lemma 4.4.** *Suppose* $\alpha$, $\beta \in \mathbb{Z}$. *Consider the structure of* $S^{\alpha,\beta}(X^0)$.

(a) *The barrier* $A^{++}$ *causes reducibility when* $\beta \geq 0$.
(b) *The barrier* $A^{+-}$ *causes reducibility when* $\alpha \geq 0$.
(c) *The barrier* $A^{-+}$ *causes reducibility when* $\beta \leq -p$.
(d) *The barrier* $A^{--}$ *causes reducibility when* $\alpha \leq -p$.

Given this information, the simplest way to organize the description of the structure of $S^{\alpha,\beta}(X^0)$ is by means of a picture describing the Hasse diagram of $S^{\alpha,\beta}(X^0)$ in various regions of the $(\alpha, \beta)$-plane (see Diagram 4.27). The diagram is symmetric across the line $\alpha + \beta = -p$ because of the natural duality between $S^{\alpha,\beta}(X^0)$ and $S^{-p-\alpha,-p-\beta}(X^0)$.

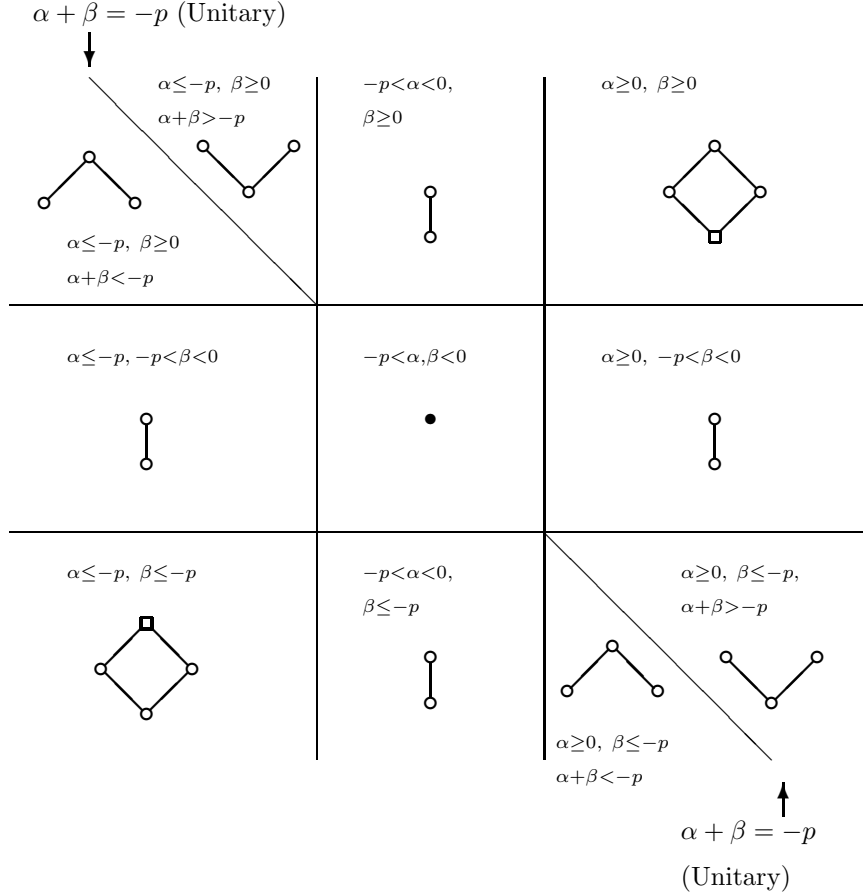

Legend:   • irreducibility and unitary
          □ finite-dimensional constituent

DIAGRAM 4.27. Hasse diagrams for $S^{\alpha,\beta}(X^0)$, $\alpha$, $\beta$ in $\mathbb{Z}$.

Unitarity is determined again by the conditions (3.4). The main novel feature of the case at hand is that some of the conditions are now vacuous because there are no $K$-types corresponding to transitions that would prevent unitarity. Thus the unitary axis occurs at $\alpha + \beta = -p$. If also $|\alpha - \beta| < p$, that is, if $\alpha$ and $\beta$ are in the central square of Diagram 4.27, then all barriers miss the $K$-type region, and if $\alpha + \beta$ varies while $\alpha - \beta$ is fixed, the barriers continue to miss the whole $K$-type region in the interval $-2p + |\alpha - \beta| < \alpha + \beta < -|\alpha - \beta|$.

Consequently, the representation remains irreducible and unitary for $\alpha + \beta$ in this range. In other words, when $|\alpha - \beta| < p$, the unitary point $\alpha + \beta = -p$ is irreducible regardless of the parity of $|\alpha - \beta|$ and the complementary series extends to a distance $p - |\alpha - \beta|$ on either side of the unitary axis:

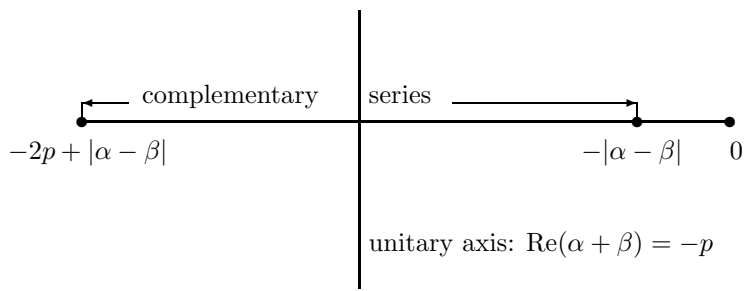

DIAGRAM 4.28. Complementary series for $U(p,1)$

In particular when $\alpha - \beta = 0$ [which is the case of the spherical principal series of $U(p,1)$], the complementary series extends all the way from the unitary axis $(\alpha + \beta = -p)$ to $\alpha + \beta = 0$, that is, $\alpha = \beta = 0$, which contains the trivial representation as its irreducible subrepresentation. This has been known since Kostant [Ko].

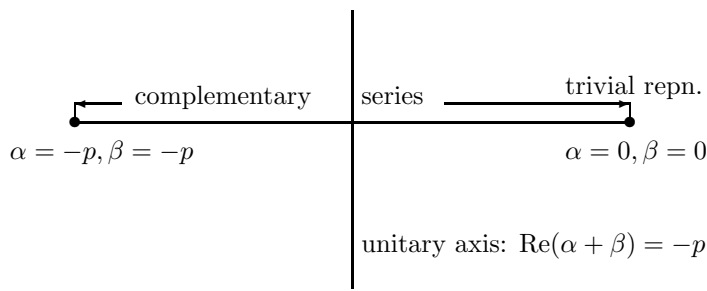

DIAGRAM 4.29. Complementary series for the spherical principal series of $U(p,1)$

At points of reducibility, some of the constituents may be unitary. When $\alpha, \beta \in \mathbb{Z}_+$, the "large" constituent that occupies the main part of the interior of the $K$-type region is unitary and the constituents living along the boundary of the $K$-type region are nonunitary (see Diagram 4.30).

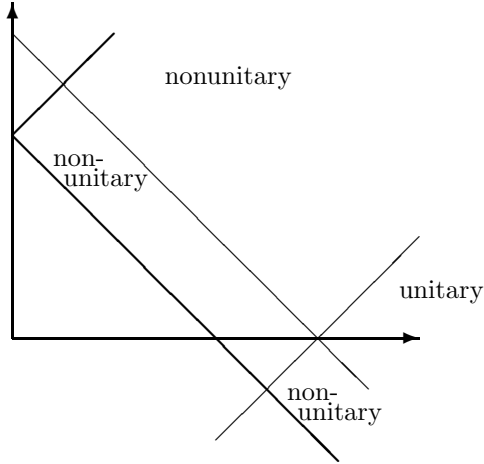

DIAGRAM 4.30

Exceptions arise, however, when a constituent consists only of $K$-types on the boundary of the $K$-type region. In this case there are no transitions in the direction transverse to the boundary, so the corresponding obstructions to unitarity do not apply and the corresponding module is unitary. In Diagram 4.31 we describe this phenomenon for $\alpha, \beta \geq 0$.

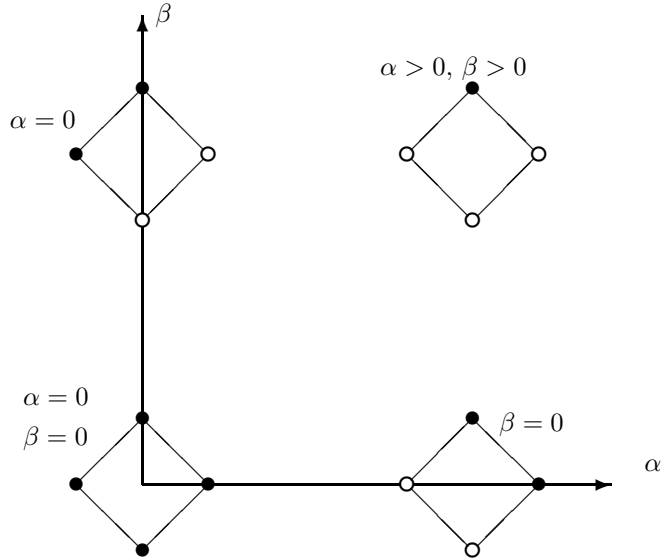

DIAGRAM 4.31. The filled-in circles indicate the unitary constituents.

Similar remarks apply to the other regions of $(\alpha, \beta)$ depicted in Diagram 4.27.

*Concluding remarks.* (a) The unitary representations with $K$-types along a boundary line of the $K$-type region have received considerable attention in the literature and are frequently called "ladder representations" [Ab, AFR, DGN, En, FR, SW].

(b) Consider the region of the $(\alpha, \beta)$-plane that is shaded in Diagram 4.32.

For a given value of $\alpha - \beta \in \mathbb{Z}$, the complementary series consists of precisely those values of $(\alpha, \beta)$ that fall inside the shaded region. It is well known that the fundamental group of $\mathrm{SU}(p, 1)$ is $\mathbb{Z}$. By considering representations of its universal covering group, one can get a family of representations parametrized by $(\alpha, \beta)$ in which $\alpha - \beta$ can be an arbitrary real number. It is perhaps plausible that the set of $(\alpha, \beta)$ for which these representations are unitary is exactly the shaded region (excluding its boundary). A. Koranyi (oral communication) has verified this for $\mathrm{SU}(2, 1)$.

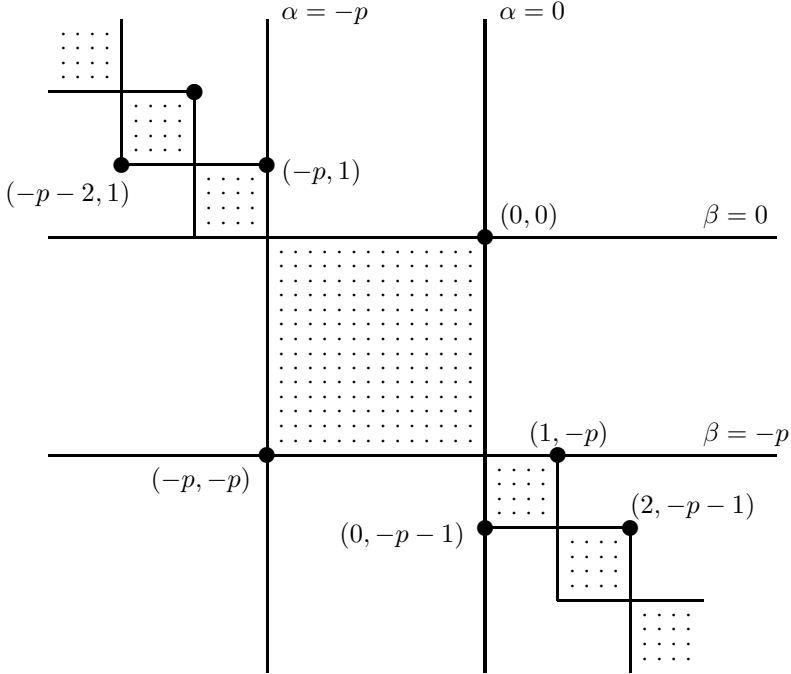

DIAGRAM 4.32

(c) When $|\alpha - \beta| \geq p$, the $K$-type region contains a parallel pair of barriers (either $A^{+-}$, $A^{-+}$ or $A^{++}$, $A^{--}$) when $\alpha + \beta = -p$ is the real point on the unitary axis and the complementary series is restricted to the same region as when $p, q > 1$ (see Diagram 3.8). In this situation, however, we have the possibility that $\beta \leq -p$, $\alpha \geq 0$ (or the reverse) so that at a point of reducibility, both the barriers $A^{+-}$ and $A^{-+}$ (or both of $A^{++}$ and $A^{--}$) intersect the $K$-type region and create submodules. These are the upper left and lower right regions in Diagram 4.27. When this happens, the constituents outside the barriers are unitary, and one of them (for example, the region above the barrier $A^{-+}$ when $\beta \leq -p \leq \alpha + \beta$) will have $K$-types contained in a diagonal strip of finite width (above $A^{-+}$ and below the boundary of the $K$-type region). Thus we obtain for $\mathrm{U}(p,1)$ a family of "ladderlike" representations whose $K$-types live not on a single diagonal line but in a diagonal strip. Also we should note that when the real point on the unitary axis is reducible (i.e., when $\alpha - \beta \equiv 0$ mod 2), the two constituents are very different: one is ladderlike, with $K$-types in a finite strip, and the other has $K$-types ranging over a full quarter plane.

## 5. The $\mathrm{Sp}(p,q)$-module structure of homogeneous functions on light cones

**5.1. The $\mathrm{Sp}(p,q)$ modules $S_j^a(X^0)$.** Let $\mathbb{H}$ denote the algebra of quaternions, and let $\mathbb{H}^{p+q} \simeq \mathbb{H}^p \oplus \mathbb{H}^q$ denote the vector space over $\mathbb{H}$ consisting of tuples

$$\omega = (u,v) = (u_1, \ldots, u_p, v_1, \ldots, v_q), \quad u_l, v_l \in \mathbb{H}.$$

We consider $\mathbb{H}^{p+q}$ to be endowed with the indefinite Hermitian form

$$(5.1) \qquad (\omega, \omega')_{p,q} = (u, u')_p - (v, v')_q = \sum_{l=1}^{p} \overline{u_l'} u_l - \sum_{l=1}^{q} \overline{v_l'} v_l.$$

Here $\overline{x}$ denotes quaternionic conjugate of $x \in \mathbb{H}$. Also, since $\mathbb{H}$ is noncommutative we should specify that we consider $\mathbb{H}^{p+q}$ as a right vector space over $\mathbb{H}$. Specifically, quaternionic scalar multiplication multiplies the coordinate of $\omega \in \mathbb{H}$ on the right:

$$\omega h = (u_1 h, u_2 h, \ldots, u_p h, v_1 h, \ldots, v_q h), \qquad h \in \mathbb{H}.$$

This is the reason we write $\overline{u_l'} u_l$ rather than $u_l \overline{u_l'}$ in formula (5.1). With this convention, the algebra $\mathrm{End}(\mathbb{H}^{p+q})$ may be realized as the space of $(p+q) \times (p+q)$ matrices with entries in $\mathbb{H}$, acting on $\omega \in \mathbb{H}^{p+q}$ by the usual rules of matrix multiplication.

Let

$$(5.2) \qquad X^0 = \left\{ \omega \in \mathbb{H}^{p+q} - \{0\} \,\big|\, (\omega, \omega)_{p,q} = 0 \right\}$$

be the light cone in $\mathbb{H}^{p+q}$. If $\boldsymbol{i}, \boldsymbol{j}$ and $\boldsymbol{k} = \boldsymbol{ij}$ denote the standard quaternion units, we can write

$$u_l = x_{4l-3} + \boldsymbol{i} x_{4l-2} + \boldsymbol{j}(x_{4l-1} + \boldsymbol{i} x_{4l}) = z_{2l-1} + \boldsymbol{j} z_{2l}$$

and

$$v_l = y_{4l-3} + \boldsymbol{i} y_{4l-2} + \boldsymbol{j}(y_{4l-1} + \boldsymbol{i} y_{4l}) = w_{2l-1} + \boldsymbol{j} w_{2l}.$$

We may use the $x_k$ and $y_k$ to identify $\mathbb{H}^{p+q}$ with $\mathbb{R}^{4p+4q}$ or the $z_k$ and $w_k$ to identify $\mathbb{H}^{p+q}$ with $\mathbb{C}^{2p+2q}$. Under these identifications the light cone $X^0$ is exactly the same light cone as previously defined for $\mathbb{R}^{4p+4q}$ or for $\mathbb{C}^{2p+2q}$.

Let $\mathrm{Sp}(p,q) \subset \mathrm{GL}(p+q, \mathbb{H})$ be the isometry group of the Hermitian form $(\cdot, \cdot)_{p,q}$. The identifications of the previous paragraph give us embeddings

$$\mathrm{Sp}(p,q) \subset \mathrm{U}(2p, 2q) \subset O(4p, 4q).$$

In particular, the group $\mathrm{Sp}(p,q)$ preserves the light cone $X^0$. In fact, it acts transitively on $X^0$ [Ja].

We can deduce in the usual fashion an action of $\mathrm{Sp}(p,q)$ on functions on $X^0$, that is,

$$(5.3) \qquad (g \cdot f)(\omega) = f(g^{-1}\omega), \qquad g \in \mathrm{Sp}(p,q), \ \omega \in \mathbb{H}^{p+q}.$$

The quaternionic scalar dilations $\mathbb{H}^\times$ also acts on $X^0$, and on functions on $X^0$, and this action of $\mathbb{H}^\times$ commutes with the action of $\mathrm{Sp}(p,q)$, so we can use $\mathbb{H}^\times$ to define invariant subspaces for $\mathrm{Sp}(p,q)$.

Since $\mathbb{H}^\times$ is not commutative, we must consider $\mathbb{H}^\times$-isotypic spaces instead of looking at only eigenspaces. Specifically, write $\mathbb{H}^\times \simeq \mathbb{R}_+^\times \cdot \mathbb{H}_1$ where $\mathbb{H}_1 \simeq \mathrm{Sp}(1) \simeq \mathrm{SU}(2)$ is the group of quaternions of norm 1. For each integer $j \geq 0$, there is a unique (up to equivalence) irreducible representation $\sigma_j$ of $\mathbb{H}_1$ of dimension $j+1$. Realize $\sigma_j$ on a space $Y_j$. For $a \in \mathbb{C}$, define a representation $\sigma_j^a$ of $\mathbb{H}^\times$ on $Y_j$ by

$$\sigma_j^a(h) = |h|^a \sigma_j(h/|h|), \qquad h \in \mathbb{H}^\times,$$

where $|h|$ denotes the quaternionic norm of $h$. The representation $\sigma_j$ of $\mathbb{H}_1$ is self-dual. Consequently, there is a nondegenerate bilinear form $\langle \cdot, \cdot \rangle_j$ (symmetric if $j$ is even, skew-symmetric if $j$ is odd) satisfying

$$\langle \sigma_j^a(\overline{h})y, y' \rangle_j = \langle y, \sigma_j^a(h)y' \rangle_j, \qquad y, y' \in Y_j, \ h \in \mathbb{H}.$$

Define a space $S_j^a(X^0; Y_j)$ of smooth $Y_j$-valued functions on $X^0$ by the recipe

$$S_j^a(X^0; Y_j) = \left\{ f \in C^\infty(X^0; Y_j) \mid f(\omega h) = \sigma_j^a(\overline{h})f(\omega), \ h \in \mathbb{H}^\times, \ \omega \in X^0 \right\}.$$

Define also a map

$$\beta : \ S_j^a(X^0; Y_j) \otimes Y_j \longrightarrow C^\infty(X^0)$$

by

$$\beta(f \otimes y)(\omega) = \langle f(\omega), y \rangle_j.$$

We may define actions of $\mathrm{Sp}(p,q)$ on $S_j^a(X^0; Y_j) \subset C^\infty(X^0; Y_j)$ and on $C^\infty(X^0)$ by means of formula (5.3). It is trivial to check that $\beta$ intertwines these actions. Furthermore, we may compute

$$(5.4) \qquad \begin{aligned} \beta(f \otimes y)(\omega h) &= \langle f(\omega h), y \rangle_j \\ &= \langle \sigma_j^a(\overline{h})f(\omega), y \rangle_j \\ &= \langle f(\omega), \sigma_j^a(h)y \rangle_j \\ &= \beta(f \otimes \sigma_j^a(h)(y))(\omega). \end{aligned}$$

Let us denote the image of $\beta$ by $S_j^a(X^0)$. It is not difficult to check that $\beta$ is a linear isomorphism from $S_j^a(X^0; Y_j) \otimes Y_j$ to $S_j^a(X^0)$. Formula (5.4) therefore shows that the joint action of $\mathrm{Sp}(p,q) \times \mathbb{H}^\times$ on $S_j^a(X^0)$ is isomorphic to the tensor product of the action of formula (5.3) of $\mathrm{Sp}(p,q)$ on $S_j^a(X^0; Y_j)$ with the action of $\sigma_j^a$ of $\mathbb{H}^\times$ on $Y_j$.

Formula (5.4) implies that $S_j^a(X^0) \subset S^a(X^0)$. Here $S^a(X^0)$ is as in formula (2.7). Standard arguments from harmonic analysis on compact groups tell us that

$$S^a(X^0) = \sum_{j \geq 0} S_j^a(X^0),$$

the sum being in the sense of topological vector spaces. We should also note that

$$S^{a+}(X^0) = \sum_{\substack{j \geq 0 \\ j \text{ even}}} S_j^a(X^0)$$

and

$$S^{a-}(X^0) = \sum_{\substack{j \geq 0 \\ j \text{ odd}}} S_j^a(X^0).$$

We call $S_j^a(X^0)$ the $\sigma_j^a$-*isotypic component* of $C^\infty(X^0)$. We want to describe the $S_j^a(X^0)$ as modules for $\mathrm{Sp}(p, q)$. According to our construction of $S_j^a(X^0)$, this is essentially equivalent to describing $S_j^a(X^0; Y_j)$ as an $\mathrm{Sp}(p, q)$ module.

Our procedure for analyzing the $\mathrm{Sp}(p, q)$-module structure of $S_j^a(X^0; Y_j)$ will follow the pattern used before in dealing with $\mathrm{O}(p, q)$ and $\mathrm{U}(p, q)$. Indeed, since $\mathrm{Sp}(p, q)$ may be regarded as a subgroup of $\mathrm{U}(2p, 2q)$ and $\mathrm{O}(4p, 4q)$, much of what we need is already done. We only give here the additional arguments that are required. In §5.2 we analyze the $K$-type structure, in §5.3 we discuss composition series and unitarity in the case $p, q > 1$, and in §5.4 we discuss the groups $\mathrm{Sp}(p, 1)$, which are particularly interesting because of degeneracy in the $K$-types similar to that already described for $\mathrm{U}(p, 1)$ in §4.5.

Before beginning our detailed description of $S_j^a(X^0; Y_j)$, we note that they can be thought of as induced representations, as for the representations of $\mathrm{U}(p, q)$ and $\mathrm{O}(p, q)$ previously studied. Precisely, fix a point $\omega_0 \in X^0$ and let $Q$ be the stabilizer in $\mathrm{Sp}(p, q)$ of $\omega_0$. Let $P$ be the stabilizer of the quaternionic line $\omega_0 \mathbb{H}$ through $\omega_0$. Define a map $\tilde{\alpha} : P \longrightarrow \mathbb{H}^\times$ by the formula $p\omega_0 = \omega_o \tilde{\alpha}(p)$. The map $\tilde{\alpha}$ is, in fact, a group isomorphism $\tilde{\alpha} : P/Q \to \mathbb{H}^\times$. Given any representation $\sigma_j^a$ of $\mathbb{H}^\times$, the composition $\sigma_j^a \circ \tilde{\alpha}$ defines a representation of $P$. It is straightforward to check that the space $S_j^a(X^0; Y_j)$ is isomorphic as an $\mathrm{Sp}(p, q)$ module to the representation induced from $\sigma_j^{-a} \circ \tilde{\alpha}$ of $P$.

*Remark.* When $p = q = 1$, the above representations constitute all principal series representations, so our analysis gives a complete account of the principal series and of the unitary dual of $\mathrm{Sp}(1, 1)$. Since $\mathrm{Sp}(1, 1) \simeq \mathrm{Spin}(4, 1)$ is the two-fold cover of the identity component of $\mathrm{O}(4, 1)$, this case reproduces the calculations of Dixmier [Di].

## 5.2. $K$-structure of $S_j^a(X^0)$.

To investigate the structure of $S_j^a(X^0)$ as $\mathrm{Sp}(p, q)$ module, we first consider its $K$-structure. Here $K$ will designate the maximal compact subgroup $\mathrm{Sp}(p) \times \mathrm{Sp}(q)$ of $\mathrm{Sp}(p, q)$ consisting of elements of $\mathrm{Sp}(p, q)$ that stabilize the decomposition $\mathbb{H}^{p+q} \simeq \mathbb{H}^p \oplus \mathbb{H}^q$. Our experience in §§2–4 suggests that we try to extract an $\mathrm{Sp}(p) \times \mathbb{H}_1$ decomposition of $\mathcal{H}^m(\mathbb{R}^{4p})$. We start with some notation.

Let $\mathcal{H}(\mathbb{H}^p) = \{f \in \mathcal{P}(\mathbb{H}^p) \mid \Delta f = 0\}$. Here $\Delta$ is the Laplacian with respect to the real coordinate system defined after formula (5.2). We observe that $\mathcal{H}(\mathbb{H}^p)$ is precisely the $\mathrm{U}(2p)$ harmonics in $\mathcal{P}(\mathbb{C}^{2p})$ and is also the $\mathrm{O}(4p)$ harmonics in $\mathcal{P}(\mathbb{R}^{4p})$.

It is well known (see [Zh]) that irreducible finite-dimensional representations of $\mathrm{Sp}(p)$ can be indexed (using the highest weight) by a $p$-tuple of integers $(\xi_1, \ldots, \xi_p)$ with $\xi_1 \geq \xi_2 \geq \cdots \geq \xi_p \geq 0$. We shall write a representation of highest weight $(\xi_1, \xi_2, 0, \ldots, 0)$ as $V_p^{(\xi_1, \xi_2)}$ [note that the subscript indicates the group $\mathrm{Sp}(p)$]. For $\mathbb{H}_1 \simeq \mathrm{Sp}(1) \simeq \mathrm{SU}(2)$, we denote by $V_1^k$, the unique irreducible representation of $\mathrm{SU}(2)$ of dimension $k+1$.

**Proposition 5.1.** *We have the following decomposition for $\mathcal{H}^m(\mathbb{R}^{4p})$ as $\mathrm{Sp}(p) \times \mathbb{H}_1 \simeq \mathrm{Sp}(p) \times \mathrm{Sp}(1)$ module:*

$$(5.5) \qquad \mathcal{H}^m(\mathbb{R}^{4p})\Big|_{\mathrm{Sp}(p) \times \mathbb{H}_1} \simeq \begin{cases} \sum_{\substack{\xi_1 \geq \xi_2 \geq 0 \\ \xi_1 + \xi_2 = m}} V_p^{(\xi_1, \xi_2)} \otimes V_1^{\xi_1 - \xi_2}, & p > 1, \\ V_1^m \otimes V_1^m, & p = 1. \end{cases}$$

*Hence, as $\mathrm{Sp}(p) \times \mathbb{H}_1$ module,*

$$(5.6) \qquad C^\infty(\mathbb{S}^{4p-1}) \simeq \begin{cases} \sum_{\xi_1 \geq \xi_2 \geq 0} V_p^{(\xi_1, \xi_2)} \otimes V_1^{\xi_1 - \xi_2}, & p > 1, \\ \sum_{m \geq 0} V_1^m \otimes V_1^m, & p = 1. \end{cases}$$

*Proof.* See §6.3. □

Assume for the moment that $p, q \geq 2$. Combining the decomposition (5.5) for $p$ and for $q$ and an application of Clebsch Gordan Formula [Zh] yields

$$(5.7) \qquad \begin{aligned} &\mathcal{H}^m(\mathbb{R}^{4p}) \otimes \mathcal{H}^n(\mathbb{R}^{4q}) \\ &= \sum_{\substack{\xi_1 \geq \xi_2 \geq 0, \ \xi_1 + \xi_2 = m \\ \eta_1 \geq \eta_2 \geq 0, \ \eta_1 + \eta_2 = n}} V_p^{(\xi_1, \xi_2)} \otimes V_q^{(\eta_1, \eta_2)} \otimes (V_1^{\xi_1 - \xi_2} \otimes V_1^{\eta_1 - \eta_2}) \\ &= \sum_{0 \leq j \leq m+n} \Big( \sum \big\{ V_p^{(\xi_1, \xi_2)} \otimes V_q^{(\eta_1, \eta_2)} \big\} \Big) \otimes V_1^j, \end{aligned}$$

where the sum of the terms in braces is over

$$(5.8) \qquad \begin{aligned} \xi_1 \geq \xi_2 &\geq 0, & \eta_1 \geq \eta_2 \geq 0, \\ \xi_1 + \xi_2 &= m, & \eta_1 + \eta_2 = n, \\ \xi_1 - \xi_2 + \eta_1 - \eta_2 &\geq j \geq |(\xi_1 - \xi_2) - (\eta_1 - \eta_2)|, \\ j &\equiv \xi_1 - \xi_2 + \eta_1 - \eta_2 \mod 2. \end{aligned}$$

Note that for a fixed $j$, there is exactly one copy of $V_p^{(\xi_1, \xi_2)} \otimes V_q^{(\eta_1, \eta_2)} \otimes V_1^j$ for each 4-tuple $(\xi_1, \xi_2, \eta_1, \eta_2)$ satisfying (5.8). It is apparent now what the $K$-structure of $S_j^a(X^0)$ is.

Define, for each 4-tuple $(\xi_1, \xi_2, \eta_1, \eta_2)$ satisfying (5.8),

$$(5.9) \qquad \begin{aligned} i_{a,j} \ : \ &V_p^{(\xi_1, \xi_2)} \otimes V_q^{(\eta_1, \eta_2)} \otimes V_1^j \to S_j^a(X^0), \\ &i_{a,j}(\phi) = \phi r_p^{2\alpha}, \end{aligned}$$

where

$$\phi \in V_p^{(\xi_1, \xi_2)} \otimes V_q^{(\eta_1, \eta_2)} \otimes V_1^j \subset \mathcal{H}^m(\mathbb{R}^{4p}) \otimes \mathcal{H}^n(\mathbb{R}^{4q}),$$

$$r_p^2 = \sum_{l=1}^{2p} z_l \overline{z_l}, \quad \text{and} \quad a = \xi_1 + \xi_2 + \eta_1 + \eta_2 + 2\alpha.$$

The following statement is clear from the discussion above.

**Lemma 5.2.** *As an* $\mathrm{Sp}(p) \times \mathrm{Sp}(q) \times \mathbb{H}_1$ *module,* $S_j^a(X^0)$ *decomposes into a direct sum*

$$S_j^a(X^0) = \sum i_{a,j}(V_p^{(\xi_1, \xi_2)} \otimes V_q^{(\eta_1, \eta_2)} \otimes V_1^j),$$

*where the sum is over the set described in* (5.8) *(see Diagram* 5.10*).*

For reasons that will be apparent after the description of the transition properties of $\mathfrak{p} \subset \mathfrak{sp}(p,q)$, we continue to organize $K$-types of $S_j^a(X^0)$ by means of the points $(m,n) = (\xi_1 + \xi_2, \eta_1 + \eta_2)$ in $(\mathbb{Z}_+)^2$. In other words, the $K$-types attached to $(m,n)$ are those satisfying (5.8). We represent the 4-tuples $(\xi_1, \xi_2, \eta_1, \eta_2)$ satisfying (5.8) as the integer points $(\xi_1 - \xi_2, \eta_1 - \eta_2)$ in an auxiliary rectangle

$$R_{m,n} = \{\,(x,y) \mid 0 \le x \le m, \ 0 \le y \le n\,\}.$$

These points fill a coset of $(2\mathbb{Z})^2$ in a strip along the diagonal line $x = y$, as pictured in Diagram 5.10. We refer to these $K$-types as the *fiber over* $(m,n)$.

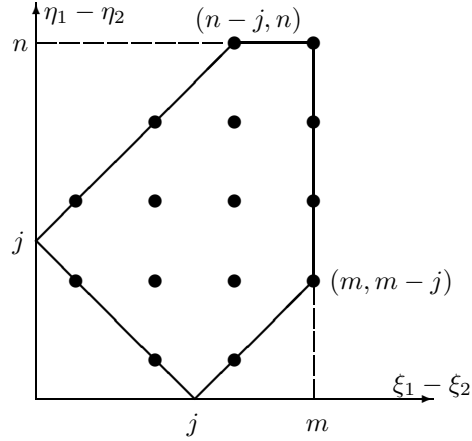

DIAGRAM 5.10. Fiber of $K$-types over $(m,n)$ (if $j \equiv m+1 \mod 2$).

Diagram 5.10 helps us to see that the fiber over $(m,n)$ is nonempty if and only if

(5.11)   (i)   $m + n \equiv j \mod 2$,
          (ii)  $m + n \ge j$.

Thus, when $p \ge q \ge 2$, the $K$-types of $S_j^a(X^0)$ fall in the region depicted in Diagram 5.12.

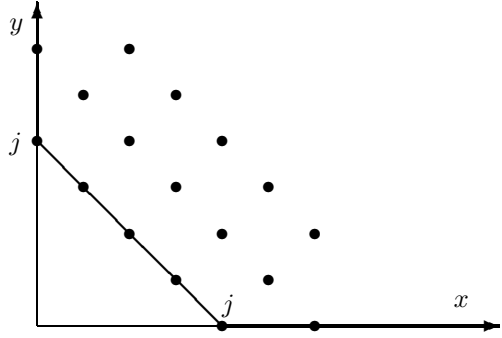

DIAGRAM 5.12. $K$-type region of $S_j^a(X^0)$ when $p, q \geq 2$ [consisting of points $(m, n)$ where $m + n \equiv j \mod 2$].

Now consider the case $p > q = 1$. The computation parallel to (5.7) is

$$
\begin{aligned}
\mathcal{H}^m(\mathbb{R}^{4p}) &\otimes \mathcal{H}^n(\mathbb{R}^4) \\
&= \sum_{\substack{\xi_1 \geq \xi_2 \geq 0, \\ \xi_1 + \xi_2 = m}} V_p^{(\xi_1, \xi_2)} \otimes V_1^n \otimes (V_1^{\xi_1 - \xi_2} \otimes V_1^n) \\
&= \sum_{0 \leq j \leq m + n} \Big( \sum \big\{ V_p^{(\xi_1, \xi_2)} \otimes V_1^n \big\} \Big) \otimes V_1^j,
\end{aligned}
$$

where the sum of the terms in braces is over

$$
\begin{aligned}
n \geq 0, \qquad \xi_1 &\geq \xi_2 \geq 0, \qquad \xi_1 + \xi_2 = m, \\
\xi_1 - \xi_2 + n &\geq j \geq |(\xi_1 - \xi_2) - n|, \\
\xi_1 - \xi_2 + n &\equiv j \mod 2.
\end{aligned}
$$

We again associate this set of $K$-types of $S_j^a(X^0)$ to $(m, n) \in (\mathbb{Z}_+)^2$. If we label these $K$-types in the fiber over $(m, n)$ by $\xi_1 - \xi_2 = l$, then we obtain all $l$ such that

$$(5.13) \qquad\qquad |j - n| \leq l \leq \min(m, j + n)$$

and congruent to $j + n$ modulo 2. This is simply the set of $K$-types on the line $\eta_1 - \eta_2 = n$ in Diagram 5.10. So if $p > q = 1$, the set of points $(m, n)$ such that the associated fiber of $K$-types is nonempty can be depicted by Diagram 5.14.

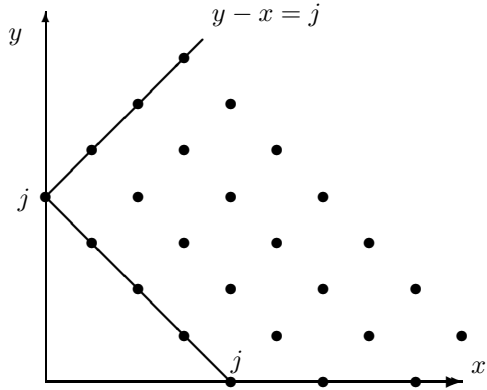

DIAGRAM 5.14. $K$-type region for $S_j^a(X^0)$ when $p > q = 1$ [consisting of points $(m, n)$ where $m + n \equiv j \mod 2$].

For the case $p = q = 1$, the same arguments show that the $K$-types of $S_j^a(X^0)$ are $V_1^m \otimes V_1^n$, where

$$(5.15) \qquad \begin{aligned} m + n \geq j \geq |m - n|, \\ m + n \equiv j \mod 2. \end{aligned}$$

The $K$-type region is then as depicted in Diagram 5.16.

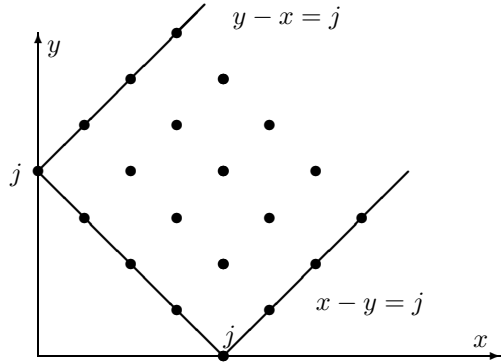

DIAGRAM 5.16. $K$-type region for $S_j^a(X^0)$ when $p = q = 1$ [consisting of points $(m, n)$ where $m + n \equiv j \mod 2$].

Furthermore, in this case the fiber over $(m, n)$ consists of a single $K$-type, viz., $V_1^m \otimes V_1^n$.

### 5.3. Composition structure and unitarity.

Since the technique involved in the consideration of unitarity is the same as in §§3 and 4.4, we will simply state the results alongside our discussion of the composition series structure. Write

$$V_{a,j}^{(\xi_1, \xi_2, \eta_1, \eta_2)} = i_{a,j} \left( V_p^{(\xi_1, \xi_2)} \otimes V_q^{(\eta_1, \eta_2)} \otimes V_1^j \right) \quad \text{if } p, q > 1,$$

$$V_{a,j}^{(\xi_1, \xi_2, \eta)} = i_{a,j} \left( V_p^{(\xi_1, \xi_2)} \otimes V_1^\eta \otimes V_1^j \right) \quad \text{if } p > q = 1,$$

$$V_{a,j}^{(\xi, \eta)} = i_{a,j} \left( V_1^\xi \otimes V_1^\eta \otimes V_1^j \right) \quad \text{if } p = q = 1,$$

and let

$$U_{a,j}^{(m,n)} = \begin{cases} \sum_{\substack{\xi_1 + \xi_2 = m \\ \eta_1 + \eta_2 = n}} V_{a,j}^{(\xi_1, \xi_2, \eta_1, \eta_2)}, & p \geq q > 1, \\ \sum_{\xi_1 + \xi_2 = m} V_{a,j}^{(\xi_1, \xi_2, n)}, & p > q = 1, \end{cases}$$

be the sum of $K$-types of $S_j^a(X^0)$ living in the fiber over $(m, n)$. We shall discuss the cases $p \geq q \geq 2$, leaving the interesting cases $p \geq q = 1$ to §5.4. Recall the functions $A_{p,q,a}^{\pm\pm}$ from §2 [see (2.20)]. In parallel with the cases of $\mathfrak{o}(p, q)$ and $\mathfrak{u}(p, q)$, we let $\mathfrak{p}$ be the space complementary to $(\mathfrak{sp}(p) \oplus \mathfrak{sp}(q))_\mathbb{C}$ in $\mathfrak{sp}(p, q)_\mathbb{C}$. The following lemma is clear from the decomposition in (5.7) and Lemma 2.3.

**Lemma 5.3.** *There exist maps, nonzero whenever the target is nonzero,*

$$T_{m,n}^{\pm\pm} : \ \mathfrak{p} \otimes U_{a,j}^{(m,n)} \to U_{a,j}^{(m\pm 1, n\pm 1)}$$

*such that the action of $z \in \mathfrak{p}$ on a $K$-type is given by*

$$\rho(z)\phi = A_{4p,4q,a}^{++}(\xi_1 + \xi_2, \ \eta_1 + \eta_2) T_{\xi_1 + \xi_2, \ \eta_1 + \eta_2}^{++}(z \otimes \phi)$$
$$+ A_{4p,4q,a}^{+-}(\xi_1 + \xi_2, \ \eta_1 + \eta_2) T_{\xi_1 + \xi_2, \ \eta_1 + \eta_2}^{+-}(z \otimes \phi)$$
$$+ A_{4p,4q,a}^{-+}(\xi_1 + \xi_2, \ \eta_1 + \eta_2) T_{\xi_1 + \xi_2, \ \eta_1 + \eta_2}^{-+}(z \otimes \phi)$$
$$+ A_{4p,4q,a}^{--}(\xi_1 + \xi_2, \ \eta_1 + \eta_2) T_{\xi_1 + \xi_2, \ \eta_1 + \eta_2}^{--}(z \otimes \phi),$$

*where $\phi \in V_{a,j}^{(\xi_1, \xi_2, \eta_1, \eta_2)} \subset U_{a,j}^{(\xi_1 + \xi_2, \ \eta_1 + \eta_2)}$.*

So the action of $\mathfrak{sp}(p, q)$ [$\subset \mathfrak{u}(2p, 2q) \subset \mathfrak{o}(4p, 4q)$] on each fiber over $(m, n)$ is compatible with Diagram 2.19. Since each fiber has several $K$-types, however, in order to have a grasp on the $\mathfrak{sp}$-module structure we must calculate how the $\mathfrak{sp}(p, q)$ action affects the various points in each fiber. This is the essence of the next lemma, which we will prove in §6.3.

**Lemma 5.4.** *Consider a fiber $U_{a,j}^{(m,n)}$ of $K$-types in $S_j^a(X^0)$. By successive applications of the maps $T_{m+1,n+1}^{--}$, $T_{m,n}^{++}$, we can move from any $K$-type in $U_{a,j}^{(m,n)}$ to any other. Hence, if the transition coefficients $A_{4p,4q,a}^{++}(m, n)$ and $A_{4p,4q,a}^{--}(m+1, \ n+1)$ are both nonzero, application of appropriate elements of the enveloping algebra $\mathfrak{U}(\mathfrak{sp})$ allows us to move from any $K$-type in $U_{a,j}^{(m,n)}$ to any other. Similar remarks apply to transitions from $U_{a,j}^{(m,n)}$ to any of the adjacent fibers $U_{a,j}^{(m\pm 1, n\pm 1)}$.*

*Remark.* The adjoint action of $K$ on $\mathfrak{p}$ makes $\mathfrak{p}$ an irreducible $K$-module. The $K$-types that can be reached from $V_{a,j}^{(\xi_1,\xi_2,\eta_1,\eta_2)}$ by one application of $\mathfrak{p}$ are isomorphic to constituents of $V^{\sim} = \mathfrak{p} \otimes V_{a,j}^{(\xi_1,\xi_2,\eta_1,\eta_2)}$. Since the $K$-weights of $\mathfrak{p}$ are of the form $\pm(e_i \pm f_j)$ (see §6.3 for notation), it is not difficult to see the $K$-types that are embeddable in $V^{\sim}$ are among the 16 whose highest weights are obtained by adding one of the vectors $(\pm 1, 0, \pm 1, 0)$, $(\pm 1, 0, 0, \pm 1)$, $(0, \pm 1, \pm 1, 0)$, or $(0, \pm 1, 0, \pm 1)$ to $(\xi_1, \xi_2, \eta_1, \eta_2)$. We will see that by moving from $U_{a,j}^{(m,n)}$ to $U_{a,j}^{(m+1,n+1)}$ by means of $\mathfrak{p}$ and then back again we can move from $V_{a,j}^{(\xi_1,\xi_2,\eta_1,\eta_2)}$ to any of its eight nearest neighbors, as illustrated in Diagram 5.17 (see also Diagram 5.10).

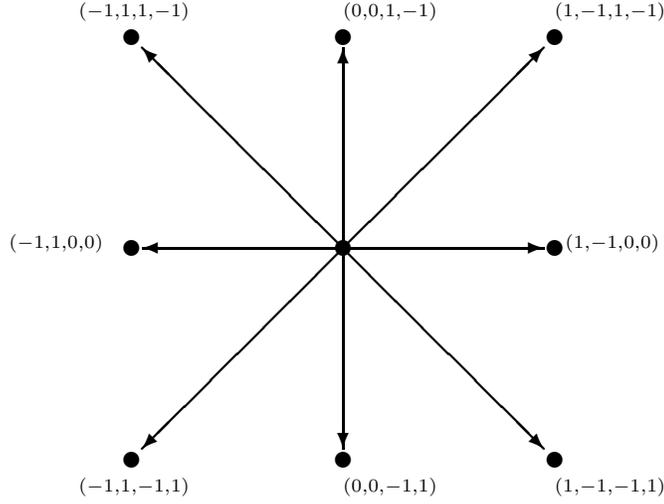

Diagram 5.17

In this diagram, the labels on the arrows indicate increments to $(\xi_1, \xi_2, [0]\eta_1, \eta_2)$ and the arrows show which increments can be achieved by applying $T_{m+1,n+1}^{--} T_{m,n}^{++}$. Thus the arrow to $(1, -1, 0, 0)$ indicates that we can move from $V_{a,j}^{(\xi_1,\xi_2,\eta_1,\eta_2)}$ to $V_{a,j}^{(\xi_1+1,\xi_2-1,\eta_1,\eta_2)}$ and so forth. Similar pictures apply to the other transitions. Of course, if $(\xi_1, \xi_2, \eta_1, \eta_2)$ lies on the boundary of the $K$-type region (see Diagram 5.10), then transitions that would move it outside the region are suppressed.

Study of Diagrams 5.10 and 5.17 will convince the reader that there are no barriers to movement between $K$-types in a given fiber, as long as one can move back and forth between the fiber and an adjacent fiber, but from the geometry of the barriers (see Diagram 2.22) we can see that the only occasion on which one cannot pass between a fiber and an adjacent one is for $U_{0,0}^{(0,0)}$—and this consists only of the trivial $K$-type, so there is nothing to prove in this case. Thus by using Lemma 5.4, we can derive the following result.

**Theorem 5.5.** (a) *The image in any* $\mathrm{O}(4p, 4q)$*-irreducible component of* $S^a(X^0)$ *of a given* $S_j^a(X^0)$ *is an irreducible module for* $\mathrm{Sp}(p, q) \times \mathbb{H}_1$. *In particular, if* $p \geq q > 1$*, then* $S_j^a(X^0)$ *is reducible if and only if* $a \in \mathbb{Z}$ *and* $a \equiv j \mod 2$.

(b) $S_j^a(X^0)$ *will have Hermitian constituents if and only if* $a$ *is on the unitary axis* (*in which case* $S_j^a(X^0)$ *will be unitary*) *or* $a$ *is real* (*in which case all its constituents will be Hermitian*). *A given* $\mathrm{Sp}(p,q)$ *constituent of* $S_j^a(X^0)$ *is unitary if and only if the* $\mathrm{O}(4p,4q)$ *transitions between the nonempty fibers of the* $\mathrm{Sp}(p,q)$ *constituent are consistent with unitarity.*

*Remark.* As in the case of $\mathrm{U}(p,q)$ modules $S^{\alpha,\beta}(X^0)$, one obtains an $\mathrm{Sp}(p,q)$ composition series for an $\mathrm{Sp}(1)$ eigenspace of $S^a(X^0)$ by intersecting it with the $\mathrm{O}(4p,4q)$ composition series for $S^a(X^0)$.

*Proof.* Again this follows essentially from Lemmas 5.3 and 5.4. The pictures when $p \geq q \geq 2$ will be similar to those under Case EE in §§2 and 3 (see Diagrams 3.13, 4.20, and 4.21). ☐

**5.4. The case of** $\mathrm{Sp}(p,1)$. The $K$-type regions are depicted in Diagrams 5.14 and 5.16. For $p > q = 1$, there is a "line" of $K$-types over each point in the $K$-type region. We note that unlike the situation of $\mathrm{U}(p,1)$, $K$-types along each fiber can be moved to adjacent ones via operators in $\mathfrak{U}(\mathfrak{sp}(p,1))$. The analogs of Lemmas 5.3 and 5.4 hold. We do not state them formally. The following result is easy.

**Theorem 5.6.** *Assume* $p \geq q = 1$. *If* $a \notin \mathbb{Z}$ *or* $a \not\equiv j \mod 2$, *then* $S_j^a(X^0)$ *is irreducible for* $\mathrm{Sp}(p,1)$.

Following §4.5, we organize our description of the structure of $\mathrm{Sp}(p,1)$, $p > 1$, by means of a picture (see Diagram 5.18) describing the Hasse diagrams of $S_j^a(X^0)$ in the $(a,j)$-plane ($a \in \mathbb{R}$). Observe that the diagram is symmetric about the axis $a = -2p-1$, because of the duality between the modules $S_j^{-2p-1+s}(X^0)$ and $S_j^{-2p-1-s}(X^0)$ (see Lemma 3.1).

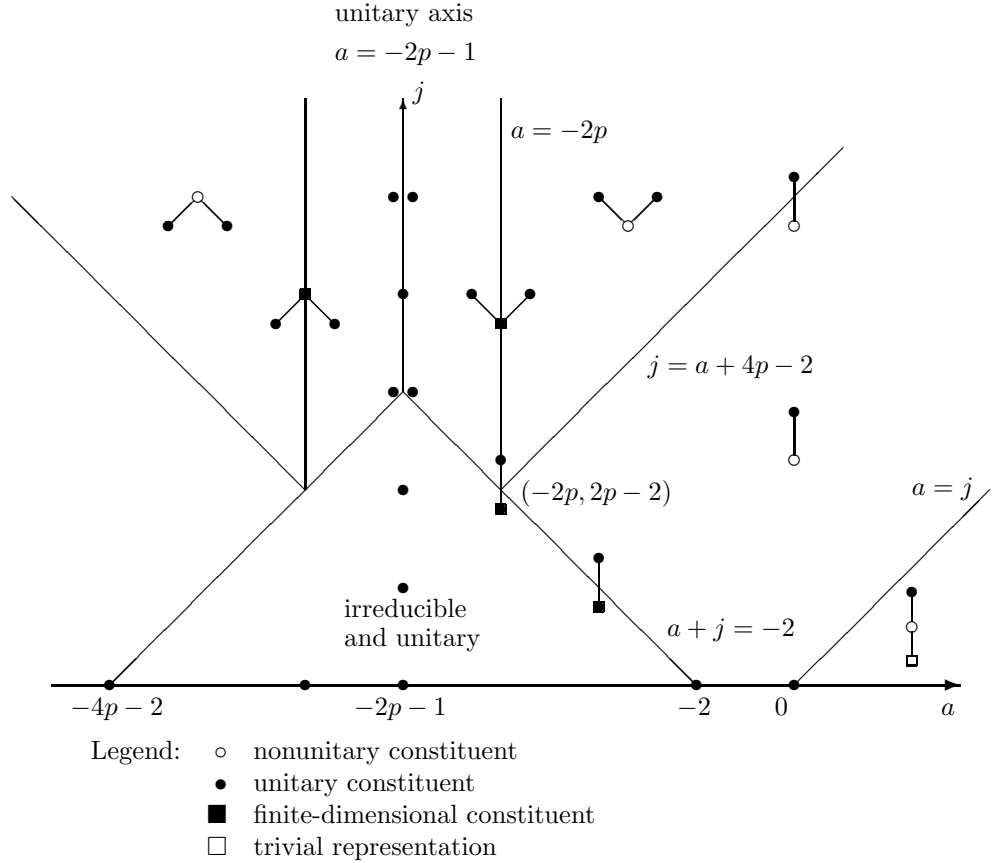

DIAGRAM 5.18. Submodule structure of $S_j^a(X^0)$ when $a \equiv j$ mod 2 for $\mathrm{Sp}(p, 1)$.

Again the existence of the complementary series is very interesting. Notice that in the triangular region given by $a - j \geq -2 - 4p$ and $a + j \leq -2$, the modules are irreducible and unitarizable. The length of a horizontal line in this triangle will be twice the length of the complementary series for a fixed $j$, that is, there is a complementary series of length $2p - j - 1$ for each $j \leq 2p - 2$. The extreme case when $j = 0$, that is, the spherical principal series, has complementary series of length $2p - 1$. Unlike the case for $\mathrm{U}(2p, 1)$ module $S^{\alpha, \beta}(X^0)$, where $\alpha = \beta$, the complementary series does not run till $a = 0$. This can be explained by observing that, for $\mathrm{U}(p, 1)$, the barriers $A^{++}$ and $A^{+-}$ always intersect on the $x$-axis, while for $\mathrm{Sp}(p, 1)$ there is a gap of length 2 between the $x$-intercepts of the $A^{++}$ and $A^{+-}$ barriers. This means that as $a$ increases from the unitary axis it is always the $A^{+-}$ barrier that first causes reducibility and prevents further unitarity, even in the extreme case $j = 0$.

Diagram 5.19 describes the situation for $\mathrm{Sp}(1, 1)$, with the same legends as in Diagram 5.18. We hope it will be self-explanatory.

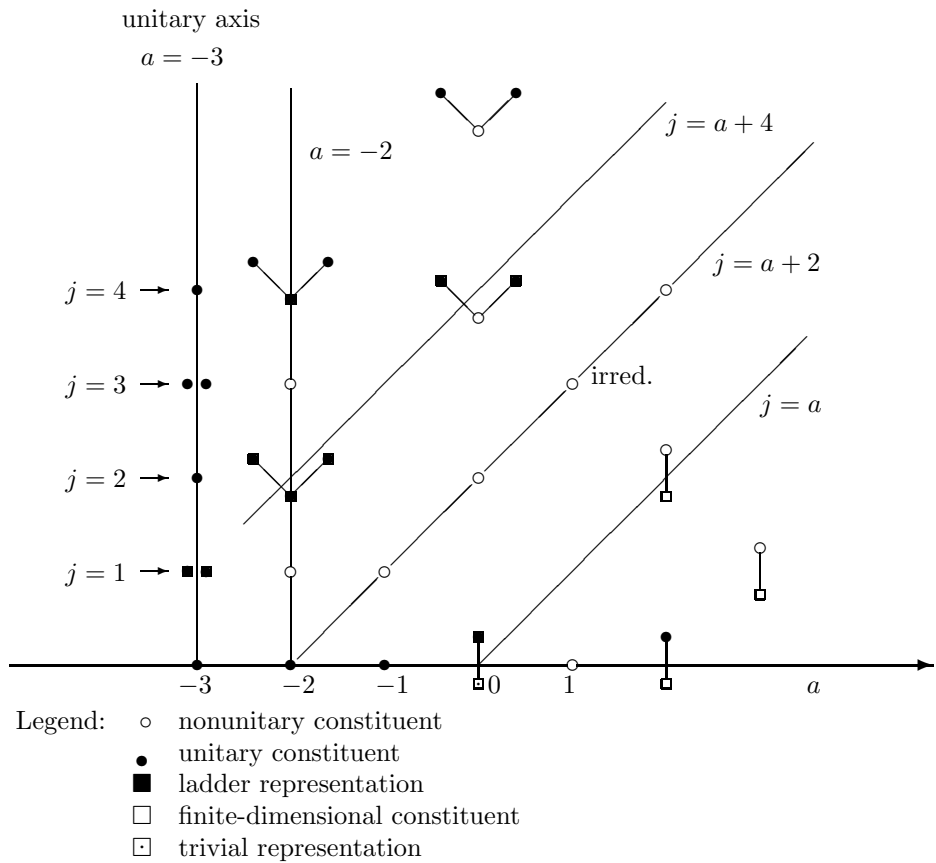

DIAGRAM 5.19. Submodule structure of $S_j^a(X^0)$ for $\mathrm{Sp}(1,1)$.

*Remarks.* We add some remarks in the case of $\mathrm{Sp}(2,1)$. Our results showed that the $K$-types of $S_j^a(X^0)$ are of the form $V_2^{(\xi_1,\xi_2)} \otimes V_1^n$ where $\xi_1, \xi_2, n$ satisfies (5.13) (see Diagram 5.12). Let us instead plot the second and third coordinates, that is, $(\xi_2, n)$, fixing $\xi_1$; this is a section of the $K$-type region parallel to the $\xi_1$ axis (see Diagram 5.20; here $\xi_1 \,(\geq j)$ is fixed and we also have the parity condition $\xi_1 + \xi_2 + n \in 2\mathbb{Z}$).

The following statements are based on the second author's thesis (see [Ta]), and we refer our readers to that for full details. We will be somewhat sketchy. If we consider arbitrary principal series representations of $\mathrm{Sp}(2,1)$, we get highest weights $(\xi_1, \xi_2, n)$ that live in regions whose convex hulls are wedges in $\mathbb{R}^3$. If we take sections parallel to the $\xi_1$-axis, we get Diagram 5.21 (note that Diagram 5.20 is a "degenerate" case of this diagram).

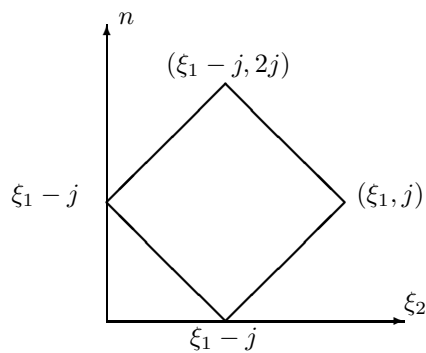

DIAGRAM 5.20

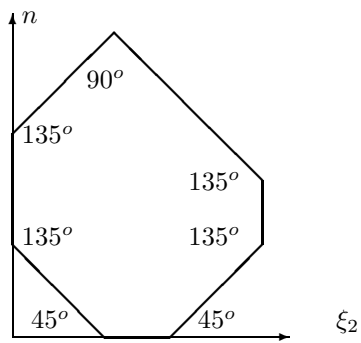

DIAGRAM 5.21

The multiplicities of $K$-types along the boundary of the figure is one, and the multiplicities "grow" linearly with the distance from the boundary, attaining a constant towards the 'central' region, much like the way the weight multiplicities of an irreducible $\mathfrak{sl}(3)$ representation behave in relation to its weight diagram. Further, if we look at $K$-type diagrams for the quotients of the principal series representations, we get diagrams like Diagram 5.22.

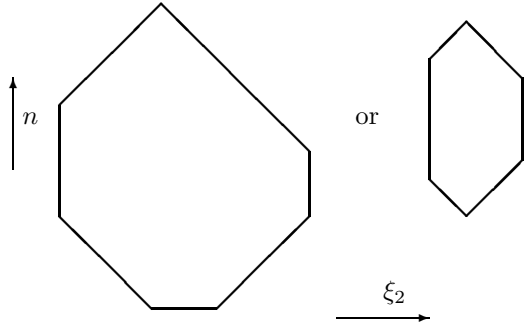

DIAGRAM 5.22

The multiplicities behave as described in the previous paragraph. Details are omitted since the results arise from painstaking computations involving Blattner's Formula, Branching Rules, and known composition series structure for $\mathrm{Sp}(2,1)$. A reason for the nice multiplicity result in this case may be the already apparent linearity in the branching rules for $\mathrm{Sp}(2) \times \mathrm{Sp}(1)$ [a maximal compact subgroup of $\mathrm{Sp}(2,1)$] to $\mathrm{Sp}(1) \times \mathrm{Sp}(1)$.

## 6. Technicalities

**6.1. Computations for $\mathrm{O}(p, q)$.** Let $x_1, \ldots, x_n$ be the standard system of coordinates in $\mathbb{R}^n$ and let

$$\Delta = \sum_{j=1}^n \frac{\partial^2}{\partial x_j^2}, \qquad r^2 = \sum_{j=1}^n x_j^2, \qquad E = \sum_{j=1}^n x_j \frac{\partial}{\partial x_j}.$$

We have the commutator

$$[\Delta, r^2] = 4E + 2n$$

and so $\{\Delta, r^2, 4E + 2n\}$ span a Lie algebra isomorphic to $\mathfrak{sl}(2, \mathbb{R})$. If $P$ is a homogeneous polynomial in $\mathcal{P}(\mathbb{R}^n)$, let $\deg P = E(P)$ be the homogeneous degree of $P$. Then

$$(6.1) \qquad \begin{aligned} \Delta(r^2 P) &= r^2 \Delta P + (4E + 2n)P \\ &= r^2 \Delta P + (4 \deg P + 2n)P. \end{aligned}$$

There are unique harmonic polynomials $h_i$ such that

$$(6.2) \qquad P = \sum_{i=1}^l h_i (r^2)^{a_i}$$

for some positive integer $l$ and $\deg P = \deg h_i + 2a_i$. To compute the $h_i$'s, we use (6.1) inductively. If $\Delta P = 0$, then $l = 1$, $h_1 = P$, and $a_1 = 0$. If $\Delta^2 P = 0$, compute

$$\Delta(r^2 \Delta P) = r^2 \Delta^2 P + (4(\deg P - 2) + 2n)\Delta P$$

$$\Rightarrow \Delta\left(P - \frac{r^2 \Delta P}{4(\deg P - 2) + 2n}\right) = 0$$

and get

(6.3)        $$P = P^\sim + r^2\left(\frac{\Delta P}{4(\deg P - 2) + 2n}\right),$$

where for convenience we will write

$$P^\sim = P - \frac{r^2 \Delta P}{4(\deg P - 2) + 2n}.$$

(It should be clear from context which space $\mathbb{R}^n$ we are working with when the definition above is being used.) One could proceed inductively using (6.1) to obtain the decomposition in (6.2) for any homogeneous P.

If $\phi = h_1 \otimes h_2 \in \mathcal{H}^m(\mathbb{R}^p) \otimes \mathcal{H}^n(\mathbb{R}^q)$, define

$$T_{ij}^{++}\phi = (x_i h_1)^\sim (y_j h_2)^\sim,$$

$$T_{ij}^{+-}\phi = c_2(x_i h_1)^\sim \frac{\partial h_2}{\partial y_j},$$

$$T_{ij}^{-+}\phi = c_1 \frac{\partial h_1}{\partial x_i}(y_j h_2)^\sim,$$

$$T_{ij}^{--}\phi = c_1 c_2 \frac{\partial h_1}{\partial x_i}\frac{\partial h_2}{\partial y_j},$$

where

$$c_1 = \frac{1}{(2\deg h_1 - 2 + p)} \quad \text{and} \quad c_2 = \frac{1}{(2\deg h_2 - 2 + q)}.$$

Extend the above definitions linearly to all elements in $\mathcal{H}^m(\mathbb{R}^p) \otimes \mathcal{H}^n(\mathbb{R}^q)$.

*Proof of Lemma* 2.3.   Suffices to show for

$$z = x_i \frac{\partial}{\partial y_j} + y_j \frac{\partial}{\partial x_i} \in \mathfrak{p}.$$

Thus, for $h = h_1 h_2 r_p^{2\gamma}$ where $h_1 \in \mathcal{H}^m(\mathbb{R}^p)$ and $h_2 \in \mathcal{H}^n(\mathbb{R}^q)$, compute

$$\left(x_i \frac{\partial}{\partial y_j} + y_j \frac{\partial}{\partial x_i}\right)(h)$$

$$= (x_i h_1)\left(\frac{\partial h_2}{\partial y_j}\right) r_p^{2\gamma} + \left(\frac{\partial h_1}{\partial x_i}\right)(y_j h_2) r_p^{2\gamma} + 2\gamma(x_i h_1)(y_j h_2) r_p^{2(\gamma-1)}.$$

From (6.3) we have

$$x_i h_1 = (x_i h_1)^\sim + c_1\left(\frac{\partial h_1}{\partial x_i}\right) r_p^2$$

and

$$y_j h_2 = (y_j h_2)^\sim + c_2 \left( \frac{\partial h_2}{\partial y_j} \right) r_q^2 = (y_j h_2)^\sim + c_2 \left( \frac{\partial h_2}{\partial y_j} \right) r_p^2$$

since $r_p^2 = r_q^2$ on $X^0$. Hence,

$$\left( x_i \frac{\partial}{\partial y_j} + y_j \frac{\partial}{\partial x_i} \right) (h)$$

$$= \left\{ (x_i h_1)^\sim + c_1 \left( \frac{\partial h_1}{\partial x_i} \right) r_p^2 \right\} \left( \frac{\partial h_2}{\partial y_j} \right) r_p^{2\gamma} + \left( \frac{\partial h_1}{\partial x_i} \right) \left\{ (y_j h_2)^\sim + c_2 \left( \frac{\partial h_2}{\partial y_j} \right) r_p^2 \right\} r_p^{2\gamma}$$

$$+ 2\gamma \left\{ (x_i h_1)^\sim + c_1 \left( \frac{\partial h_1}{\partial x_i} \right) r_p^2 \right\} \left\{ (y_j h_2)^\sim + c_2 \left( \frac{\partial h_2}{\partial y_j} \right) r_p^2 \right\} r_p^{2(\gamma-1)},$$

which gives the required result if $p \geq q \geq 2$. $\quad\square$

*Proof of Lemma* 3.1. We will first show that if $dx$, $dy$ are the rotation-invariant probability measures on $\mathbb{S}^{p-1}$ and $\mathbb{S}^{q-1}$ respectively, then

$$J(f) = \int_{\mathbb{S}^{p-1} \times \mathbb{S}^{q-1}} f(x, y) \, dx \, dy$$

is an $\mathrm{O}(p, q)$-invariant linear functional on $S^{2-p-q}(X^0)$. To this effect, let

$$X_s^0 = \left\{ (x, y) \in X^0 \,\middle|\, \tfrac{1}{s} \leq r_p(x) \leq s \right\}, \qquad s > 0.$$

The $\mathrm{O}(p, q)$-invariant measure on $X^0 \simeq \mathbb{S}^{p-1} \times \mathbb{S}^{q-1} \times \mathbb{R}_+^\times$ is $t^{p+q-3} \, dx \, dy \, dt$. Observe that for $f \in S^{2-p-q}(X^0)$

$$\frac{1}{2 \log s} \int_{X_s^0} f(tx, ty) \, t^{p+q-3} \, dx \, dy \, dt = \frac{1}{2 \log s} \left( \int_{\frac{1}{s}}^{s} \frac{1}{t} \, dt \right) J(f) = J(f).$$

If $g \in \mathrm{O}(p, q)$, the image of $X_s^0$ under $g$ is

$$g X_s^0 = \left\{ (t'x', t'y') \in X^0 \,\middle|\, x' \in \mathbb{S}^{p-1}, \ y' \in \mathbb{S}^{q-1}, \ \frac{\phi(x', y')}{s} \leq t' \leq s\phi(x', y') \right\},$$

where $\phi(x', y')$ is a smooth function on $\mathbb{S}^{p-1} \times \mathbb{S}^{q-1}$. Thus,

$$J(g^{-1} \cdot f) = \frac{1}{2 \log s} \int_{g \cdot X_s^0} f(t'x', t'y') \, t'^{(p+q-3)} \, dx' \, dy' \, dt'$$

$$= \frac{1}{2 \log s} \int_{\mathbb{S}^{p-1} \times \mathbb{S}^{q-1}} f(x', y') \left( \int_{\frac{\phi(x', y')}{s}}^{s\phi(x', y')} \frac{dt'}{t'} \right) dx' \, dy'$$

$$= \int_{\mathbb{S}^{p-1} \times \mathbb{S}^{q-1}} f(x', y') \, dx' \, dy'$$

$$= J(f),$$

which proves our earlier assertion.

For $f \in S^{-a+2-p-q}(X^0)$ and $h \in S^{\overline{a}}(X^0)$ define

$$\langle f, h \rangle = J(f \overline{h}) = \int_{\mathbb{S}^{p-1} \times \mathbb{S}^{q-1}} f(x, y) \overline{h(x, y)} \, dx \, dy.$$

Then it is easy to see that this is an $\mathrm{O}(p,q)$-invariant Hermitian pairing between $S^{-a+2-p-q}(X^0)$ with $S^{\overline{a}}(X^0)$. It gives a Hermitian inner product on $S^a(X^0)$ if

$$-a+2-p-q=\overline{a} \quad \text{or} \quad \mathrm{Re}\, a = -\frac{(p+q)}{2}+1.$$

This proves Lemma 3.1. $\quad\square$

**6.2. Computations for** $U(p,q)$**.** Define for $S_{ij}$ and $T_{ij}$ as in (4.18) and $\phi = h_1 \otimes h_2 \in \mathcal{H}^{m_1,m_2}(\mathbb{C}^p) \otimes \mathcal{H}^{n_1,n_2}(\mathbb{C}^q)$,

$$S^{1001}(S_{ij} \otimes (h_1 \otimes h_2)) = (z_i h_1)^{\sim}(\overline{w_j}h_2)^{\sim},$$

$$S^{0-101}(S_{ij} \otimes (h_1 \otimes h_2)) = c_1\frac{\partial h_1}{\partial \overline{z_i}}(\overline{w_j}h_2)^{\sim},$$

$$S^{10-10}(S_{ij} \otimes (h_1 \otimes h_2)) = c_2(z_i h_1)^{\sim}\frac{\partial h_2}{\partial w_j},$$

(6.4)
$$S^{0-1-10}(S_{ij} \otimes (h_1 \otimes h_2)) = c_1 c_2\frac{\partial h_1}{\partial \overline{z_i}}\frac{\partial h_2}{\partial w_j},$$

$$T^{0110}(T_{ij} \otimes (h_1 \otimes h_2)) = (\overline{z_i}h_1)^{\sim}(w_j h_2)^{\sim},$$

$$T^{-1010}(T_{ij} \otimes (h_1 \otimes h_2)) = c_1\frac{\partial h_1}{\partial z_i}(w_j h_2)^{\sim},$$

$$T^{010-1}(T_{ij} \otimes (h_1 \otimes h_2)) = c_2(\overline{z_i}h_1)^{\sim}\frac{\partial h_2}{\partial \overline{w_j}},$$

$$T^{-100-1}(T_{ij} \otimes (h_1 \otimes h_2)) = c_1 c_2\frac{\partial h_1}{\partial z_i}\frac{\partial h_2}{\partial \overline{w_j}}$$

where

$$c_1 = \frac{1}{(p+m_1+m_2-1)},$$

$$c_2 = \frac{1}{(q+n_1+n_2-1)},$$

$$P^{\sim} = P - \frac{r_n^2 \Delta_n P}{n + \deg_{n,z} P + \deg_{n,\bar{z}} P - 1}.$$

[Note that $P^{\sim}$ can be shown to be harmonic if $\Delta_n^2 P = 0$ through computations similar to those before (6.3). Here $\deg_{n,z} P$ and $\deg_{n,\bar{z}} P$ are the degrees in $z$ and $\bar{z}$ of $P$ in $\mathcal{P}(\mathbb{C}^n)$.] Extend the definitions (6.4) linearly to all of $\mathcal{H}^{m_1,m_2}(\mathbb{C}^p) \otimes \mathcal{H}^{n_1,n_2}(\mathbb{C}^q)$. Lemma 4.1 then follows from a straightforward computation as in Lemma 2.3.

We still need to verify that one could move freely within each fiber (see the discussion after Diagram 4.19). Write [see (4.11)]

$$V^{\alpha,\beta}_{(m_1,m_2,n_1,n_2)} = j_{\alpha,\beta}(\mathcal{H}^{m_1,m_2}(\mathbb{C}^p) \otimes \mathcal{H}^{n_1,n_2}(\mathbb{C}^q)).$$

Let us show that there is a nonzero map $U^+$ from $\mathfrak{U}(\mathfrak{u}(p,q))$ with

(6.5)
$$U^+\colon\; V^{\alpha,\beta}_{(m_1,m_2,n_1,n_2)} \to V^{\alpha,\beta}_{(m_1+1,m_2-1,n_1-1,n_2+1)}.$$

This is a movement in the south-westerly direction along the fiber depicted in Diagram 4.15. Let $\Pi^{(m_1,m_2,n_1,n_2)}$ be the projection to the $K$-type $V^{\alpha,\beta}_{(m_1,m_2,n_1,n_2)}$, that is,

$$(6.6) \qquad \Pi^{(m_1,m_2,n_1,n_2)}\colon\ S^{\alpha,\beta}(X^0) \to V^{\alpha,\beta}_{(m_1,m_2,n_1,n_2)}.$$

This projection certainly arises from the action of $K$ on $S^{\alpha,\beta}(X^0)$. Let

$$m = m_1 + m_2, \qquad n = n_1 + n_2$$

and write

$$(6.7) \qquad \begin{aligned} S^{\varepsilon_1,\varepsilon_2,\eta_1,\eta_2}_{ij}(\phi) &= j_{\alpha,\beta}\left(S^{\varepsilon_1,\varepsilon_2,\eta_1,\eta_2}(S_{ij}\otimes\phi)\right), \\ T^{\varepsilon_1,\varepsilon_2,\eta_1,\eta_2}_{ij}(\phi) &= j_{\alpha,\beta}\left(T^{\varepsilon_1,\varepsilon_2,\eta_1,\eta_2}(T_{ij}\otimes\phi)\right) \end{aligned}$$

for $\phi \in \mathcal{H}^{m_1,m_2}(\mathbb{C}^p)\otimes\mathcal{H}^{n_1,n_2}(\mathbb{C}^q)$. Then the maps

(a) $\Pi(\rho(S_{ij})\{\Pi^{(m_1+1,\ m_2,\ n_1,\ n_2+1)}\rho(S_{ij})(j_{\alpha,\beta}(\phi))\})$
$\quad = A^{--}(m+1,\ n+1)A^{++}(m,n)S^{0-1-10}_{ij}(S^{1001}_{ij}(\phi))$,

(b) $\Pi(\rho(S_{ij})\{\Pi^{(m_1+1,\ m_2,\ n_1-1,\ n_2)}\rho(S_{ij})(j_{\alpha,\beta}(\phi))\})$,
$\quad = A^{-+}(m+1,\ n-1)A^{+-}(m,n)S^{0-101}_{ij}(S^{10-10}_{ij}(\phi))$,

(c) $\Pi(\rho(S_{ij})\{\Pi^{(m_1,\ m_2-1,\ n_1-1,\ n_2)}\rho(S_{ij})(j_{\alpha,\beta}(\phi))\})$
$\quad = A^{++}(m-1,\ n-1)A^{--}(m,n)S^{1001}_{ij}(S^{0-1-10}_{ij}(\phi))$,

(d) $\Pi(\rho(S_{ij})\{\Pi^{(m_1,\ m_2-1,\ n_1,\ n_2+1)}\rho(S_{ij})(j_{\alpha,\beta}(\phi))\})$
$\quad = A^{+-}(m-1,\ n+1)A^{-+}(m,n)S^{10-10}_{ij}(S^{0-101}_{ij}(\phi))$,

where

$$\Pi = \Pi^{(m_1+1,\ m_2-1,\ n_1-1,\ n_2+1)},$$

have the mapping property as in (6.5). It suffices to show that at least one of them has nonzero coefficients. A simple check using formulas (2.20) verifies that all four coefficients vanish only when $m = n = 0 = \alpha = \beta$, but when this happens, $V^{\alpha,\beta}_{(m_1\pm1,m_2\mp1,n_1\mp1,n_2\pm1)}$ vanishes, that is, the fiber at $(0,0)$ has only one $K$-type. Thus, we can always move in the south-westerly direction along the fiber. The proof for the north-easterly movement is similar.

## 6.3. Computations for $\mathrm{Sp}(p,q)$.

Choose $z_1,\ldots,z_{2p},w_1,\ldots,w_{2q},\overline{z_1},\ldots,\overline{z_{2p}},\overline{w_1},\ldots,\overline{w_{2q}}$ as a system of coordinates in $(\mathbb{H}^{p+q})_{\mathbb{C}}$. The complexified Lie algebra of $\mathrm{Sp}(p,q)$ may be described as the set of matrices in $M_{2p,2q}(\mathbb{C})$

$$\begin{bmatrix} A_{11} & A_{12} & A_{13} & A_{14} \\ A_{21} & -A^{\mathsf{t}}_{11} & A_{23} & A_{24} \\ A^{\mathsf{t}}_{24} & -A^{\mathsf{t}}_{14} & A_{33} & A_{34} \\ -A^{\mathsf{t}}_{23} & A^{\mathsf{t}}_{13} & A_{43} & -A^{\mathsf{t}}_{33} \end{bmatrix}$$

where $A_{11}$, $A_{12}$, $A_{21} \in M_p(\mathbb{C})$, $A_{33}$, $A_{34}$, $A_{43} \in M_q(\mathbb{C})$, $A_{13}$, $A_{14}$, $A_{23}$, $A_{24} \in M_{p,q}(\mathbb{C})$ with $A_{12} = A^{\mathsf{t}}_{12}$, $A_{21} = A^{\mathsf{t}}_{21}$, $A_{34} = A^{\mathsf{t}}_{34}$ and $A_{43} = A^{\mathsf{t}}_{43}$. Note that $\mathfrak{sp}(p)_{\mathbb{C}}$ and $\mathfrak{sp}(q)_{\mathbb{C}}$ are embedded block diagonally in $\mathfrak{sp}(p,q)_{\mathbb{C}}$.

It is not difficult to check that the following first-order differential operators span the action of $\mathfrak{sp}(p)_{\mathbb{C}}$ on the space of polynomial functions on $(\mathbb{H}^{p+q})_{\mathbb{C}}$ [$1 \le i,j \le p$ for (a) and $1 \le i \le j \le p$ for (b) and (c)]:

$$\text{(a)} \quad A^1_{ij} = z_i \frac{\partial}{\partial z_j} - \overline{z_j} \frac{\partial}{\partial \overline{z_i}} + \overline{z_{i+p}} \frac{\partial}{\partial \overline{z_{j+p}}} - z_{j+p} \frac{\partial}{\partial z_{i+p}},$$

$$(6.8) \quad \text{(b)} \quad B^1_{ij} = z_i \frac{\partial}{\partial z_{j+p}} - \overline{z_{j+p}} \frac{\partial}{\partial \overline{z_i}} + z_j \frac{\partial}{\partial z_{i+p}} - \overline{z_{i+p}} \frac{\partial}{\partial \overline{z_j}},$$

$$\text{(c)} \quad C^1_{ij} = z_{j+p} \frac{\partial}{\partial z_i} - \overline{z_i} \frac{\partial}{\partial \overline{z_{j+p}}} + z_{i+p} \frac{\partial}{\partial z_j} - \overline{z_j} \frac{\partial}{\partial \overline{z_{i+p}}}.$$

Similarly, the following describes the action of $\mathfrak{sp}(q)_{\mathbb{C}}$ [$1 \le i, j \le q$ for (a) and $1 \le i \le j \le q$ for (b) and (c)]:

$$\text{(a)} \quad A^2_{ij} = w_i \frac{\partial}{\partial w_j} - \overline{w_j} \frac{\partial}{\partial \overline{w_i}} + \overline{w_{i+q}} \frac{\partial}{\partial \overline{w_{j+q}}} - w_{j+q} \frac{\partial}{\partial w_{i+q}},$$

$$(6.9) \quad \text{(b)} \quad B^2_{ij} = w_i \frac{\partial}{\partial w_{j+q}} - \overline{w_{j+q}} \frac{\partial}{\partial \overline{w_i}} + w_j \frac{\partial}{\partial w_{i+q}} - \overline{w_{i+q}} \frac{\partial}{\partial \overline{w_j}},$$

$$\text{(c)} \quad C^2_{ij} = w_{j+q} \frac{\partial}{\partial w_i} - \overline{w_i} \frac{\partial}{\partial \overline{w_{j+q}}} + w_{i+q} \frac{\partial}{\partial w_j} - \overline{w_j} \frac{\partial}{\partial \overline{w_{i+q}}}.$$

Finally, the space $\mathfrak{p}$ complementary to $\mathfrak{sp}(p)_{\mathbb{C}} \oplus \mathfrak{sp}(q)_{\mathbb{C}}$ in $\mathfrak{sp}(p,q)_{\mathbb{C}}$ is spanned by the following ($1 \le i \le p$ and $1 \le j \le q$):

$$\text{(a)} \quad P_{ij} = z_i \frac{\partial}{\partial w_j} + w_{q+j} \frac{\partial}{\partial z_{p+i}} + \overline{w_j} \frac{\partial}{\partial \overline{z_i}} + \overline{z_{p+i}} \frac{\partial}{\partial \overline{w_{q+j}}},$$

$$(6.10) \quad \text{(b)} \quad Q_{ij} = z_i \frac{\partial}{\partial w_{q+j}} - w_j \frac{\partial}{\partial z_{p+i}} + \overline{w_{q+j}} \frac{\partial}{\partial \overline{z_i}} - \overline{z_{p+i}} \frac{\partial}{\partial \overline{w_j}},$$

$$\text{(c)} \quad U_{ij} = z_{p+i} \frac{\partial}{\partial w_j} - w_{q+j} \frac{\partial}{\partial z_i} + \overline{w_j} \frac{\partial}{\partial \overline{z_{p+i}}} - \overline{z_i} \frac{\partial}{\partial \overline{w_{q+j}}},$$

$$\text{(d)} \quad V_{ij} = z_{p+i} \frac{\partial}{\partial w_{q+j}} + w_j \frac{\partial}{\partial z_i} + \overline{w_{q+j}} \frac{\partial}{\partial \overline{z_{p+i}}} + \overline{z_i} \frac{\partial}{\partial \overline{w_j}}.$$

Let

$$E_{11} = \sum_{l=1}^{2p} z_l \frac{\partial}{\partial z_l}, \qquad E_{22} = \sum_{l=1}^{2p} \overline{z_l} \frac{\partial}{\partial \overline{z_l}},$$

$$E_{12} = \sum_{l=1}^{p} (z_{l+p} \frac{\partial}{\partial \overline{z_l}} - z_l \frac{\partial}{\partial \overline{z_{l+p}}}), \quad \text{and} \quad E_{21} = \sum_{l=1}^{p} \left( \overline{z_l} \frac{\partial}{\partial z_{l+p}} - \overline{z_{l+p}} \frac{\partial}{\partial z_l} \right).$$

It is not difficult to check that $\{\Delta_p, r_p^2, E_{11}, E_{12}, E_{21}, E_{22}\}$ commutes with $\mathfrak{sp}(p)$ on the space of polynomial functions

$$\mathcal{P}(\mathbb{H}_{\mathbb{C}}^p) \simeq \mathbb{C}[z_1, \ldots, z_{2p}, \overline{z_1}, \ldots, \overline{z_{2p}}].$$

They span a Lie algebra isomorphic to

$$\mathfrak{o}(4)_{\mathbb{C}}^* = \left\{ X \in \mathfrak{gl}(4) \ \middle| \ X\tilde{J} + \tilde{J}X = 0 \text{ and } XJ + JX = 0 \right\}_{\mathbb{C}},$$

where

$$\tilde{J} = \begin{bmatrix} 0 & I_2 \\ I_2 & 0 \end{bmatrix} \quad \text{and} \quad J = \begin{bmatrix} I_2 & 0 \\ 0 & -I_2 \end{bmatrix}.$$

Here $I_2$ is the $2 \times 2$ identity matrix. Observe that $\{E_{11}, E_{12}, E_{21}, E_{22}\}$ spans the Lie algebra $\mathfrak{gl}(2)$ where $E_{ij}$ corresponds to the standard matrix unit at $(i, j)$ in $\mathfrak{gl}(2)$. In particular, one sees that the Lie algebra spanned by $\{E_{11} - E_{22}, \ E_{12}, \ E_{21}\}$ corresponds to the infinitesimal action of $\mathbb{H}_1 \simeq \mathrm{Sp}(1)$.

With the notation as in §5.2, we have the following result, which implies Proposition 5.1 immediately.

**Lemma 6.1.** *The decomposition (for $p > 1$)*

$$\mathcal{H}(\mathbb{H}^p)\Big|_{\mathrm{Sp}(p) \times \mathbb{H}_1} = \sum V_p^{(\xi_1, \xi_2)} \otimes V_1^{\xi_1 - \xi_2}$$

*runs through all representations $V_p^{(\xi_1, \xi_2)}$ of $\mathrm{Sp}(p)$. The $\mathrm{Sp}(p)$ module $V_p^{(\xi_1, \xi_2)} \otimes V_1^{\xi_1 - \xi_2}$ is generated by the joint highest weight vector*

$$(6.11) \qquad h_{p,(\xi_1, \xi_2)} = z_1^{\xi_1 - \xi_2} \begin{vmatrix} z_1 & \overline{z_{1+p}} \\ z_2 & \overline{z_{2+p}} \end{vmatrix}^{\xi_2}.$$

*For $p = 1$ the isotypic decomposition of $\mathcal{H}(\mathbb{H})$ is*

$$\mathcal{H}(\mathbb{H})\Big|_{\mathrm{Sp}(1) \times \mathrm{Sp}(1)} = \sum_{m=0}^{\infty} V_1^m \otimes V_1^m,$$

*where $V_1^m \otimes V_1^m$ is generated by $h_{1,m} = z_1^m$.*

*Remark.* If we apply the lowering operator $E_{21}$ in $\mathfrak{sp}(1)$ to the joint highest weight vector (6.11), we find that a spanning set of $\mathrm{Sp}(p)$ highest weight vectors in $\mathcal{H}(\mathbb{H}^p)$ are given by

$$(6.12) \qquad \begin{cases} h_{p,(\xi_1,\xi_2),j} = \overline{z_{p+1}}^j h_{p,(\xi_1-j,\xi_2)}, & j = 0, 1, ..., \xi_1 - \xi_2, \\ h_{1,m,j} = \overline{z_2}^j h_{1,m-j}, & j = 0, 1, ..., m. \end{cases}$$

*Proof.* Choose the positive roots of $\mathfrak{sp}(p)_{\mathbb{C}}$ to be those described by (6.8)(a) (for $1 \leq i \leq j \leq p$) and (6.8)(b). It is easy to check that $h_{n,(\xi_1,\xi_2)}$ are $\mathfrak{sp}(p)_{\mathbb{C}} \times \mathfrak{sp}(1)_{\mathbb{C}}$ highest weight vectors of weight $(\xi_1, \xi_2, 0, 0, \ldots, 0) \otimes (\xi_2 - \xi_2)$. For the remaining part of the first assertion, we refer to [Ho2]. $\square$

*Proof of Lemma 5.4.* We proceed to compute the action of the $\mathfrak{p}$-part of $\mathrm{Sp}(p, q)$ on $S_j^a(X^0)$. From Lemma 6.1, the $\mathrm{Sp}(p) \times \mathrm{Sp}(q)$ highest weight vectors in $\mathcal{H}^m(\mathbb{R}^{4p}) \otimes \mathcal{H}^n(\mathbb{R}^{4q})$ are

$$h_{p,(a_1,a_2),j} h_{q,(b_1,b_2),k} = \overline{z_{p+1}}^j h_{p,(a_1-j,a_2)} \overline{w_{q+1}}^k h_{q,(b_1-k,b_2)}.$$

For fixed $(a_1, a_2)$ and $(b_1, b_2)$, these vectors span an $\mathfrak{sp}(1)$ module isomorphic to $V_1^{a_1-a_2} \otimes V_1^{b_1-b_2}$. The Clebsch-Gordan formula tells us that the $\mathfrak{sp}(1)$ highest weight vectors in this module are the functions

$$h_{p,(a_1-j,a_2)} h_{q,(b_1-j,b_2)} \begin{vmatrix} z_1 & \overline{z_{p+1}} \\ w_1 & \overline{w_{q+1}} \end{vmatrix}^j.$$

To adjust the homogeneity of these functions, we should multiply by an appropriate power of $r_p^2$. Thus we wish to study the functions

$$\gamma(a_1, a_2, b_1, b_2, c, d) = z_1^{a_1} \begin{vmatrix} z_1 & \overline{z_{p+1}} \\ z_2 & \overline{z_{p+2}} \end{vmatrix}^{a_2} w_1^{b_1} \begin{vmatrix} w_1 & \overline{w_{q+1}} \\ w_2 & \overline{w_{q+2}} \end{vmatrix}^{b_2} \begin{vmatrix} z_1 & \overline{z_{q+1}} \\ w_1 & \overline{w_{q+1}} \end{vmatrix}^c (r_p^2)^d.$$

Here $a_j, b_j$, and $c$ are positive integers and $d$ is a complex number. Note that our use here of $a_j$, $b_j$ is slightly different from that of a few lines above. A check shows that $\gamma(a_1, a_2, b_1, b_2, c, d)$ is an $\mathrm{Sp}(p) \times \mathrm{Sp}(q)$ highest weight vector of $\mathrm{Sp}(p)$ weight $(a_1 + a_2 + c, \; a_2, 0, ..., 0)$ and $\mathrm{Sp}(q)$ weight $(b_1 + b_2 + c, \; b_2, 0, ..., 0)$. Also it is an $\mathfrak{sp}(1)$ highest weight vector of weight $a_1 + b_1$.

Recall the notation in §5.3. The possible components in $\mathfrak{p} \otimes V_p^{(\xi_1, \xi_2)} \otimes V_q^{(\eta_1, \eta_2)}$ are the sixteen components obtained by adding $(\pm 1, 0, \pm 1, 0)$, $(\pm 1, 0, 0, \pm 1)$, $(0, \pm 1, \pm 1, 0)$ and $(0, \pm 1, 0, \pm 1)$ to $(\xi_1, \xi_2, \eta_1, \eta_2)$. Our intention will be to show that we can move from the fiber over $(m, n)$ to a neighboring fiber $(m \pm 1, \; n \pm 1)$ and then back, and doing this a sufficient number of times to get to any other points in the original fiber.

Because the moves between, say $(\xi_1, \xi_2, \eta_1, \eta_2)$ to $(\xi_1 + 1, \xi_2, \; \eta_1 + 1, \; \eta_2)$, and in the reverse direction are essentially adjoint to one another (exactly so on the unitary axis), we can do one move if and only if we can do the reverse move. Thus it suffices to check if we can do eight of the sixteen moves.

To show that we can do the (1,0,1,0) move, we apply $Q_{11}$ [see (6.10)] to one of our functions $\gamma$ :

(6.13)      $$Q_{11}\gamma(a_1, a_2, b_1, b_2, c, d) = d\gamma(a_1, a_2, b_1, b_2, \; c+1, \; d-1).$$

To show that we can do the (1,0,0,1) move, consider

(6.14)      $$\begin{aligned} &(A_{21}^2 Q_{11} - (b_1 + c + 1)Q_{12})\gamma(a_1, a_2, b_1, b_2, c, d) \\ &= -b_1 d\gamma(a_1 + 1, \; a_2, \; b_1 - 1, \; b_2 + 1, \; c, \; d-1). \end{aligned}$$

Similarly for the (0,1,1,0) move, we appeal to the formula

(6.15)      $$\begin{aligned} &(A_{21}^1 Q_{11} - (a_1 + c + 1)Q_{21})\gamma(a_1, a_2, b_1, b_2, c, d) \\ &= -a_1 d\gamma(a_1 - 1, \; a_2 + 1, \; b_1 + 1, \; b_2, \; c, \; d-1). \end{aligned}$$

Finally, for the (0,1,0,1) move, we can verify that the operator

$$\begin{aligned} &A_{21}^1 A_{21}^2 Q_{11} - (a_1 + c + 1)A_{21}^2 Q_{21} \\ &\quad - (b_1 + c + 1)A_{21}^1 Q_{12} + (a_1 + c + 1)(b_1 + c + 1)Q_{22} \end{aligned}$$

applied to $\gamma(a_1, a_2, b_1, b_2, c, d)$ yields

(6.16)      $$c(a_1 + b_1 + c + 1)d\gamma(a_1, \; a_2 + 1, \; b_1, \; b_2 + 1, \; c-1, \; d-1).$$

Thus in passing from the fiber $U_{a,j}^{(m,n)}$ to the fiber $U_{a,j}^{(m+1,n+1)}$ (see §5.3, before Lemma 5.3) by means of the action of the $Q_{ij}$, we can add all four of the conceivable vectors $(1,0,1,0)$, $(0,1,1,0)$, $(1,0,0,1)$, and $(0,1,0,1)$ to a given highest weight of $U_{a,j}^{(m,n)}$. By taking adjoints, we can add the negatives of these vectors when passing from $U_{a,j}^{(m+1,n+1)}$ to $U_{a,j}^{(m,n)}$. Thus passing from $U_{a,j}^{(m,n)}$ to $U_{a,j}^{(m+1,n+1)}$ and back again allows us to move from a given $K$-type in $U_{a,j}^{(m,n)}$ to any of its eight nearest neighbors, as illustrated in Diagram 5.17.

If the vanishing of either of the transition coefficients $A_{4p,4q,a}^{++}(m,n)$ or $A_{4p,4q,a}^{--}(m+1, n+1)$ prevents us from moving to $U_{a,j}^{(m+1,n+1)}$ or back to $U_{a,j}^{(m,n)}$, we can equally well move from $U_{a,j}^{(m,n)}$ down to $U_{a,j}^{(m-1,n-1)}$ and back up to achieve the same end. Study of the diagrams of the $K$-type regions (Diagrams 5.12, 5.14, 5.16) and the positions of the barriers (Diagram 2.34) shows that the only time neither of these strategies works is when $(m,n)$ is on the border of the $K$-type region; but in this case, we can see from Diagram 5.10 that the fiber $U_{a,j}^{(m,n)}$ contains only one $K$-type, so there is nothing to show. Although we do not need them, similar calculations may be done for the transitions from $U_{a,j}^{(m,n)}$ to $U_{a,j}^{(m+1,n-1)}$ and back. This concludes the proof of Lemma 5.4. $\quad\square$

*Concluding remark.* Although in this paper we have been exclusively concerned with questions of infinitesimal structure, our results have some implications for harmonic analysis on hyperbolic spaces (also describable as isotropic (reductive) symmetric spaces [Sc]). To describe these, consider the complex projective hyperboloid $Y = \mathrm{U}(p,q)/(\mathrm{U}(p-1,q) \times \mathrm{U}(1))$. Similar remarks apply to the quaternionic hyperboloid. As noted in [Sc], the real hyperboloid

$$X_1 = \{w \in \mathbb{R}^{2p+2q} \mid (w,w)_{p,q} = 1\}$$

fibers over $Y$ with fiber $\mathrm{U}(1)$. In fact, this is a principle $\mathrm{U}(1)$-bundle. It is well known (see the references on the first page, especially [Sc]) that all the representations occurring in the harmonic analysis of the $\mathrm{O}(2p,2q)$ action on $X_1$ are among those considered in §§2 and 3. Our results imply that the restriction of any one of these representations to $\mathrm{U}(p,q)$ decomposes discretely; more precisely, for any of these representations of $\mathrm{O}(2p,2q)$, the eigenspaces of the fiber $\mathrm{U}(1)$ of $X_1$ over $Y$ will be irreducible modules for $\mathrm{U}(p,q)$. It follows that the $\mathrm{U}(p,q)$ harmonic analysis of any line bundle constructed from the fibration $X_1 \longrightarrow Y$ reduces immediately to the $\mathrm{O}(2p,2q)$ harmonic analysis on $X_1$. In particular, harmonic analysis on $Y$ reduces to harmonic analysis on $X_1$.

## 7. Afterword of caution

In this article we have tried to exhibit some of the elegance of representation theory of semisimple groups. We hope some readers will be inspired to go further, in the texts [Kn], [Wa], or the expository accounts [Ho4], [Vo4]. Here at the end, we should warn the reader otherwise unfamiliar with representation theory that the examples we have studied are still only a small, and in some ways misleadingly comprehensible, sample of the full spectrum of representations of semisimple groups.

One way to locate our representations in relation to the whole is to recall their description as degenerate principal series [see, e.g., formulas (2.9)–(2.11)]. They are representations induced from certain finite-dimensional representations of the group $P_1$, the stabilizer of an isotropic line. This group is a maximal parabolic subgroup, and it is its maximality, and the corresponding smallness of the coset space $G/P_1$, that leads to the use of the adjective " degenerate". If one considers representations induced from finite-dimensional irreducible representations of a minimal, rather than maximal, parabolic subgroup, then one has the (ordinary or nondegenerate) principal series. [For $\mathrm{O}(p,q)$ a general parabolic subgroup is describable as the stabilizer of a nested sequence $X_1 \subset X_2 \subset X_3 \subset \cdots \subset X_k$ of isotropic subspaces. If the sequence contains a subspace of every dimension from one to $\min(p,q)$, then the parabolic stabilizing it is minimal. Similar remarks apply to $\mathrm{U}(p,q)$ and $\mathrm{Sp}(p,q)$.] Principal series have the following properties:

(1) They come in continuous families, parametrized by several (in place of the single parameter of our examples) complex variables.
(2) For generic values of the parameters, they are irreducible.
(3) In general, they have finite composition series.
(4) Any (reasonable) irreducible representation is equivalent (via an appropriate notion of equivalence) to a constituent of some principal series.

Properties (1)–(3) are qualitatively similar to what we saw in this paper, and property (4) tells us that if we understand the principal series, then we understand an essential part, in particular the building blocks, of semisimple representation theory.

It might seem that properties (1)–(4) virtually finish the subject, but no one has so far managed to achieve control of general principal series comparable to the picture given for the examples in this paper. The point is that the combinatorial structure of the composition series of a general principal series, and of its constituents, is rather spectacularly more complex than what we encountered in our examples. (See [CC] for a recent discussion.) Some of the most important advances in the subject (Langlands classification [Kn, Wa, Vo2]; Kazhdan-Lusztig conjectures [BB, BK, KL, Vo3]) have involved obtaining significant information about the composition series of principal series.

Let us remark that when $q = 1$, our maximal parabolic is the minimal and only proper parabolic subgroup, and in this case our representations are ordinary, not degenerate, principal series. When $p$ is also small, our representations in fact include all, or essentially all, principal series, so that we find all the irreducible representations, and since we know which are unitary, we classify the unitary dual. Specifically we obtain all principal series for $\mathrm{U}(1,1)$ and for $\mathrm{Sp}(1,1)$ ($\simeq \mathrm{Spin}(4,1)$), and we obtain all principal series up to twisting by a quasi-character of the whole group for $\mathrm{O}(2,1)$ and $\mathrm{U}(2,1)$.

Since, at least when $p$ and $q$ are greater than one, the representations we have considered are quite special, one may ask if any interest attaches to them, aside from the inherent attractiveness of the picture we have tried to present. In fact, these representations arise in studying harmonic analysis on hyperboloids. The full import of the term "harmonic analysis" has never been completely decided, but many would agree that it certainly includes the decomposition of $L^2$-spaces of homogeneous spaces. Among the homogeneous spaces of semisimple Lie groups, the class of "semisimple symmetric spaces" or "affine symmetric spaces" [F-J] recommends itself as particularly interesting. The hyperboloids

$$X_t = \{ w \in \mathbb{R}^{p+q} \mid (w, w)_{p,q} = t \}$$

are very closely related to (precisely, are principal fibre bundles with compact fibre over) some of the simplest examples of semisimple symmetric spaces [Sc], and the representations we have studied enable one to decompose the $L^2$-spaces of hyperboloids. The connection with harmonic analysis on semisimple symmetric spaces was the motivation for many of the papers cited in the introduction.

DEPARTMENT OF MATHEMATICS, YALE UNIVERSITY, NEW HAVEN, CONNECTICUT 06520
*E-mail address*: `howe@lom1.math.yale.edu`

DEPARTMENT OF MATHEMATICS, NATIONAL UNIVERSITY OF SINGAPORE, KENT RIDGE, SINGA-PORE 0511, SINGAPORE
*E-mail address*: `mattanec@nusvm.bitnet`